# Structures de contact en dimension trois et bifurcations des feuilletages de surfaces

Emmanuel GIROUX*

Août 1999

Le but principal de cet article est de classifier les structures de contact sur quelques variétés de dimension trois, à savoir les espaces lenticulaires, la plupart des fibrés en tores au-dessus du cercle, le tore plein et le tore épais. Cette classification complète divers travaux antérieurs [Et, Gi3, Gi4, Ka, Ma, Th] et résulte de la confrontation entre deux techniques : la chirurgie qui fournit de nombreuses structures de contact et la tomographie qui permet d'analyser une structure de contact donnée *a priori* et d'en former une image combinatoire. Les méthodes de chirurgie reposent sur un théorème de Y. ELIASHBERG [El2, El5] — revisité par R. GOMPF dans [Go] — et produisent des structures de contact holomorphiquement remplissables sur des variétés closes. L'étude tomographique — développée dans les parties 2 et 3 de ce texte — s'appuie sur les notions introduites dans [Gi1] et dégage un petit nombre de modèles possibles pour les structures de contact tendues sur chacune des variétés citées plus haut.

Le texte qui suit s'organise ainsi. La partie 1 donne les énoncés précis des résultats et présente les outils qui servent à leur démonstration. La partie 2 étudie systématiquement les propriétés des familles à un paramètre de feuilletages sur une surface $F$ qui proviennent des structures de contact sur le produit de $F$ par l'intervalle ; cette étude conduit à de nouvelles démonstrations des théorèmes de D. BENNEQUIN [Be] et de Y. ELIASHBERG [El4] concernant les structures de contact sur la sphère $\mathbf{S}^3$. La partie 3 prolonge la partie 2 dans le cas où $F$ est un tore et fournit des formes normales pour les structures de contact tendues sur le tore épais. La partie 4 complète la démonstration des théorèmes de classification.

Je tiens beaucoup à remercier ici Bruno SÉVENNEC ; sans sa culture mathématique et sa disponibilité infinies, ce travail n'aurait probablement pas vu le jour. Le texte qui suit a été en partie écrit lors d'un séjour à l'American Institute of Mathematics ; je remercie vivement cet institut pour son accueil et le confort qu'il m'a offert. Charles THOMAS et Yasha ELIASHBERG m'avaient auparavant donné l'occasion de présenter des ébauches de ce travail à Cambridge en juillet 94, à Park City en juillet 97 et à Stanford en novembre 98 ; je les remercie, ainsi que John ETNYRE pour l'intérêt qu'il a porté à mes exposés. D'autre part, certains des résultats démontrés ici ont été très récemment annoncés par Ko HONDA.

## 1 Présentation des résultats

Pour énoncer les résultats, quelques précisions sont nécessaires. Tout d'abord, on se donne une variété orientée $V$ de dimension trois et les structures de contact qu'on

---

considère sur $V$ sont orientables et directes. Autrement dit, chacune d'elles est le noyau $\xi$ d'une 1-forme $\alpha$ dont le produit avec $d\alpha$ est partout non nul et positif pour l'orientation choisie. Ensuite, on classifie les structures de contact $\xi$ qui sont *tendues*, c'est-à-dire telles qu'aucun disque plongé dans $\operatorname{Int} V = V \setminus \partial V$ ne soit tangent à $\xi$ en tous les points de son bord. Les autres structures de contact, dites *vrillées*, ont été classifiées par Y. ELIASHBERG dans [El1].

Les exemples principaux de structures de contact tendues sont issus de la géométrie complexe : le bord orienté $V$ d'une surface complexe compacte $W$ est strictement pseudo-convexe si et seulement si le champ $\xi$ des droites complexes tangentes à $V$ est une structure de contact directe. Les structures de contact ainsi obtenues sont dites *holomorphiquement remplissables* et sont tendues d'après un théorème de M. GROMOV [Gr].

D'autres exemples viennent de l'observation suivante : si une structure de contact sur $V$ induit une structure tendue sur un revêtement de $V$, elle est tendue. Le fait que l'assertion analogue pour les structures vrillées soit fausse conduit à une nouvelle dichotomie : on dit qu'une structure de contact tendue sur $V$ est *virtuellement vrillée* si elle induit une structure vrillée sur un revêtement fini de $V$ ; à l'opposé, une structure de contact est *universellement tendue* si elle induit une structure tendue sur le revêtement universel. Lorsque le groupe fondamental de $V$ est résiduellement fini — ce qui est le cas pour les variétés étudiées ici —, toute structure de contact tendue sur $V$ est soit virtuellement vrillée, soit universellement tendue. Il est à noter qu'on trouve des structures de contact holomorphiquement remplissables dans les deux catégories.

## A  Les espaces lenticulaires

Un des nombres qui mesurent la complexité d'une variété $V$ de dimension trois est son genre, défini comme le plus petit genre d'une surface connexe plongée dans $V$ et séparant $V$ en deux morceaux dont chacun se rétracte sur un bouquet de cercles. Les variétés de genre 1 sont les espaces lenticulaires. Si on exclut $\mathbf{S}^1 \times \mathbf{S}^2$, ils forment une famille $\mathbf{L}_{p,q}$ indexée par des entiers premiers entre eux tels que $0 < q < p$. L'espace $\mathbf{L}_{p,q}$ est le quotient orienté du produit $\mathbf{T}^2 \times [0,1]$ par la relation d'équivalence dont les classes non ponctuelles sont les cercles de direction $(0,1)$ sur $\mathbf{T}^2 \times \{0\}$ et ceux de direction $(p,q)$ sur $\mathbf{T}^2 \times \{1\}$.

**Théorème 1.1.** *Soit $p$ et $q$ des entiers premiers entre eux tels que $0 < q < p$. Les structures de contact tendues et orientées sur $\mathbf{L}_{p,q}$ sont toutes holomorphiquement remplissables et forment un nombre fini de classes d'isotopie égal à $\prod_{i=0}^{n}(a_i - 1)$ où*

$$\frac{p}{q} = a_0 - \cfrac{1}{a_1 - \cfrac{1}{\cdots - \cfrac{1}{a_n}}} \quad \text{avec } a_i \geq 2 \text{ pour } 0 \leq i \leq n.$$

*En outre, ces classes d'isotopie occupent des classes d'homotopie distinctes dans l'espace des champs de plans orientés et toutes, sauf deux si $q < p - 1$ et une si $q = p - 1$, sont constituées de structures de contact virtuellement vrillées.*

Si on s'intéresse aux structures de contact tendues orientables mais pas orientées — ce qui est le plus souvent le cas dans la suite —, le nombre de leurs classes d'isotopie sur $\mathbf{L}(p,q)$ est $\lceil \frac{1}{2} \prod_{i=0}^{n}(a_i - 1) \rceil$, où $\lceil x \rceil$ désigne le plus petit entier supérieur ou égal à $x$.

# B  Les fibrés en tores au-dessus du cercle

Les variétés fibrées en tores au-dessus du cercle sont, comme les espaces lenticulaires, des quotients simples du tore épais $\mathbf{T}^2 \times [0,1]$. Cela dit, leur genre vaut 2 et, surtout, elles contiennent un tore incompressible. Ce tore vient compliquer la classification des structures de contact en faisant apparaître un phénomène de torsion.

**Définition 1.2.** *Soit $\xi$ une structure de contact sur une variété $V$ de dimension trois. On appelle $2\pi$-torsion de $\xi$ le supremum des entiers $n \geq 1$ pour lesquels on peut plonger dans $(\operatorname{Int} V, \xi)$ le produit $\mathbf{T}^2 \times [0,1]$ muni de la structure de contact d'équation*

$$\cos(2n\pi z)\,dx - \sin(2n\pi z)\,dy = 0, \qquad (x,y,z) \in \mathbf{T}^2 \times [0,1].$$

*Si de tels entiers n'existent pas, la torsion est décrétée nulle.*

*On définit de même la $\pi$-torsion en prenant sur $\mathbf{T}^2 \times [0,1]$ la structure de contact d'équation $\cos(n\pi z)\,dx - \sin(n\pi z)\,dy = 0$.*

Le calcul de cet invariant est délicat. Il n'est par exemple pas établi que toute structure de contact tendue sur une variété compacte ait une torsion finie (cf. conjecture 1.4). Sur les variétés envisagées ici, le calcul complet de la torsion ressort de l'étude tomographique. D'autre part, il découle de [El1] que toute structure de contact vrillée a une torsion infinie.

Pour toute matrice $A \in \operatorname{SL}_2(\mathbf{Z})$, on note $\mathbf{T}_A^3$ le quotient orienté du produit $\mathbf{T}^2 \times \mathbf{R}$ par la transformation $(x,t) \mapsto (Ax, t+1)$. Le type topologique de $\mathbf{T}_A^3$ est déterminé par la classe de conjugaison de $A$ et dépend fortement de la nature dynamique de $A$. Celle-ci transparaît aussi dans la classification des structures de contact, ce qui crée une multitude de situations différentes. La plus intéressante est sans doute celle où $A$ est hyperbolique ($|\operatorname{tr}(A)| > 2$). Dans ce cas, $A$ laisse deux droites invariantes ; la pente de sa droite instable (resp. stable) dans les diverses bases entières du plan décrit l'orbite d'un nombre quadratique sous l'action de $\operatorname{PSL}_2(\mathbf{Z})$ et la queue de son développement en fraction continue est donc périodique et indépendante de la base.

**Théorème 1.3.** *Soit $A$ une matrice dans $\operatorname{SL}_2(\mathbf{Z})$. Lorsque $A$ est hyperbolique, on note $z$ la pente de sa direction instable et $a_m, \ldots, a_{m+n}$ les entiers formant une période dans le développement de $z$ en fraction continue du type*

$$z = a_0 - \cfrac{1}{\cdots - \cfrac{1}{a_m - \cfrac{1}{a_{m+1} - \cdots}}} \quad \text{avec } a_i \geq 2 \text{ pour tout } i \geq 1.$$

**a)** *Les structures de contact universellement tendues sur $\mathbf{T}_A^3$ forment un nombre infini de classes de conjugaison où la $2\pi$-torsion prend toutes les valeurs entières positives ou nulles. De plus, si $\operatorname{tr}(A) \neq -2$, chaque classe est uniquement déterminée par sa $2\pi$-torsion.*

**b)** *Les structures de contact virtuellement vrillées et orientées sur $\mathbf{T}_A^3$ ont toutes une torsion nulle et forment un nombre fini $N(A)$ de classes d'isotopie. En outre, $N(A)$ vaut :*

1) $\prod_{i=0}^{n}(a_{m+i} - 1)^k$ *si $\operatorname{tr}(A) < -2$ et si $A = P^k$ où $P \in \operatorname{SL}_2(\mathbf{Z})$ est primitive ;*
2) $\prod_{i=0}^{n}(a_{m+i} - 1)^k - 2$ *si $\operatorname{tr}(A) > 2$ et si $A = P^k$ où $P \in \operatorname{SL}_2(\mathbf{Z})$ est primitive ;*
3) $|k| - 1$ *si $A = \begin{pmatrix} 1 & 0 \\ k & 1 \end{pmatrix}$ avec $k < 0$ ;*

4) 0 si $A = \begin{pmatrix} 1 & 0 \\ k & 1 \end{pmatrix}$ avec $k > 0$ impair ou si $A$ est d'ordre fini mais non conjuguée à $\begin{pmatrix} -k & -1 \\ 1 & 0 \end{pmatrix}$, $k \in \{0, 1\}$ ;

5) au plus 2 dans tout autre cas où $\operatorname{tr}(A) \neq -2$.

La partie a) de ce théorème précise le théorème 6 de [Gi4] ; en fait, l'étude tomographique dit que toute structure de contact universellement tendue sur $\mathbf{T}_A^3$ est conjuguée à l'une des structures étudiées dans [Gi4]. La partie b) laisse en suspens le calcul de $N(A)$ dans quelques cas. Un autre problème non résolu est la détermination des structures de contact remplissables sur $\mathbf{T}_A^3$. Il est probable que toutes les structures de torsion nulle sont holomorphiquement remplissables mais que, comme dans le cas du tore $\mathbf{T}^3$, les autres ne sont que symplectiquement remplissables [El4, Gi3].

Pour conclure cette partie consacrée aux variétés closes, on propose une conjecture suggérée par les théorèmes 1.1, 1.3 et les résultats de [Co1, Co2, Gi4, Gi5] :

**Conjecture 1.4.** Soit $V$ une variété close orientée de dimension trois. Les structures de contact virtuellement vrillées sur $V$ forment toujours un nombre fini de classes d'isotopie et les structures de contact universellement tendues forment un nombre infini de classes d'isotopie — et de conjugaison — si et seulement si $V$ est *toroïdale*, c'est-à-dire contient un tore incompressible plongé.

## C   Le tore épais

L'ingrédient clé dans la démonstration des théorèmes 1.1 et 1.3 est l'analyse des structures de contact tendues sur le tore épais $\mathbf{T}^2 \times [0, 1]$. Sur cette variété, on classifie les structures modulo l'action des $\partial$-isotopies, c'est-à-dire des isotopies relatives au bord. Plutôt que de simplement compter les classes, on cherche des invariants qui les déterminent.

La trace d'une structure de contact sur le bord du tore épais est un premier invariant de sa classe de $\partial$-isotopie. Cette trace, lorsque la structure est transversale au bord, est un champ de droites orientable qui engendre un feuilletage non singulier. Pour $i = 0, 1$, on fixe désormais un feuilletage non singulier $\sigma_i$ de $\mathbf{T}^2 \times \{i\}$ et on s'intéresse à l'espace $\mathcal{SCT}(\mathbf{T}^2 \times [0, 1]; \sigma_0, \sigma_1)$ des structures de contact tendues sur $\mathbf{T}^2 \times [0, 1]$ qui impriment $\sigma_0 \sqcup \sigma_1$ sur le bord. D'après un théorème classique de J. GRAY, les composantes connexes de cet espace pour la topologie $\mathcal{C}^\infty$ sont exactement les classes de $\partial$-isotopie cherchées. On note aussi $\mathcal{SCT}_u(\dots)$ et $\mathcal{SCT}_v(\dots)$ les sous-espaces de $\mathcal{SCT}(\dots)$ constitués des structures de contact respectivement universellement tendues et virtuellement vrillées.

Dans $\mathcal{SCT}(\mathbf{T}^2 \times [0, 1]; \sigma_0, \sigma_1)$, il y a de nouveau un invariant évident de $\partial$-isotopie, pour peu que l'on choisisse une orientation de $\sigma_0$. En effet, ce choix détermine une unique orientation compatible pour chaque structure de contact $\xi$ dans $\mathcal{SCT}(\mathbf{T}^2 \times [0, 1]; \sigma_0, \sigma_1)$ (cf. 2.A). Dès lors, $\xi$ possède une classe d'Euler relative

$$\chi_\partial(\xi) \in H_1\big(\mathbf{T}^2 \times [0, 1]; \mathbf{Z}\big) \cong \mathbf{Z}^2 \, .$$

Cette classe est la classe d'homologie de l'entrelacs orienté formé par les zéros d'une section générique de $\xi$ qui, sur le bord, est partout non nulle et engendre $\sigma_0 \sqcup \sigma_1$.

On suppose dorénavant que chaque feuilletage $\sigma_i$ est une suspension (c'est-à-dire qu'une certaine courbe fermée simple coupe transversalement toutes les orbites) et est de l'un des types suivants :

0) $\sigma_i$ est topologiquement linéarisable, c'est-à-dire que ses orbites sont toutes fermées ou toutes denses.

1) $\sigma_i$ a deux orbites fermées qui sont non dégénérées (le germe d'holonomie n'est pas tangent à l'identité) et entre lesquelles spiralent les autres orbites.

Un tel feuilletage $\sigma_i$ a une *direction asymptotique* $D(\sigma_i) = D_i \subset H_1(\mathbf{T}^2; \mathbf{R})$ : les cycles réels $[L_t]/t$, où $L_t$ est un segment d'orbite de longueur $t$ complété en un lacet par un segment géodésique minimal, ont — au signe près — une limite commune non nulle dans $H_1(\mathbf{T}^2; \mathbf{R})$ [Sc] et la droite $D_i$ qu'engendre cette limite (la demi-droite si $\sigma_i$ est orienté) est indépendante de la métrique. On note alors :
- $C = C(\sigma_0, \sigma_1) \subset H_1(\mathbf{T}^2; \mathbf{R}) \setminus \{0\}$ l'un des deux cônes convexes bordés à gauche par $D_0$, à droite par $D_1$ et ouvert (resp. fermé) du côté de $D_i$ si $\sigma_i$ est de type 0 (resp. 1) ;
- $E = E(\sigma_0, \sigma_1)$ l'enveloppe convexe des points entiers de $C$ ;
- $B = B(\sigma_0, \sigma_1)$ l'ensemble des points entiers primitifs de $\partial E$, ordonné par la relation « $w \preceq w'$ si $w \wedge w' \leq 0$ ».

L'ensemble $B$ est un « intervalle de points entiers » sur $\partial E$ (on exclut simplement les points non primitifs qu'il peut y avoir sur $D_i$ si $D_i \subset C$). On note dans la suite $\partial B$ l'ensemble formé des extrémités éventuelles de $B$ et $\operatorname{Int} B$ l'ensemble $B \setminus \partial B$. En outre, pour toute partie finie $Q$ de $B$, on pose

$$\gamma(Q) = \begin{cases} 0 & \text{si } Q = \varnothing, \\ \sum_{j=1}^{k} (-1)^{j-1} w_j & \text{si } Q = \{w_1, \ldots, w_k\} \text{ avec } w_1 \prec \cdots \prec w_k. \end{cases}$$

**Théorème 1.5.** *Soit $\sigma_0$ (resp. $\sigma_1$) une suspension sur $\mathbf{T}^2$ de type 0 ou 1. On oriente $\sigma_0$ de telle sorte que la demi-droite de ses cycles asymptotiques soit le bord gauche du cône $C$.*

*La classe d'Euler relative et la $\pi$-torsion définissent une application*

$$\chi_\partial \times \tau_\pi \colon \mathcal{SCT}(\mathbf{T}^2 \times [0,1]; \sigma_0, \sigma_1) \longrightarrow H_1(\mathbf{T}^2 \times [0,1]; \mathbf{Z}) \times \mathbf{N}$$

*dont les fibres sont toutes connexes, sauf dans un cas décrit plus bas. De plus, les sous-espaces $\mathcal{SCT}_u$ et $\mathcal{SCT}_v$ ont pour images respectives $X_u \times \mathbf{N}$ et $X_v \times \{0\}$ où*

$$X_u = \{2\gamma(Q) \mid Q \cap \operatorname{Int} B = \varnothing,\ Q \subset B\},$$
$$X_v = \{2\gamma(Q) \mid Q \cap \operatorname{Int} B \neq \varnothing,\ Q \subset B\}.$$

*L'exception est la suivante : lorsque $\sigma_0$ et $\sigma_1$ sont tous deux de type 1 et ont la même direction asymptotique, deux fibres ont une infinité de composantes connexes supplémentaires qui, sous l'action des difféomorphismes relatifs au bord, ne forment cependant qu'une orbite. Ces fibres exceptionnelles sont celles de $(0,0)$ et de $(2w, 0)$, où $w$ est le point entier primitif situé sur le bord gauche de $C$.*

En particulier, lorsque l'une des directions asymptotiques $D_0$, $D_1$ est une droite irrationnelle, $\mathcal{SCT}(\mathbf{T}^2 \times [0,1]; \sigma_0, \sigma_1)$ compte une infinité de composantes connexes qui, dans l'espace des champs de plans, occupent des classes d'homotopie distinctes. Sur les variétés closes, il est au contraire très probable que seul un nombre fini de classes d'homotopie dans l'espace des champs de plans contiennent des structures de contact tendues [KM].

## D  Le tore plein

La classification qu'on donne sur le tore plein $\mathbf{S}^1 \times \mathbf{D}^2$ est semblable à celle obtenue sur le tore épais mais plus simple car, en l'absence de tout tore incompressible, la torsion disparaît.

On désigne par $\mathcal{SCT}(\mathbf{S}^1 \times \mathbf{D}^2; \sigma)$ l'espace des structures de contact tendues sur $\mathbf{S}^1 \times \mathbf{D}^2$ qui impriment sur le bord $\mathbf{S}^1 \times \partial \mathbf{D}^2 \cong \mathbf{T}^2$ une suspension fixée $\sigma$ de type 0 ou 1. On note $\sigma_0$ le feuilletage de $\mathbf{T}^2$ par les cercles orientés de direction $(0, -1)$ (qui sont des méridiens du tore plein) et on considère le cône

$$C = C(\sigma_0, \sigma) \subset H_1(\mathbf{T}^2; \mathbf{R}) \setminus \{0\}$$

(cf. C) dont le bord gauche est la demi-droite $\mathbf{R}_+(0, -1)$. On note en outre $\widehat{B}$ l'ensemble des points entiers primitifs de $\partial E(\sigma_0, \sigma)$ qui appartiennent à une arête non verticale.

**Théorème 1.6.** *Soit $\sigma$ une suspension sur $\mathbf{T}^2$ de type 0 ou 1. On oriente $\sigma$ de telle sorte que la demi-droite de ses cycles asymptotiques soit le bord droit du cône $C$.*

*La classe d'Euler relative définit une application*

$$\chi_\partial \colon \mathcal{SCT}(\mathbf{S}^1 \times \mathbf{D}^2; \sigma) \longrightarrow H_1(\mathbf{S}^1 \times \mathbf{D}^2; \mathbf{Z}) \cong \mathbf{Z}$$

*dont les fibres sont toutes connexes. De plus, les sous-espaces $\mathcal{SCT}_u$ et $\mathcal{SCT}_v$ ont pour images respectives $\widehat{X}_u$ et $\widehat{X}_v$ où*

$$\widehat{X}_u = \{1 + 2i_*\gamma(Q) \mid Q \cap \operatorname{Int} \widehat{B} = \varnothing \text{ et } \operatorname{Card} Q = 0 \pmod{2},\ Q \subset \widehat{B}\},$$
$$\widehat{X}_v = \{1 + 2i_*\gamma(Q) \mid Q \cap \operatorname{Int} \widehat{B} \neq \varnothing \text{ et } \operatorname{Card} Q = 0 \pmod{2},\ Q \subset \widehat{B}\},$$

*et $i \colon \mathbf{Z}^2 = H_1(\mathbf{T}^2; \mathbf{Z}) \to H_1(\mathbf{S}^1 \times \mathbf{D}^2; \mathbf{Z}) = \mathbf{Z}$ est la projection sur le premier facteur.*

## E  Chirurgie

Soit $p$ et $q$ des entiers premiers entre eux tels que $0 < q < p$. Le nombre $p/q$ a un unique développement (fini) en fraction continue du type

$$\frac{p}{q} = a_0 - \cfrac{1}{a_1 - \cfrac{1}{\cdots - \cfrac{1}{a_n}}} \quad \text{avec } a_i \geq 2 \text{ pour } 0 \leq i \leq n.$$

On explique ci-dessous comment construire $\prod_{i=0}^n (a_i - 1)$ structures de contact holomorphiquement remplissables sur $\mathbf{L}_{p,q}$ qui sont non homotopes dans l'espace des champs de plans orientés. Le socle de cette construction est le fait [Ro] que $\mathbf{L}_{p,q}$ s'obtient à partir de $\mathbf{S}^3$ par chirurgie sur un entrelacs pondéré $(K_0, \ldots, K_n; -a_0, \ldots, -a_n)$ où les nœuds $K_i$ sont individuellement triviaux et forment une chaîne : chaque $K_i$ borde un disque plongé dont l'intersection avec $K_j$, $j \neq i$, est vide si $j \neq i \pm 1$ et réduite à un point sinon. On rappelle brièvement en quoi consiste la chirurgie.

Si $K$ est un nœud orienté dans $\mathbf{S}^3$ et $W$ un voisinage tubulaire de $K$, le groupe $H_1(\partial W; \mathbf{Z})$ possède une base canonique $(\mu, \lambda)$ définie comme suit :
 – $\mu$ (resp. $\lambda$) est la classe d'une courbe fermée simple sur $\partial W$ qui borde une surface dans $W$ (resp. dans $\mathbf{S}^3 \setminus \operatorname{Int} W$) ;
 – l'image de $\lambda$ dans $H_1(W; \mathbf{Z})$ est la classe de $K$ et l'intersection de $\mu$ avec $\lambda$ vaut 1 (le tore $\partial W$ est orienté comme bord de $W$).
Étant donné un nombre $r \in \mathbf{Q} \cup \{\infty\}$, la variété obtenue à partir de $\mathbf{S}^3$ par chirurgie sur le nœud pondéré $(K, r)$ est

$$(\mathbf{S}^3 \setminus W) \cup_\partial W'$$

où $W'$ est un tore plein et où le difféomorphisme de recollement $\partial W \to \partial W'$ envoie $r\mu + \lambda$ dans le noyau de l'application $H_1(\partial W'; \mathbf{Z}) \to H_1(W'; \mathbf{Z})$. La chirurgie sur un entrelacs est la chirurgie simultanée sur ses composantes.

Pour la suite, on oriente arbitrairement les nœuds $K_i$, $0 \le i \le n$, on en prend des voisinages tubulaires disjoints $W_i$ et on note $(\mu_i, \lambda_i)$ la base canonique de $H_1(\partial W_i; \mathbf{Z})$ décrite plus haut. À la description de $\mathbf{L}_{p,q}$ par chirurgie sur $(K_0, \ldots, K_n; -a_0, \ldots, -a_n)$, correspond une présentation du groupe $H_1(\mathbf{L}_{p,q}; \mathbf{Z}) \simeq \mathbf{Z}/p\mathbf{Z}$ : c'est le groupe engendré par les $\mu_i$ et défini par les relations

$$a_0 \mu_0 - \mu_1 = 0,$$
$$-\mu_{i-1} + a_i \mu_i - \mu_{i+1} = 0, \quad 1 \le i \le n-1,$$
$$-\mu_{n-1} + a_n \mu_n = 0.$$

On désigne alors par $\overline{\mu_i}$ l'image de $\mu_i$ dans $H_1(\mathbf{L}_{p,q}; \mathbf{Z})$ ; il est clair que $\overline{\mu_0}$ (de même que $\overline{\mu_n}$) est un générateur de ce groupe.

Soit maintenant $\zeta_0$ la structure de contact ordinaire sur $\mathbf{S}^3$, *i.e.* le champ des droites complexes tangentes à la sphère unité de $\mathbf{C}^2$. À tout nœud legendrien orienté $L$ dans $(\mathbf{S}^3, \zeta_0)$, on associe classiquement deux nombres :
 – l'*invariant de Thurston-Bennequin* $\mathtt{tb}(L)$ de $L$ est l'enlacement de $L$ avec $L + \nu$ où $\nu$ est la normale à $\zeta_0$ le long de $L$ ;
 – l'*invariant de Maslov* $\mathtt{m}(L)$ de $L$ est le nombre d'enroulement de la tangente orientée à $L$ dans une trivialisation quelconque de $\zeta_0$.

Ces nombres sont liés par (au moins) deux relations, dont la première est facile et la seconde est une forme du théorème de Bennequin :

$$\mathtt{tb}(L) + \big|\mathtt{m}(L)\big| = 1 \pmod{2},$$
$$\mathtt{tb}(L) + \big|\mathtt{m}(L)\big| \le -\chi(S),$$

où $S$ est une surface de Seifert quelconque de $L$. De plus, si deux entiers satisfont ces relations avec $\chi(S) = 1$, ce sont les invariants $\mathtt{tb}(L)$ et $\mathtt{m}(L)$ d'un nœud legendrien topologiquement trivial.

Dans ce contexte, le théorème de chirurgie de Y. ELIASHBERG [El4] s'énonce ainsi : si $L_0 \sqcup \cdots \sqcup L_n$ est un entrelacs legendrien orienté dans $(\mathbf{S}^3, \zeta_0)$, la variété $V$ obtenue à partir de $\mathbf{S}^3$ par chirurgie sur l'entrelacs pondéré $\big(L_0, \ldots, L_n; \mathtt{tb}(L_0) - 1, \ldots, \mathtt{tb}(L_n) - 1\big)$ porte une structure de contact orientée qui est holomorphiquement remplissable et dont la classe d'Euler vaut

$$\sum_{i=0}^{n} \mathtt{m}(L_i) \overline{\mu_i} \in H_1(V; \mathbf{Z}).$$

La proposition qui suit découle directement de ce résultat :

**Proposition 1.7.** *Soit $p$ et $q$ des entiers premiers entre eux tels que $0 < q < p$. Sur l'espace lenticulaire orienté $\mathbf{L}_{p,q}$, il existe au moins $\prod_{i=0}^{n}(a_i - 1)$ structures de contact orientées holomorphiquement remplissables qui, dans l'espace des champs de plans orientés, appartiennent à des classes d'homotopie distinctes.*

*Démonstration.* Étant donné des entiers $b_i \in [0, a_i - 2]$, $0 \le i \le n$, la chaîne $K_0 \sqcup \cdots \sqcup K_n$ est isotope à un entrelacs legendrien $L = L_0 \sqcup \cdots \sqcup L_n$ pour lequel

$$\begin{aligned}\mathtt{tb}(L_i) &= 1 - a_i, \\ \mathtt{m}(L_i) &= 2 - a_i + 2b_i,\end{aligned} \quad 0 \le i \le n.$$

La structure de contact $\zeta(b_1,\ldots,b_n)$ obtenue dans $\mathbf{L}_{p,q}$ par chirurgie de contact sur l'entrelacs legendrien pondéré $(L_0,\ldots,L_n;-a_0,\ldots,-a_n)$ est alors holomorphiquement remplissable. En outre, sa classe d'Euler est — en dualité de Poincaré avec —

$$\sum_{i=0}^{n}(2-a_i+b_i)\overline{\mu_i} \in H_1(\mathbf{L}_{p,q};\mathbf{Z}).$$

Mieux, le théorème 4.12 de [Go] montre que la classe

$$\sum_{i=0}^{n} b_i\overline{\mu_i} \in H_1(\mathbf{L}_{p,q};\mathbf{Z})$$

est un invariant homotopique du champ de plans orienté $\zeta(b_1,\ldots,b_n)$. Il reste donc à montrer que cette classe détermine la valeur des coefficients $b_i$. À cette fin, on observe que les relations qui définissent l'homologie de $\mathbf{L}_{p,q}$ à partir des méridiens $\mu_i$ se réécrivent sous la forme

$$\mu_i = p_i\mu_0, \qquad 0 \le i \le n,$$
$$0 = p_{n+1}\mu_0 \qquad \text{avec } p_0 = 0 \text{ et } p_{i+1} > (a_i-1)p_i \text{ pour } 0 \le i \le n.$$

Par suite,

$$\sum_{i=0}^{n} b_i\overline{\mu_i} = \Bigl(\sum_{i=0}^{n} b_i p_i\Bigr)\overline{\mu_0}.$$

Or les inégalités $p_{i+1} > (a_i-1)p_i$ montrent que les $(n+1)$-uplets d'entiers $b_i \in [0, a_i-2]$ sont lexicographiquement rangés dans le même ordre que les nombres $\sum_{i=0}^{n} b_i p_i \in \mathbf{Z}$. □

La proposition ci-dessus fournit indirectement beaucoup de structures de contact tendues sur le tore plein, le tore épais et les fibrés en tores au-dessus du cercle (qui sont revêtus par l'intérieur du tore épais). On cherche dès lors à montrer qu'on les a toutes obtenues.

## F  Tomographie

On discute ici brièvement les idées de l'étude tomographique sur le tore épais. Une structure de contact sur $\mathbf{T}^2 \times [0,1]$ imprime un feuilletage sur toutes les tranches $\mathbf{T}^2 \times \{t\}$ et fournit ainsi un « film » sur le tore, c'est-à-dire une famille à un paramètre $t$ de feuilletages. Ce film décrit entièrement la structure (lemme 2.1) mais, en contrepartie, est extrêmement sensible aux perturbations. Le but de l'étude tomographique est de dégager un petit nombre de films simples, déterminés par quelques données combinatoires, tels que toute structure de contact tendue sur $\mathbf{T}^2 \times [0,1]$ soit $\partial$-isotope à une structure qui imprime l'un de ces films.

Avant de présenter les modèles, on en donne deux exemples. Soit $\zeta$ et $\eta$ les structures de contact d'équations respectives

$$\begin{aligned} \sin(2\pi t)\,dx + \cos(2\pi t)\,dy &= 0, \\ \cos(2\pi x)\,dy - \sin(2\pi x)\,dt &= 0, \end{aligned} \qquad (x,y,t) \in \mathbf{T}^2 \times [-1,1].$$

Sur chaque tore $\mathbf{T}^2 \times \{t\}$, la structure $\zeta$ imprime un feuilletage linéaire de pente $-\tan(2\pi t)$ tandis que $\eta$ trace un feuilletage fixe qui a deux cercles de singularités — $\{x = 1/4\}$ et $\{x = 3/4\}$ — et dont les autres orbites sont les segments $\{y = \text{Const}\}$.

Soit maintenant $G_s$, $s \in [-1, 1]$, le graphe de la fonction $(x, y) \mapsto s(1/4 + \varepsilon \cos(2\pi x))$. Pour $s = \pm 1$, les feuilletages $\zeta G_s$ et $\eta G_s$ respectivement induits sur $G_s$ par $\zeta$ et $\eta$ sont topologiquement isotopes : ce sont des suspensions ayant deux orbites fermées hyperboliques — les cercles $\{x = 1/4\}$ et $\{x = 3/4\}$ — dont l'orientation dépend du signe de $s$. Les films $\zeta G_s$ et $\eta G_s$, $s \in [-1, 1]$, sont eux très différents. Le premier est une famille de suspensions dont la direction asymptotique tourne d'un demi-tour. Dans le second, les feuilletages $\eta G_s$, $s < 0$ (resp. $s > 0$), ont un type topologique constant et c'est le feuilletage singulier $\eta G_0$ qui assure la transition d'un type à l'autre ; autrement dit, le retournement de chaque orbite fermée $\{x = 1/2 \pm 1/4\}$ s'effectue par l'intermédiaire d'un cercle de singularités.

On reprend ci-dessous les notations et hypothèses du théorème 1.5 (cf. C). On pose en outre $F_t = \mathbf{T}^2 \times \{t\}$ pour tout $t \in [0, 1]$. On dira qu'une structure de contact $\xi$ sur $V = \mathbf{T}^2 \times [0, 1]$ est *sous forme normale*[1] s'il existe des points $t_1 < \cdots < t_n \in [0, 1]$ tels que les conditions suivantes soient remplies :
- pour tout $t \in [0, 1] \setminus \{t_1, \ldots, t_n\}$, le feuilletage $\xi F_t$ imprimé par $\xi$ sur $F_t$ est une suspension ;
- pour $1 \leq i \leq n$, le feuilletage $\xi F_{t_i}$ n'a aucune orbite fermée et ses singularités forment deux cercles le long desquels $\xi$ coïncide transversalement avec le champ des plans tangents aux tores $F_t$ ;
- sur chaque intervalle $]t_i, t_{i+1}[$, $1 \leq i \leq n-1$, les suspensions $\xi F_t$ ont une direction asymptotique $D_t$ qui varie.

Chaque feuilletage $\xi F_{t_i}$ apparaît ainsi dans le film $\xi F_t$ comme étape intermédiaire dans le retournement de deux orbites fermées. La droite $D_{t_i} \subset H_1(\mathbf{T}^2; \mathbf{R})$ qui contient la classe des cercles de singularités coïncide en particulier avec $D_t$ pour $t$ voisin de $t_i$. On appelle *lieu de retournement* de $\xi$ un ensemble $R(\xi)$ qui contient un vecteur entier primitif de chaque droite $D_{t_i}$. Il n'est pas très difficile de voir que cet ensemble, avec l'amplitude des séquences de rotation $[t_i, t_{i+1}]$ (c'est-à-dire l'angle total dont $D_t$ tourne sur $\mathbf{P}^1$ pour $t \in [t_i, t_{i+1}]$), déterminent $\xi$ à isotopie près. De plus, les parties 2 et 3 conduisent au résultat suivant :

**Proposition 1.8.** *Toute structure de contact tendue dans $\mathcal{SCT}(V; \sigma_0, \sigma_1)$ est $\partial$-isotope à une structure de contact sous forme normale $\xi$ dont le lieu de retournement $R(\xi)$ est inclus dans $B$. De plus :*
- *$\chi_\partial(\xi) = 2\gamma(R(\xi))$ ;*
- *$\xi$ est universellement tendue si et seulement si $R(\xi) \cap \operatorname{Int} B = \varnothing$ ;*
- *$\xi$ est virtuellement vrillée si et seulement si $R(\xi) \cap \operatorname{Int} B \neq \varnothing$ et, dans ce cas, toutes les droites $D_t$, $t \in ]0, 1[$, rencontrent $C$.*

*Enfin, si $\xi$ est virtuellement vrillée, tout ensemble fini $Q \subset B$ tel que $2\gamma(Q) = \chi_\partial(\xi)$ est le lieu de retournement d'une structure de contact sous forme normale $\partial$-isotope à $\xi$.*

## G  Fractions continues et enveloppes convexes

Il reste à relier les expressions numériques des théorèmes 1.1 et 1.3 aux configurations géométriques des théorèmes 1.5 et 1.6. Le lien vient des relations bien connues entre les fractions continues et les enveloppes convexes de points entiers [Ar] ; on en rappelle ci-dessous quelques aspects pertinents.

---

[1] Une définition techniquement plus satisfaisante sera donnée au terme de la proposition 3.22.

Soit $D_0$ et $D_1$ deux droites de $\mathbf{R}^2$ et $C$ un des deux cônes convexes bordés à gauche par $D_0$ et à droite par $D_1$. Ce cône peut être pris ouvert ou fermé du côté de $D_i$ mais ne contient pas 0. On pose :

$$C' = \{w' \in \mathbf{R}^2 \mid w' \wedge w > 0 \ \forall w \in C\}.$$

Ainsi, $C'$ est un cône convexe bordé à gauche par $D_1$ et à droite par $D_0$ ; en outre, $C'$ est ouvert (resp. fermé) du côté de $D_i$ si et seulement si $C$ est fermé (resp. ouvert) de ce même côté. On note aussi $E$ (resp. $E'$) l'enveloppe convexe des points entiers de $C$ (resp. $C'$) et $B$ (resp. $B'$) l'ensemble des points entiers primitifs de $\partial E$ (resp. $\partial E'$).

La dualité entre $C$ et $C'$ donne la propriété suivante (cf. lemme 3.37) : un point entier $w' \in C'$ est sur $\partial E'$ (resp. est un sommet de $\partial E'$) si et seulement si l'équation $w' \wedge w = 1$ a une solution entière (resp. plusieurs) dans $C$. Comme $C = (-C')'$, ces solutions sont sur $\partial E$ et on note $S(w') \subset \partial E$ leur ensemble (c'est un sommet ou l'ensemble des points entiers d'une arête de $\partial E$).

**Lemme 1.9.** *On pose $B' = \{w'_i, i \in I\}$ où $I$ est un intervalle de $\mathbf{Z}$ et $w'_i \wedge w'_{i+1} > 0$ pour tout $i \in I \setminus \sup I$. On note d'autre part $a_i - 1$ le cardinal de $S(w'_i)$ pour $i \in \operatorname{Int} I = I \setminus \partial I$.*
**a)** *La pente de $D_1$ dans chaque base entière $(w'_i, -w'_{i-1})$ vaut*

$$a_i - \cfrac{1}{a_{i+1} - \cfrac{1}{a_{i+2} - \cdots}}.$$

**b)** *Les combinaisons linéaires finies $\sum_{i \in I} b_i w'_i$, où $0 \leq b_i \leq a_i - 2$ pour tout $i \in \operatorname{Int} I$, sont deux à deux distinctes.*

*Démonstration.*
**a)** La dualité décrite plus haut montre que les extrémités de $S(w'_i)$ sont les points $w'_i - w'_{i-1}$ et $w'_{i+1} - w'_i$. Par suite,

$$w'_{i+1} + w'_{i-1} - 2w'_i = (a_i - 2)w'_i, \quad i.e. \quad w'_{i+1} + w'_{i-1} = a_i w'_i.$$

La matrice qui exprime la base $(w'_i, w'_{i+1})$ dans la base $(w'_{i-1}, w'_i)$ est donc $A_i = \begin{pmatrix} 0 & -1 \\ 1 & a_i \end{pmatrix}$. Celle qui exprime la base $(w'_{i+k}, w'_{i+k+1})$ dans la base $(w'_{i-1}, w'_i)$ est alors $A_i A_{i+1} \cdots A_{i+k}$ et l'homographie correspondante s'écrit

$$x \longmapsto h_k(x) = -\cfrac{1}{a_i - \cfrac{1}{\cdots - \cfrac{1}{a_{i+k} + x}}}.$$

Dans cette expression, $x$ et $h_k(x)$ sont les inverses des pentes dans la base $(w'_{i-1}, w'_i)$, ou encore les opposés des inverses des pentes dans la base $(w'_i, -w'_{i-1})$. Par ailleurs, lorsque $k$ tend vers l'infini, l'image de la droite dirigée par $w'_i$ tend vers $D_1$. On obtient ainsi la formule promise.

**b)** Les projections des points $w'_i$ sur $D_1$ parallèlement à $D_0$ donnent, après choix d'un vecteur sur $D_1$, des réels $x_i$ distincts, ordonnés et strictement positifs pour $i > \inf I$. Ces nombres vérifient de plus

$$a_i x_i = x_{i-1} + x_{i+1} \quad \text{pour tout } i \in \operatorname{Int} I.$$

Ainsi, $x_{i+1} > (a_i - 1)x_i$ pour tout $i \in \operatorname{Int} I$. Comme dans la preuve de la proposition 1.7, ces relations montrent que les familles finies de coefficients $b_i$, $i > \inf I$, vérifiant $0 \leq b_i \leq a_i - 2$ pour $i < \sup I$ sont lexicographiquement rangés dans le même ordre que les nombres $\sum_{i > \inf I} b_i x_i$. Le résultat en découle immédiatement. $\square$

**Corollaire 1.10.** *Soit $w'_0, \ldots, w'_n$, $0, n \in \operatorname{Int} I$, des points consécutifs sur $B'$. Le nombre de valeurs prises par la fonction $\gamma$ sur les parties $Q$ de cardinal pair (resp. impair) dans $\widehat{B} = \bigcup_{i=0}^n S(w'_i) \subset B$ est $\prod_{i=0}^n (a_i - 1)$.*

*Démonstration.* Soit $w_n = w'_{n+1} - w'_n$ l'extrémité droite de $S(w'_n)$. L'application qui, à toute partie $Q \subset \widehat{B}$, associe l'ensemble

$$Q' = \begin{cases} Q \setminus \{w_n\} & \text{si } w_n \in Q, \\ Q \cup \{w_n\} & \text{si } w_n \notin Q, \end{cases}$$

définit une bijection de l'ensemble des parties de cardinal impair dans celui des parties de cardinal pair. De plus, $\gamma(Q') = \gamma(Q) - w_n$ de sorte que $\gamma$ prend autant de valeurs sur les parties de cardinal impair que sur celles de cardinal pair. D'autre part, pour toute partie $Q$ de cardinal pair dans $\widehat{B}$, la relation de Chasles montre que $\gamma(Q) = \sum_{i=0}^n b_i w'_i$ avec $0 \leq b_i \leq a_i - 2$. Le corollaire résulte donc du lemme 1.9-b). $\square$

**Exemple 1.11.** Soit $D_0$ et $D_1$ les droites respectivement engendrées par les vecteurs $(0,1)$ et $(p,q)$, où $p$ et $q$ sont des entiers premiers entre eux tels que $0 < q < p$. On note $\widehat{B}$ l'ensemble des points de $B$ situés sur une arête de longueur finie de $\partial E$. Le nombre de valeurs prises par la fonction $\gamma$ sur les parties de cardinal pair dans $\widehat{B}$ est $\prod_{i=0}^n (a_i - 1)$ où

$$\frac{p}{q} = a_0 - \cfrac{1}{a_1 - \cfrac{1}{\cdots - \cfrac{1}{a_n}}} \quad \text{avec } a_i \geq 2 \text{ pour } 0 \leq i \leq n.$$

## 2 Tomographie des structures de contact

### A Lemme de reconstruction

$F$ désigne une surface close orientée, $V$ le produit $F \times [-1, 1]$ et $F_t$, $t \in [-1, 1]$, la surface $F \times \{t\}$. Les structures de contact qu'on considère sur $V$ sont toutes directes pour l'orientation produit et (co)orientables. En particulier, chacune d'elles est définie globalement par une équation de Pfaff $\alpha = 0$. Dans la suite, lorsqu'on choisit une telle équation, on munit la structure de contact de l'orientation induite par $d\alpha$.

Les lemmes de cette partie reposent sur l'interprétation de la condition $\alpha \wedge d\alpha > 0$ lorsqu'on exprime $\alpha$ sous la forme

$$\alpha = \beta_t + u_t \, dt, \qquad t \in [-1, 1],$$

où les $\beta_t$ et les $u_t$ sont respectivement des 1-formes et des fonctions sur $F$. La *condition de contact* s'écrit

$$u_t \, d\beta_t + \beta_t \wedge \left( du_t - \dot{\beta}_t \right) > 0. \tag{$*$}$$

Le signe fait ici référence à l'orientation de $F$ et $\dot\beta_t$ est la dérivée de $\beta_t$ par rapport à $t$.

Dans ce contexte, le feuilletage caractéristique $\xi F_t$ de $F_t$ est défini par l'équation $\beta_t = 0$. D'autre part, la famille des formes $\beta_t$ étant donnée, les familles de fonctions $u_t$ qui vérifient l'inégalité (∗) constituent un ensemble convexe. Le premier lemme ci-dessous découle alors directement du théorème de J. Gray :

**Lemme 2.1.** *Si deux structures de contact sur $V$ impriment les mêmes feuilletages sur chaque surface $F_t$, elles sont isotopes relativement au bord.* □

Ainsi, une structure de contact sur $V$ est entièrement déterminée par le *film* des feuilletages $\xi F_t$.

Cet énoncé recèle cependant une petite subtilité. On réserve en effet le terme de *feuilletage* pour désigner un objet géométrique, à savoir l'ensemble des courbes intégrales d'une 1-forme. Mais il résulte facilement du théorème des fonctions implicites que deux 1-formes qui définissent le même feuilletage et ont une différentielle non nulle en leurs singularités sont multiples l'une de l'autre par une fonction partout non nulle.

Dans la suite, lorsque $\xi$ est (co)orientée par le choix d'une équation de Pfaff $\alpha = 0$ où $\alpha = \beta_t + u_t\,dt$, on regarde les orbites des feuilletages $\xi F_t$ comme des courbes (co)orientées : la relation $\beta = Y \lrcorner \omega$ met de fait en dualité les 1-formes et les champs de vecteurs sur $F$ moyennant le choix d'une forme d'aire $\omega$. Dans cette situation, les singularités de $\xi F_t$ — points où la structure $\xi$ est tangente à $F_t$ — sont affectées d'un signe, positif ou négatif selon que les orientations de $\xi$ et $F$ coïncident ou non. Ce signe est aussi celui de $d\beta_t$ et de la divergence du champ de vecteurs $Y$ dual de $\beta_t$. La condition de contact force donc les singularités à avoir une divergence non nulle. Or un peu de théorie des formes normales montre que ces singularités ne présentent — à un ensemble de codimension infinie près — que trois types topologiques déterminés par leur indice, $1$, $-1$ ou $0$. Les singularités *positives* d'indice $1$ sont ainsi des foyers (ou des nœuds) *répulsifs*.

L'influence de la géométrie complexe sur l'étude des structures de contact [El3] est à l'origine d'un vocabulaire particulier. Une singularité de $\xi F_t$ est ainsi appelée *point complexe* de $F_t$ et, selon que son indice vaut $1$, $-1$ ou $0$, le point est dit *elliptique*, *hyperbolique* ou *parabolique*. Lorsque $\xi F_t$ est non singulier — ce qui suppose que $F$ soit un tore —, $F_t$ est dit *totalement réel*.

Une autre conséquence immédiate du théorème de Gray est que deux films qui sont $\mathcal{C}^1$-proches et ont même début et même fin représentent la même structure. On peut donc imposer au film $\xi F_t$ toute propriété qui est générique pour les familles à un paramètre de feuilletages sur une surface.

**Propriétés 1.** Génériquement, les singularités et les orbites fermées de chaque feuilletage $\xi F_t$ sont isolées — c'est vrai pour les familles indexées par un nombre fini quelconque de paramètres. Les phénomènes suivants jouent alors un rôle important :

$1_k$) $\xi F_t$ possède exactement $k$ orbites fermées *dégénérées*, c'est-à-dire dont l'application de premier retour est tangente à l'identité ;

$2_k$) $\xi F_t$ possède exactement $k$ orbites reliant l'une à l'autre deux singularités d'indices négatifs ou nuls (topologiquement, des selles ou des nœuds-selles).

Dans un film, c'est-à-dire une famille à un paramètre de feuilletages, $1_k$ et $2_k$ n'apparaissent génériquement jamais avec $k > 1$. On observe en revanche $1_0$ et $2_0$ durant un ouvert dense $\Omega \subset [-1, 1]$. Le reste du temps, on observe $1_1$ ou $2_1$ mais jamais simultanément et sous des formes bénignes : pour les orbites fermées dégénérées, l'application de premier retour a une dérivée seconde non nulle et les orbites reliant des singularités d'indices négatifs ou nuls sont de simples connexions de selles non dégénérées.

# B  Lemme de réalisation

Pour comprendre la spécificité des films de contact, on commence par montrer qu'un film dont toutes les images vérifient la propriété discutée ci-dessous s'intègre en une (unique) structure de contact.

Reprenant la terminologie de [Gi1], on dit qu'une surface $F_t$ est *convexe* dans $(V, \xi)$ si elle possède un voisinage tubulaire trivialisé $U \simeq F_t \times \mathbf{R}$ dans lequel la structure $\xi$ est $\mathbf{R}$-invariante — la trivialisation de $U$ est étrangère à la structure produit sur $V$.

Cette propriété ne dépend que du germe de $\xi$ le long de $F_t$, *i.e.* du feuilletage $\xi F_t$ ou de la forme $\beta_t$ qui le définit. Par ailleurs, une structure de contact $\mathbf{R}$-invariante sur $F_t \times \mathbf{R}$ admet une équation tout aussi invariante. Il s'ensuit alors de la condition de contact $(*)$ qu'une surface $F_t$ est convexe si et seulement s'il existe une fonction $v_t$ sur $F_t$ qui vérifie

$$v_t \, d\beta_t + \beta_t \wedge dv_t > 0, \qquad (**)$$

condition dont on donne ci-dessous la traduction géométrique.

**Définition 2.2.** On dira que le feuilletage $\xi F_t$ est *expansif* sur un domaine $G$ de $F_t$ s'il se laisse définir par une 1-forme qui est positive sur le bord orienté $\partial G$ et dont la différentielle est partout positive sur $G$. Dualement dit, $\xi F_t$ est expansif s'il est dirigé par un champ de vecteurs sortant dont le flot dilate l'aire.

On dira qu'une multi-courbe $\Gamma_t$ (sous-variété close de dimension 1) *scinde* le feuilletage $\xi F_t$ si elle évite les singularités, coupe transversalement les orbites et découpe $F_t$ en régions où $\xi F_t$ est expansif. Une telle multi-courbe, si elle existe, est unique à isotopie près parmi les multi-courbes qui scindent. En outre, l'espace des feuilletages scindés par une multi-courbe donnée est contractile.

Si $F_t$ est une surface convexe, les zéros de n'importe quelle fonction $v_t$ vérifiant $(**)$ forment une multi-courbe qui scinde $\xi F_t$ car

$$\begin{cases} \text{là où } v_t = 0, & \beta_t \wedge dv_t > 0, \\ \text{là où } v_t \neq 0, & v_t^2 \, d\Big(\dfrac{\beta_t}{v_t}\Big) > 0 \,. \end{cases}$$

Inversement, s'il existe une multi-courbe $\Gamma_t$ qui scinde $\xi F_t$, il est facile de construire une solution $v_t$ de $(**)$ qui s'annule exactement sur $\Gamma_t$ [Gi1].

**Lemme 2.3.** *Toute famille à un paramètre de feuilletages scindés par des multi-courbes est le film d'une structure de contact. Précisément, soit $\beta_t$, $t \in [-1, 1]$, une famille de 1-formes sur $F$. On suppose que, pour tout $t$, il existe une fonction $v_t \colon F \to \mathbf{R}$ vérifiant la condition $(**)$ :*

$$v_t \, d\beta_t + \beta_t \wedge dv_t > 0 \,.$$

*(On ne demande a priori aucune continuité de $v_t$ en fonction de $t$.) Il existe alors sur $V = F \times [-1, 1]$ une structure de contact d'équation $\beta_t + u_t \, dt = 0$, où $u_t$ est une famille (lisse) de fonctions $F \to \mathbf{R}$.*

*Démonstration.* Si les fonctions $v_t$ forment une famille lisse en $t$, on pose simplement $u_t = \lambda v_t$ où $\lambda$ est un nombre positif assez grand. Sinon, une partition de l'unité sur $[-1, 1]$ permet de gagner cette régularité car la condition $(**)$ définit, $\beta_t$ étant donnée, un cône convexe ouvert de fonctions $v_t$. □

Ce lemme et ceux qui suivent en C sont utiles car les surfaces convexes abondent. L'énoncé ci-dessous — dont la preuve est donnée dans [Gi1, proposition II.2.6] — montre qu'elles sont génériques et les caractérise sous une hypothèse faible (cf. lemme 2.9).

On dit ci-dessous qu'un feuilletage $\xi F_t$ vérifie la *propriété de Poincaré-Bendixson* si l'ensemble limite de toute demi-orbite est soit une singularité, soit une orbite fermée (cycle), soit une réunion finie de singularités et d'orbites les reliant (poly-cycle). Ainsi, le théorème de Poincaré-Bendixson affirme qu'un feuilletage de la sphère ou d'une région plane vérifie cette propriété dès que ses singularités sont isolées.

**Proposition 2.4.** *Soit $F_t$ une surface dont le feuilletage caractéristique $\xi F_t$ vérifie la propriété de Poincaré-Bendixson. $F_t$ est convexe si et seulement si les conditions suivantes sont remplies[2] :*
- *aucune orbite fermée de $\xi F_t$ n'est dégénérée ;*
- *aucune orbite de $\xi F_t$ ne va d'une singularité négative à une singularité positive.*

Cette proposition décrit les deux causes principales de non convexité. Il y a cependant un autre exemple très intéressant de surfaces non convexes (dont le feuilletage ne vérifie pas la propriété de Poincaré-Bendixson) :

**Définition 2.5.** On dit qu'une surface contenue dans une variété de contact est *pré-lagrangienne* si son feuilletage caractéristique se laisse définir globalement par une 1-forme fermée non singulière. — La terminologie vient de ce qu'une telle surface se relève en une surface lagrangienne dans la variété symplectisée. —

Une telle surface, si elle est connexe et close, est un tore ou une bouteille de Klein. Le feuilletage caractéristique d'un tore pré-lagrangien est différentiablement conjugué à un feuilletage linéaire. Le tore possède donc un voisinage tubulaire produit $\mathbf{R}^2/\mathbf{Z}^2 \times [-\varepsilon, \varepsilon]$, muni de coordonnées $(x_1, x_2, \theta)$, dans lequel la structure a une équation du type

$$\cos(\theta_0 + \theta)\, dx_1 - \sin(\theta_0 + \theta)\, dx_2 = 0\,.$$

Ce voisinage est ainsi entièrement feuilleté par des tores pré-lagrangiens dont les feuilletages caractéristiques ont des directions qui varient dans tout un intervalle.

Pour une bouteille de Klein pré-lagrangienne, le modèle local est plus rigide. C'est le quotient de $\mathbf{T}^2 \times [-\varepsilon, \varepsilon]$ muni de la structure d'équation

$$\cos\theta\, dx_1 - \sin\theta\, dx_2 = 0$$

par la transformation $(x, \theta) \mapsto (x_1 + 1/2, -x_2, -\theta)$. En effet, comme tout feuilletage non singulier de $\mathbf{K}^2$, le feuilletage caractéristique d'une bouteille pré-lagrangienne a des orbites fermées et est donc formé de cercles coorientables. Les tores d'équation $|\theta| = \text{Const} \neq 0$ qui entourent la bouteille de Klein sont aussi pré-lagrangiens.

## C Lemme d'unicité

Il s'agit encore d'un corollaire direct du théorème de Gray qui permet de substituer par un modèle toute séquence du film formée uniquement d'images convexes.

**Lemme 2.6.** *Soit $\xi_0$ et $\xi_1$ deux structures de contact qui impriment les mêmes feuilletages sur le bord. Si chaque surface $F_t$ est convexe pour les deux structures et contient une multi-courbe $\Gamma_t$ qui varie continûment avec $t$ et scinde à la fois $\xi_0 F_t$ et $\xi_1 F_t$, alors $\xi_0$ et $\xi_1$ sont isotopes relativement au bord.*

---

[2]Contrairement aux apparences, ces conditions sont indépendantes de l'orientation du feuilletage $\xi F_t$ choisie implicitement.

*Démonstration.* Si $\xi_i$ a pour équation
$$\beta_t^i + u_t^i \, dt = 0, \qquad i \in \{0, 1\},$$
les hypothèses fournissent, pour tout $t \in [-1, 1]$, des fonctions $v_t^0$ et $v_t^1$ qui s'annulent exactement sur $\Gamma_t$ et vérifient
$$v_t^i \, d\beta_t^i + \beta_t^i \wedge dv_t^i > 0, \qquad i = 0, 1\,.$$
Ces fonctions sont alors multiples l'une de l'autre par une fonction positive $h_t = v_t^1/v_t^0$ de sorte que, quitte à multiplier initialement l'équation de $\xi_1$ par $h_t$, on peut prendre $v_t^1 = v_t^0 = v_t$. Les équations
$$\beta_t^i + \left((1-s)u_t^i + s\lambda v_t\right) dt = 0, \quad \text{pour } i = 0, 1, \ s \in [0,1] \text{ et } \lambda > 0,$$
définissent toutes des structures de contact. De même pour les équations
$$(1-s)\beta_t^0 + s\beta_t^1 + \lambda v_t \, dt = 0$$
— si du moins $\lambda$ est assez grand —, ce qui fournit l'homotopie cherchée. □

**Remarque 2.7.** Étant donné une structure de contact $\xi_0$ sur $V$ pour laquelle toutes les surfaces $F_t$ sont convexes, on regarde l'espace des structures de contact $\xi$ ayant les propriétés suivantes :
– $\xi$ coïncide avec $\xi_0$ sur le bord ;
– chaque surface $F_t$ est convexe pour $\xi$ et contient une multi-courbe qui varie continûment avec $t$ et scinde à la fois $\xi_0 F_t$ et $\xi F_t$.

Le lemme 2.6 affirme que cet espace est connexe mais sa preuve montre qu'il est en fait contractile. La version « sans bord » qu'on énonce et démontre ci-dessous supporte tout aussi bien les paramètres.

**Lemme 2.8.** *Soit $\xi_0$ et $\xi_1$ deux structures de contact sur $F \times \mathbf{R}$ pour lesquelles toutes les surfaces $F_t = F \times \{t\}$ sont convexes. Si les multi-courbes $\Gamma^0$ et $\Gamma^1$ qui scindent respectivement $\xi_0 F_0$ et $\xi_1 F_0$ sont isotopes, alors $\xi_0$ et $\xi_1$ sont isotopes.*

*Démonstration.* Quitte à bouger chaque structure par une isotopie horizontale, on peut supposer que $\Gamma^0 = \Gamma^1 = \Gamma$ et que, pour tout réel $t$, la courbe $\Gamma \times \{t\}$ scinde à la fois $\xi_0 F_t$ et $\xi_1 F_t$. Des arguments semblables à ceux invoqués dans le lemme 2.6 montrent que, à de nouvelles isotopies horizontales près laissant fixe $\Gamma \times \mathbf{R}$, les structures $\xi_0$ et $\xi_1$ ont des équations du type
$$\beta_t^i + u \, dt = 0, \qquad i = 0, 1,$$
où $u$ est une fonction $F \to \mathbf{R}$ qui s'annule exactement sur $\Gamma$ et vérifie
$$u \, d\beta_t^i + \beta_t^i \wedge du > 0 \quad \text{pour } i = 0, 1 \text{ et tout } t \in \mathbf{R}.$$

Comme les solutions $\beta$ sur $F$ de l'inégalité $u \, d\beta + \beta \wedge du > 0$ forment un convexe, il est facile d'en trouver une famille $\beta_t$, $t \in \mathbf{R}$, telle qu'en outre
$$\begin{cases} \beta_t = \beta_t^0 & \text{pour tout } t \text{ proche d'un entier pair,} \\ \beta_t = \beta_t^1 & \text{pour tout } t \text{ proche d'un entier impair.} \end{cases}$$

En multipliant au besoin $\beta_t$ par $\lambda(t)$, où $\lambda$ est une fonction égale à 1 au voisinage de **Z**, on peut imposer

$$u\, d\beta_t + \beta_t \wedge \left(du - \dot\beta_t\right) > 0,$$

ce qui assure que l'équation

$$\beta_t + u\, dt = 0$$

définit une structure de contact $\xi$ sur $F \times \mathbf{R}$. Toutes les surfaces $F_t$ sont convexes et, pour tout entier $k \in \mathbf{Z}$, la structure $\xi$ coïncide avec $\xi_0$ près de $F_{2k}$ et avec $\xi_1$ près de $F_{2k+1}$. Le lemme 2.6, combiné à un argument d'exhaustion, donne le résultat voulu. □

## D  Lemme de banalisation dynamique

La propriété de Poincaré-Bendixson est satisfaite sur un ensemble qui est de codimension exactement 1 dans l'espace des feuilletages sur une surface autre que la sphère. Le lemme qui suit permet cependant de l'imposer pour les feuilletages caractéristiques des surfaces $F_t$ non convexes.

**Lemme 2.9.** *Soit $\xi_0$ une structure de contact sur $V = F \times [-1, 1]$. Si $F_{-1}$ et $F_1$ sont convexes, $\xi_0$ est isotope relativement au bord à une structure de contact $\xi$ pour laquelle le feuilletage caractéristique de toutes les surfaces $F_t$ non convexes vérifie la propriété de Poincaré-Bendixson.*

*Démonstration.* Pour toute partie $P$ de $F$, on note comme d'habitude $P_t$, $t \in [-1, 1]$, l'ensemble $P \times \{t\}$.

On suppose d'abord que $\xi_0$, arbitrairement orientée, vérifie les deux propriétés suivantes :

1) $F_t$ est une surface convexe pour tout $t \in [-1, -1/2] \cup [1/2, 1]$ ;
2) le feuilletage $\xi_0 F_t$ est indépendant de $t$ sur chacun des intervalles $[-3/4, -1/2]$ et $[1/2, 3/4]$ ;
3) $F$ est découpée en deux surfaces planes $F^+$ et $F^-$, éventuellement non connexes, dont le bord commun, relevé dans $F_{\pm 1/2}$, est transversal au feuilletage $\xi_0 F_{\pm 1/2}$, les orbites sortant de $F^+_{\pm 1/2}$ pour entrer dans $F^-_{\pm 1/2}$.

On désigne par $G^+$ et $G^-$ des rétrécissements respectifs de $F^+$ et $F^-$ dont les bords sont isotopes à $\partial F^\pm$ parmi les courbes ayant, dans $F_{\pm 1/2}$, des relevés transversaux à $\xi_0 F_{\pm 1/2}$.

On choisit alors une fonction strictement croissante $g \colon [0, 1] \to [1/2, 1]$ qui induit l'identité sur $[3/4, 1]$ puis on pose $\xi_1 = \phi_1^* \xi_0$ où $\phi_s$, $s \in [0, 1]$, est une isotopie à support dans $F \times [-3/4, 3/4]$ qui déplace les points verticalement et vérifie

$$\phi_1(G_t^+) = G_{g(t)}^+ \qquad \text{pour tout } t \geq 0,$$
$$\phi_1(G_t^-) = G_{-g(t)}^- \qquad \text{pour tout } t \leq 0.$$

Le théorème de Poincaré-Bendixson assure que, pour tout $t \in [-3/4, 3/4]$, le feuilletage $\xi_1 F_t \simeq \xi_0 \phi_1(F_t)$ vérifie la propriété de Poincaré-Bendixson. En effet, si $t \in [0, 3/4]$ par exemple, la courbe $\partial G_t^+$ découpe $F_t$ en régions planes qui, vu les propriétés 2) et 3), sont invariantes par $\xi_1 F_t$ dans le passé ou le futur. D'autre part, $\xi_1$ coïncide avec $\xi_0$ hors de $F \times [-3/4, 3/4]$ si bien que, d'après la propriété 1), toutes les surfaces $F_t$ avec $|t| \geq 3/4$ sont convexes pour $\xi_1$.

Reste à mettre en place les propriétés 1), 2) et 3) en déformant $\xi_0$ par une isotopie relative au bord. On trace d'abord sur $F$ une multi-courbe $K$ dont le complémentaire est formé de régions planes. Comme les surfaces $F_{-1}$ et $F_1$ sont convexes, elles sont découpées en régions où leur feuilletage caractéristique est expansif. Quitte à bouger $K$ sur $F$ par une isotopie, on peut supposer qu'aucune composante connexe de $K_{\pm 1}$ n'est incluse dans une seule région d'expansion de $F_{\pm 1}$. Les lemmes sur les surfaces convexes (voir aussi la proposition II.3.6 et l'exemple II.3.7 de [Gi1] fournissent alors une isotopie de $V$ qui met en place les propriétés 1), 2) et :

4) la multi-courbe $K_{\pm 1/2}$ est réunion d'orbites et de singularités de $\xi_0 F_{\pm 1/2}$ — plus précisément, chaque arc d'intersection avec une région en expansion contient une selle ou un puits.

Il suffit dès lors de prendre pour $F^-$ un voisinage tubulaire de $K$ et pour $F^+$ l'adhérence du complémentaire. □

**Remarque 2.10.** Soit $\mathcal{P}$ une propriété générique quelconque pour les familles à un paramètre de feuilletages sur une surface. Dans le lemme ci-dessus, on peut aussi imposer $\mathcal{P}$ à tous les feuilletages $\xi_1 F_t$. En effet, la convexité est une propriété stable par perturbation et la propriété de Poincaré-Bendixson est garantie par une propriété également stable, à savoir l'existence de courbes transversales qui cantonnent la dynamique dans des régions planes.

## E  Lemme de naissance-mort

Le défaut de convexité d'une surface $F_0$ indique, par définition-même, une variation du type *différentiable* des feuilletages $\xi F_t$ quand $t$ traverse la valeur 0. Les deux prochains lemmes montrent que la condition de contact force en fait une bifurcation beaucoup plus nette : un changement du type *topologique* des feuilletages.

On suppose ici que le feuilletage $\xi F_0$ possède une orbite fermée $C_0$ *dégénérée*, c'est-à-dire dont l'application de premier retour est tangente à l'identité. On dira que $C_0$ est *positive* (resp. *négative*) si la dérivée seconde de cette application est non nulle et a le signe annoncé. On note que ce signe ne dépend pas de l'orientation de $\xi$ choisie car, lorsqu'on la renverse, on change à la fois l'orientation et la coorientation de $\xi F_0$.

**Lemme 2.11.** *Toute orbite fermée dégénérée $C_0$ qui est positive (resp. négative) indique la naissance (resp. la mort) instantanée de deux orbites fermées non dégénérées dans la famille $\xi F_t$. Autrement dit, au voisinage de $C_0$ et pour $t$ petit, $\xi F_t$ présente deux orbites fermées — qui sont de plus non dégénérées — si $\operatorname{sgn}(C_0) t > 0$ et n'en possède aucune si $\operatorname{sgn}(C_0) t < 0$.*

*Démonstration.* Comme $C_0$ est une orbite fermée *dégénérée*, $\xi$ admet une équation du type $\beta_t + u_t \, dt = 0$ pour laquelle $d\beta_0$ s'annule identiquement sur $C_0$. On choisit alors sur $F \simeq F_0$ un arc $A$ transversal à $\xi F_0$ et coupant $C_0$ en un point $a_0$ de son intérieur. Pour $t$ assez petit, on note $\pi_t$ l'application de premier retour de $\xi F_t$ sur $A$ et on pose $a_t = \pi_t(a_0)$. L'observation clé est que la condition de contact impose l'inégalité

$$\beta_0(\dot a_0) < 0 \quad \text{où} \quad \dot a_0 = \frac{da_t}{dt}\bigg|_{t=0}. \tag{\dag}$$

Soit en effet $C_t$ la courbe fermée simple orientée obtenue en suivant l'orbite de $\xi F_t$ entre $a_0$ et $a_t$ puis en bouclant par $A$. Le cycle $C_t \cup (-C_0)$ est le bord d'une chaîne singulière $D_t$ contenue dans un $\delta t$-voisinage de $C_0$ pour une certaine constante $\delta$. Comme $d\beta_t$ est nulle sur $C_0$ pour $t = 0$,

$$\left| \int_{D_t} d\beta_t \right| < \text{Const } t^2.$$

Par ailleurs,
$$\int_{D_t} d\beta_t = \int_{C_t} \beta_t - \int_{C_0} \beta_t$$
et le premier de ces deux termes se réduit à une intégrale sur le segment orienté $A_t$ de $A$ qui va de $a_t$ à $a_0$. Ainsi,
$$0 = \frac{d}{dt}\Big(\int_{D_t} d\beta_t\Big)\Big|_{t=0} = \frac{d}{dt}\Big(\int_{A_t} \beta_t\Big)\Big|_{t=0} - \frac{d}{dt}\Big(\int_{C_0} \beta_t\Big)\Big|_{t=0}$$
donc
$$\beta_0(\dot{a}_0) = -\frac{d}{dt}\Big(\int_{A_t} \beta_t\Big)\Big|_{t=0} = -\frac{d}{dt}\Big(\int_{C_0} \beta_t\Big)\Big|_{t=0} = -\int_{C_0} \dot{\beta}_0.$$
Or, le long de $C_0$, la condition de contact s'écrit
$$\beta_0 \wedge (du_0 - \dot{\beta}_0) > 0,$$
ce qui entraîne, compte tenu de l'orientation choisie pour les feuilletages caractéristiques,
$$\int_{C_0} \dot{\beta}_0 > \int_{C_0} du_0 = 0.$$

Reste à voir que l'inégalité (†) force la bifurcation. Pour fixer les idées, on suppose que l'orbite $C_0$ est positive — c'est-à-dire que la dérivée seconde de $\pi_0$ est positive en $a_0$ — et on identifie $A$ à un intervalle autour de $0$ en intégrant $\beta_0$ à partir de $a_0$. On peut ainsi écrire
$$\pi_t(x) = x - t\lambda(x,t) + x^2\mu(x,t)$$
où les fonctions $\lambda$ et $\mu$ sont strictement positives sur un voisinage $[-\varepsilon, \varepsilon] \times [-\varepsilon, \varepsilon]$. Pour $t < 0$, l'équation
$$t\lambda(x,t) - x^2\mu(x,t) = 0$$
n'a pas de solutions si bien que $\pi_t$ n'a pas de points fixes. En revanche, pour $t > 0$ assez petit, on a, comme en $t = 0$,
$$\pi_t(-\varepsilon) > -\varepsilon \quad \text{et} \quad \pi_t(\varepsilon) > \varepsilon$$
tandis que
$$\pi_t(0) = -t\lambda(0,t) < 0.$$
Par suite, $\pi_t$ a un point fixe $x_-$ dans l'intervalle $]-\varepsilon, 0[$ et un autre $x_+$ dans $]0, \varepsilon[$. De plus, comme la fonction $\pi_t - \text{id}$ est strictement convexe, sa dérivée première s'annule au plus une fois. Or elle s'annule forcément entre deux zéros consécutifs, donc $x_-$ et $x_+$ sont les seuls points fixes de $\pi_t$ et sont non dégénérés, respectivement attractif et répulsif.  $\square$

**Remarque 2.12.** Le fait que l'orbite fermée $C_0$ soit peu dégénérée (*i.e.* ait une application de premier retour dont la dérivée seconde est non nulle) n'intervient ni pour établir l'inégalité (†), ni pour écrire l'équation des points fixes de $\pi_t$ sous la forme
$$t\lambda(x,t) - x^2\mu(x,t) = 0.$$
Autrement dit, si $C_0$ est une orbite fermée dégénérée quelconque, $\mu$ peut s'annuler à un ordre arbitraire en $(0,0)$ mais $\lambda(0,0)$ est strictement positif de sorte que l'équation ci-dessus reste régulière et définit encore une sous-variété locale. Ainsi, la réunion des orbites fermées proches de $C_0$ sur les différents tores $F_t$, $t \in [-\varepsilon, \varepsilon]$, forme une surface au voisinage de $C_0$, quel que soit le degré de dégénérescence de $C_0$.

# F  Lemme de croisement

On suppose ici que $\xi F_0$ possède une *connexion de selles* $C_0$, c'est-à-dire une orbite (orientée) qui relie deux selles. On dira que la connexion $C_0$ est *rétrograde* si elle part d'une selle $b_0^-$ qui est négative et aboutit à une selle $b_0^+$ qui est positive. Encore une fois, cette notion est indépendante de l'orientation de $\xi$ choisie.

Pour $t$ proche de 0, chaque feuilletage $\xi F_t$ possède exactement une selle $b_t^\pm$ au voisinage de $b_0^\pm$. A l'instant $t = 0$, la présence de la connexion $C_0$ signifie qu'une des séparatrices instables de $b_0^-$ coïncide avec une des séparatrices stables de $b_0^+$. Lorsque $t$ varie au voisinage de 0, on peut suivre continûment les séparatrices en question et un arc $A$, transversal à $\xi F_0$ et coupant $C_0$ en un unique point $a_0$, coupe encore la séparatrice stable (resp. instable) de $b_t^+$ (resp. $b_t^-$) en un seul point $a_t^+$ (resp. $a_t^-$). Grâce à la coorientation de $\xi$, on peut comparer sur $A$ les positions de $a_t^+$ et $a_t^-$.

**Lemme 2.13.** *Toute connexion de selles rétrograde $C_0$ indique le croisement instantané des séparatrices concernées dans la famille $\xi F_t$. Précisément, quand $t$ est positif (resp. négatif), la séparatrice stable passe au-dessus (resp. au-dessous) de la séparatrice instable.*

*Démonstration.* Si $\beta_t + u_t\, dt = 0$ est une équation de $\xi$, la fonction $u_0$ est négative en $b_0^-$ et positive en $b_0^+$ ; elle s'annule donc en un point $a_0$ de $C_0$. On va observer le mouvement des séparatrices sur un arc $A$ de $F \cong F_0$ transversal à $\xi F_0$ et coupant $C_0$ au seul point $a_0$. Pour cela, on désigne par $S_t^\pm$ la partie de la séparatrice de $b_t^\pm$ comprise entre l'arc $A$ et $b_t^\pm$ — en particulier, $C_0 = S_0^+ \cup S_0^-$ — et on pose $a_t^\pm = A \cap S_t^\pm$. On va montrer que

$$\beta_0\left(\dot{a}_0^+\right) > 0 \quad \text{où} \quad \dot{a}_0^+ = \frac{da_t^+}{dt}\bigg|_{t=0}. \tag{$\ddagger$}$$

L'inégalité

$$\beta_0\left(\dot{a}_0^-\right) < 0 \quad \text{où} \quad \dot{a}_0^- = \frac{da_t^-}{dt}\bigg|_{t=0}$$

en résultera formellement par un simple renversement d'orientation sur $\xi$.

On supprime désormais les exposants $+$, pour alléger les notations, mais on pose

$$A_t = \left\{a_s = a_s^+ \mid s \in [0, t]\right\},$$
$$B_t = \left\{b_s = b_s^+ \mid s \in [0, t]\right\},$$

ces arcs étant respectivement orientés de $a_t$ vers $a_0$ et de $b_0$ vers $b_t$. Ainsi,

$$\beta_0\left(\dot{a}_0^+\right) = -\frac{d}{dt}\left(\int_{A_t} \beta_t\right)\bigg|_{t=0}$$

et $S_0 \cup B_t \cup (-S_t) \cup A_t$ borde une chaîne singulière $D_t$. Malheureusement, lorsqu'on l'applique directement à $\beta_t$, la formule de Stokes s'écrit

$$-\int_{A_t} \beta_t = \int_{B_t} \beta_t + \int_{S_0} \beta_t - \int_{D_t} d\beta_t$$

et ne fait pas apparaître ostensiblement la condition de contact. La raison principale est que le terme $\int_{D_t} d\beta_t$ est difficile à évaluer. L'idée est donc de le supprimer en multipliant la forme $\beta_t$ par une fonction positive $h_t$ afin de la rendre fermée. La fonction $h_t$ est construite plus loin ; on montre qu'elle est définie sur un domaine $G_t$ convenable (contenant $S_0$ pour $t = 0$), qu'on peut lui imposer $h_t(a_0) = 1$ et qu'alors

$$\beta_0(\dot{a}_0) = \int_{S_0} h_0 \dot{\beta}_0 \,.$$

Auparavant, on déduit l'inégalité (‡) de cette formule.

Le long de $S_0$, la forme $h_0 \beta_0$ est fermée, c'est-à-dire que

$$d\beta_0 = \beta_0 \wedge d\log h_0,$$

et on peut réécrire la condition de contact

$$u_0 \, d\beta_0 + \beta_0 \wedge (du_0 - \dot{\beta}_0) > 0$$

sous la forme

$$\beta_0 \wedge (u_0 (d\log h_0) + du_0 - \dot{\beta}_0) > 0\,.$$

En multipliant par $h_0$, on obtient

$$\beta_0 \wedge (d(u_0 h_0) - h_0 \dot{\beta}_0) > 0\,.$$

Soit alors $x$ un point courant de $S_0$ (qu'on va faire tendre vers $b_0$) et $S_0(x)$ l'arc orienté de $S_0$ qui va de $a_0$ à $x$. Comme $u_0$ est nulle en $a_0$, l'inégalité ci-dessus entraîne que la fonction

$$f(x) = \int_{S_0(x)} h_0 \dot{\beta}_0 - h_0(x) u_0(x)$$

croît quand $x$ se rapproche de $b_0$ et est ainsi minorée par une constante positive $\varepsilon$ hors d'un voisinage de $a_0$. Par ailleurs, $u_0$ est positive en $b_0$, donc, pour $x$ assez loin sur $S_0$,

$$\int_{S_0(x)} h_0 \dot{\beta}_0 = f(x) + h_0(x) u_0(x) > \varepsilon,$$

et, à la limite,

$$\int_{S_0} h_0 \dot{\beta}_0 \geq \varepsilon > 0\,.$$

Reste à trouver la fonction $h_t$ qui, par définition, doit vérifier

$$d\beta_t = \beta_t \wedge d\log h_t\,.$$

Cette équation n'a pas de solutions aux singularités de $\beta_t$ (par exemple en $b_t$) puisque $d\beta_t$ est non nulle. Elle est en revanche facile à résoudre sur le domaine $G_t$ formé des segments d'orbites de $\xi F_t$ issus de $A$ et arrêtés (si besoin est) avant qu'ils ne reviennent dans $A$. On pose $h_t = 1$ sur $A$ et on prolonge $\log h_t$ à $G_t$ en intégrant sur les orbites de $\xi F_t$, à partir de $A$, la 1-forme $\eta_t$ définie hors des singularités de $\beta_t$ par

$$d\beta_t = \beta_t \wedge \eta_t\,.$$

Évidemment, le domaine $G_t$ ne couvre pas l'arc $B_t$, donc pas non plus la 2-chaîne $D_t$, mais les observations qui suivent permettent de contourner cette difficulté :

- comme $b_t$ est une singularité non dégénérée de $\beta_t$, la norme de $db_t/dt$ est bornée par une constante $\delta$ et l'arc $B_t$ reste ainsi dans le $\delta t$-voisinage $U_t$ de $b_0$ ;
- pour $t$ assez petit, le domaine $G_t$ couvre une bande qui contient à la fois $S'_0 = S_0 \setminus U_t$ et $S'_t = S_t \setminus U_t$ ;
- la fonction $h_t$ est bornée sur $G_t \cap U_t$. En effet, $d\beta_t$ est positive en $b_t$ tandis que $\beta_t$ s'annule ; la relation $d\beta_t = \beta_t \wedge d\log h_t$ montre alors que, lorsque $x$ s'approche de $b_t$ sur $S_t$, la valeur de $d\log h_t(x)$ sur le vecteur unitaire tangent à $S_t$ tend vers $-\infty$ si bien que $h_t(x)$ tend vers $0$.

On substitue alors à $D_t$ une 2-chaîne $D'_t$ contenue dans $G_t$ et bordée par $S'_0 \cup B'_t \cup (-S'_t) \cup A_t$, où $B'_t$ est l'arc de $\partial U_t$ qui joint $S'_0$ à $S'_t$ en restant dans $G_t$. La formule de Stokes, appliquée cette fois à $h_t\beta_t$, s'écrit

$$-\int_{A_t} \beta_t = \int_{B'_t} h_t\beta_t + \int_{S'_0} h_t\beta_t \,.$$

Or les coefficients de $\beta_t$ sur $U_t$, ainsi que la longueur de $B'_t$, sont majorés par $\operatorname{Const} t$ donc

$$\left| \int_{B'_t} \beta_t \right| < \operatorname{Const} t^2$$

et de même, puisque $h_t$ est bornée sur $G_t \cap U_t$,

$$\left| \int_{S_0 \setminus S'_0} \beta_t \right| < \operatorname{Const} t^2 \,.$$

L'égalité

$$\beta_0(\dot{a}_0) = \frac{d}{dt}\left(\int_{S_0} h_t\beta_t\right)\bigg|_{t=0} = \int_{S_0} h_0\dot{\beta}_0$$

en découle immédiatement. $\square$

## G  Lemme d'élimination

On suppose que $\xi F_0$ possède deux singularités non dégénérées *en position d'élimination*, c'est-à-dire un foyer $a_0$ et une selle $b_0$ qui ont le même signe et sont liés par une orbite. Cette configuration stable survit en chaque feuilletage $\xi F_t$, pour $t$ dans un intervalle $[-\varepsilon_0, \varepsilon_0]$, et on note $C_t$ l'adhérence de l'orbite qui joint $a_t$ à $b_t$ sur $F_t$.

Lorsque la surface $F_0$ est convexe, le lemme II.3.3 de [Gi1] montre qu'une isotopie permet d'éliminer les points complexes $a_0$ et $b_0$. On donne ci-dessous une version de ce lemme qui ne suppose plus la convexité et contrôle les modifications apportées aux surfaces voisines.

**Lemme 2.14.** *On se donne un nombre positif $\delta < \varepsilon_0$ et un voisinage $U$ de la nappe $\bigcup_{|t| \le \delta} C_t$ dans lequel les seuls points complexes de $F_t$, pour $t \in [-\varepsilon_0, \varepsilon_0]$, sont $a_t$ et $b_t$. Par ailleurs, toujours pour $t \in [-\varepsilon_0, \varepsilon_0]$, on note $C'_t$ la séparatrice de $b_t$ opposée à $C_t$ et $G_t$ la réunion des segments d'orbites de $\xi F_t$ dont les deux extrémités sont dans $U$.*
**a)** *On peut déformer $\xi$ par une isotopie à support dans $U$ afin que $F_t$ n'ait plus de points complexes dans $U$ pour $t \in [-\delta, \delta]$.*
**b)** *Pour $t \in [-\delta, \delta]$, on choisit continûment dans $U$ un point $b'_t$ sur $C'_t$ et un point $a'_t$ relié à $a_t$ par son orbite sur $\xi F_t$. On peut alors réaliser l'élimination de telle sorte qu'à la fin, pour tout $t \in [-\delta, \delta]$, l'orbite du point $b'_t$ dans $U \cap F_t$ passe arbitrairement près de $a'_t$.*

**c)** *On suppose que, pour tout $t$ dans un compact $\Lambda \subset\, ]-\varepsilon_0, \varepsilon_0[$, le feuilletage $\xi F_t|_{G_t}$ se laisse définir par une 1-forme dont la différentielle ne s'annule pas. On peut alors réaliser l'élimination de telle sorte que, pour $t \in \Lambda$, chaque surface $F_t$ qui est initialement convexe et dont le feuilletage $\xi F_t$ est scindé par une courbe disjointe de $U$ reste convexe.*

*Démonstration.* Quitte à tout bouger par une isotopie horizontale (*i.e.* une isotopie de $V$ respectant chaque $F_t$) et à prendre $U$ plus petit, on s'arrange pour avoir d'une part :
- $U = E \times [-\varepsilon, \varepsilon] \subset F \times [-1,1] = V$, où $E$ est un ouvert de $F$ muni de coordonnées $(x_1, x_2) \in \mathbf{R}^2$ et où $[-\delta, \delta] \cup \Lambda \subset\, ]-\varepsilon, \varepsilon[\, \subset [-\varepsilon_0, \varepsilon_0]$.

On suppose d'autre part que, pour tout $t \in [-\varepsilon, \varepsilon]$, la situation se présente comme suit :
- les points $a_t$ et $b_t$ ont pour coordonnées respectives $(0,0)$ et $(-1,0)$ dans le plan $E_t = E \times \{t\}$ ;
- les points $a'_t$ et $b'_t$ ont pour coordonnées respectives $(1,0)$ et $(-2,0)$ ;
- le feuilletage $\xi F_t$ est transversal à tous les cercles de rayons inférieurs à 1 autour de $a_t$ ;
- hors du disque de rayon $1/2$ autour de $a_t$, les droites d'équations $x_2 = \text{Const}$ dans $E_t \simeq \mathbf{R}^2 \times \{t\}$ sont transversales à $\xi F_t$, sauf l'axe $x_2 = 0$ qui (toujours hors du disque) est réunion de feuilles et contient en particulier $C_t$ ainsi que la composante de $C'_t \cap E_t$ qui porte $b'_t$.

D'autre part, pour fixer les idées, on suppose que les singularités $a_t$ et $b_t$ sont positives et on se donne une équation de $\xi$ sur $V$ qu'on écrit, comme d'habitude, $\beta_t + u_t\, dt = 0$.

**Assertion.** *Il existe sur $V$ une forme de contact $\beta'_t + u'_t\, dt$ qui satisfait aux conditions suivantes :*

*a1) les formes $\beta_t + u_t\, dt$ et $\beta'_t + u'_t\, dt$ sont homotopes, relativement à $V \setminus U$, parmi les formes de contact sur $V$ ;*

*a2) pour $t \in [-\delta, \delta]$, chaque forme $\beta'_t$ est non singulière sur $E_t$ ;*

*b1) pour $t \in [-\varepsilon, \varepsilon]$, chaque différence $\beta'_t - \beta_t$ est un multiple positif ou nul de $dx_2$ le long du segment $[b'_t, a'_t] = [-2,1]$ ;*

*b2) pour $t \in [-\delta, \delta]$, chaque $\beta'_t$ est grande devant $\beta_t$ sur le sous-segment $(-1, 1/2]$ de $[b'_t, a'_t]$ ;*

*c1) pour $t \in [-\varepsilon, \varepsilon]$, chaque $\beta'_t$ s'écrit $\beta_t + u'_t\, dh_t$ où $h_t$ est une fonction $F_t \to \mathbf{R}$ à support dans $E_t$ ;*

*c2) si $t \in \Lambda$, si $F_t$ est convexe et si $\xi F_t$ est scindé par une multi-courbe disjointe de $E_t$, il existe une fonction $v'_t : F_t \to \mathbf{R}$ qui coïncide avec $u'_t$ sur le support de $dh_t$ et vérifie la condition de convexité $v'_t\, d\beta_t + \beta_t \wedge dv'_t > 0$.*

Le lemme découle immédiatement de cette assertion.

**a)** C'est une conséquence directe des conditions a1) et a2).

**b)** Les conditions b1) et b2) entraînent que, pour $t \in [-\delta, \delta]$, l'orbite de $b'_t$ sur le feuilletage défini dans $E_t$ par $\beta'_t$ est proche du segment $[b'_t, a'_t]$.

**c)** Les conditions c1) et c2) montrent que, si $t \in \Lambda$, si $F_t$ est convexe et si $\xi F_t$ est scindé par une multi-courbe qui évite $E_t$, alors

$$v'_t\, d\beta'_t + \beta'_t \wedge dv'_t = v'_t\, d\beta_t + \beta_t \wedge dv'_t > 0$$

car $v'_t$ coïncide avec $u'_t$ sur le support de $dh_t$. □

Avant de prouver l'assertion, on pose encore quelques notations. Pour tout réel positif $r$, on désigne par

$$\begin{aligned} D(r) &\quad \text{le disque de rayon } r \text{ autour de } (0,0), \\ R(r) &\quad \text{le rectangle } [-1-r, 1/2+r] \times [-r, r], \\ Q(r) &\quad \text{le carré } [-1-r, -1+r] \times [-r, r]. \end{aligned}$$

*Preuve de l'assertion.* Il s'agit de construire les fonctions $u'_t$ et $h_t$. On choisit d'abord :
- une fonction $\rho \colon [-\varepsilon, \varepsilon] \to [0,1]$ qui est nulle près du bord et vaut 1 sur $[-\delta, \delta] \cup \Lambda$ ;
- une fonction $\sigma_r \colon E = \mathbf{R}^2 \to [0,1]$ qui, pour tout $r \in \,]0, 1/4[$, est nulle hors de $D(1-r) \cup R(2r)$, vaut 1 sur $D(1-2r) \cup R(r)$ et a une différentielle bornée par $\mathrm{Const}\, r^{-1}$.

On se donne aussi une famille de fonctions $v_t \colon G_t \to \,]0, \infty[$, $t \in [-\varepsilon, \varepsilon]$, vérifiant

$$\begin{aligned} d\Bigl(\frac{\beta_t}{v_t}\Bigr) &> 0 &&\text{pour tout } t \in \Lambda, \\ v_t(b_t) &= u_t(b_t) &&\text{pour tout } t \in [-\varepsilon, \varepsilon]. \end{aligned}$$

On pose alors

$$h_t = \rho(t)\, h \quad \text{et} \quad u'_t = u_t + \rho(t)\bigl(\sigma_r\,(v_t - u_t) + w\bigr).$$

Les nouvelles inconnues $h$ et $w$ sont :
- pour $h$, une fonction $\mathbf{R}^2 \to \mathbf{R}$ à support dans le rectangle $R(r)$, nulle sur l'axe réel et égale, sur le rectangle $R(r/2)$, à une fonction du type $(x_1, x_2) \mapsto c x_2$ où $c$ est un nombre positif ;
- pour $w$, une fonction radiale $\mathbf{R}^2 \to [0, \infty[$ de support un disque $D(1-s)$, $s < r$, et logarithmiquement décroissante (à l'intérieur de $D(1-s)$) le long des orbites de $\xi F_t$.

Reste à voir que les conditions énumérées dans l'assertion sont remplies lorsque $r$ et $h$ sont petits et que $\log w$ décroît fortement. Dans les calculs qui suivent, on pose

$$w_t = \sigma_r\,(v_t - u_t) + w, \quad \text{donc} \quad u'_t = u_t + \rho(t)\, w_t.$$

Tout d'abord, la condition de contact pour la 1-forme $\beta'_t + u'_t\, dt$, où $\beta'_t = \beta_t + u'_t\, dh_t$, est que la 2-forme suivante soit positive sur $F$ pour tout $t \in [-\varepsilon, \varepsilon]$ :

$$\begin{aligned} u'_t\, d\beta'_t &+ \beta'_t \wedge \bigl(du'_t - \dot\beta'_t\bigr) \\ &= u_t\, d\beta_t + \beta_t \wedge \bigl(du_t - \dot\beta_t\bigr) + dh \wedge \bigl(\rho(t)\dot u_t \beta_t + \rho'(t) u_t \beta_t - \rho(t) u_t \dot\beta_t\bigr) \\ &\quad + \rho(t)\Bigl(w_t\, d\beta_t + \beta_t \wedge dw_t + w_t\, dh \wedge \bigl(2\rho'(t)\beta_t - \rho(t)\dot\beta_t\bigr)\Bigr). \end{aligned}$$

On réécrit plus synthétiquement cette expression

$$\mu_t + dh \wedge \nu_t + \rho(t)\bigl(w_t\, d\beta_t + \beta_t \wedge dw_t + w_t\, dh \wedge \lambda_t\bigr)$$

où les formes $\lambda_t$, $\mu_t$ et $\nu_t$ ne dépendent ni de $h$ ni de $w$.

Comme $|d\sigma_r| \leq \mathrm{Const}\, r^{-1}$, les formes $\beta_t \wedge d\sigma_r$ sont bornées indépendamment de $r$ au voisinage de $a_t$ et $b_t$. On choisit alors $r$ assez petit pour que les 2-formes

$$d\beta_t \quad \text{et} \quad \tfrac{1}{3}\mu_t + (v_t - u_t)\,\beta_t \wedge d\sigma_r + \sigma_r\bigl((v_t - u_t)\, d\beta_t + \beta_t \wedge d(v_t - u_t)\bigr)$$

soient positives respectivement sur $D(2r) \cup Q(2r)$ et sur $Q(2r)$ (en $b_t$, $\beta_t$ et $v_t - u_t = w_t - w$ s'annulent). On prend ensuite pour $h$ une fonction $R(r) \to \mathbf{R}$ du type indiqué plus haut et suffisamment $\mathcal{C}^1$-petite pour que les formes

$$\tfrac{1}{3}\mu_t + (w_t - w)\, dh \wedge \lambda_t, \quad d\beta_t + dh \wedge \lambda_t \quad \text{et} \quad \tfrac{1}{3}\mu_t + dh \wedge \nu_t$$

soient positives respectivement sur $Q(r)$, sur $D(r)$ et sur $F_t$ tout entière. On prend enfin $w$ grande sur $D(r)$ et logarithmiquement très décroissante le long des orbites de $\xi F_t$ sur $D(1-s) \setminus D(r)$. Ainsi,

$$w\left(\beta_t \wedge \frac{dw}{w} + d\beta_t + dh \wedge \lambda_t\right)$$

est positive sur $\operatorname{Int} D(1-s)$ (et arbitrairement grande sur tout compact de ce disque).

Ce qui précède montre en fait que la condition a1) est vérifiée : les segments joignant respectivement $\beta_t + u_t\, dt$ à $\beta_t + u'_t\, dt$ et $\beta_t + u'_t\, dt$ à $\beta'_t + u'_t\, dt$ donnent le chemin voulu parmi les formes de contact.

Pour voir que la condition a2) est satisfaite lorsque $h$ est assez petite et $w$ assez grande ($r$ étant fixé), on observe que :
– sur la région $R(r) \setminus R(r/2)$, la forme $\beta_t$ n'approche pas de $0$ ;
– sur $R(r/2) \setminus D(1/2)$, les niveaux de $h$ sont transversaux aux orbites de $\xi F_t$ ;
– sur $D(1/2)$, la fonction $w$ est arbitrairement grande.

Ces remarques permettent aussi de réaliser b1) et b2). Enfin, la condition c1) est vérifiée par construction-même et c2) l'est pour la raison suivante : si le feuilletage $\xi F_t$ est scindé par une multi-courbe disjointe de $E_t$, il est expansif sur la région $G_t$. Comme toute orbite de $\xi F_t$ qui quitte $G_t$ n'y revient plus, la fonction $v'_t = v_t + w$, qui coïncide avec $u'_t$ sur le support de $h_t$, se prolonge à $F_t$ en une fonction qui vérifie la condition de convexité. □

**Remarque 2.15.** Le procédé décrit ci-dessus permet aussi de supprimer d'un coup, pour tout $t \in [-\delta, \delta]$, plusieurs paires de singularités $(a^i_t, b^i_t)$ qui sont en position d'élimination dans $\xi F_t$ à chaque instant $t \in [-\varepsilon_0, \varepsilon_0]$. Le support de l'isotopie à effectuer est un voisinage $U$ arbitrairement petit de la réunion des arcs d'élimination $C^i_t$. De plus, les surfaces convexes $F_t$, $t \in \Lambda$, dont le feuilletage est scindé par une courbe disjointe de $U$ restent convexes à la condition suivante : pour tout $t \in \Lambda$, le feuilletage $\xi F_t |_{G_t}$ est défini par une forme dont la différentielle ne s'annule pas, où $G_t$ désigne cette fois la région formée des segments d'orbites dont les deux extrémités sont dans des composantes de $U$ contenant des singularités d'un même signe.

## H  Lemme de préparation

On suppose maintenant que $\xi$ est une structure de contact *tendue* sur $V$, c'est-à-dire qu'aucun disque plongé dans $V$ n'est tangent à $\xi$ en tous les points de son bord. Le lemme ci-dessous permet de déformer $\xi$ pour que les surfaces $F_t$ non convexes aient le moins de points complexes possible.

**Lemme 2.16.** *Sur $V$, toute structure de contact tendue et à bord convexe se laisse déformer, par une isotopie relative au bord, en une structure pour laquelle les surfaces $F_t$ non convexes ont toutes exactement $|\chi(F)|$ points complexes.*

**Remarque 2.17.** Si $F$ est la sphère $\mathbf{S}^2$, le lemme ci-dessus donne une isotopie au terme de laquelle toutes les sphères $F_t$ sont convexes. Les lemmes 2.6 et 2.8 permettent alors de redémontrer facilement les résultats de Y. ELIASHBERG selon lesquels $\mathbf{S}^2 \times \mathbf{R}$, $\mathbf{S}^2 \times \mathbf{S}^1$, $\mathbf{S}^3$ et $\mathbf{R}^3$ portent une unique structure de contact tendue.

*Démonstration.* On déforme d'abord la structure par une isotopie qui banalise la dynamique (lemme 2.9) tout en réalisant un film générique (remarque 2.10). On considère ensuite l'ensemble $\Sigma$ des instants $t$ dans $[-1, 1]$ où la surface $F_t$ n'est pas convexe. C'est un fermé d'intérieur vide (propriétés 1) qui, si $\xi$ désigne la structure, est la réunion disjointe des ensembles suivants (proposition 2.4) :
  – l'ensemble $\Sigma_1$ des $t$ pour lesquels $\xi F_t$ présente une orbite fermée dégénérée — celle-ci est positive ou négative (cf. E), d'où une sous-partition de $\Sigma_1$ en $\Sigma_1^+$ et $\Sigma_1^-$ ;
  – l'ensemble $\Sigma_2$ des $t$ pour lesquels $\xi F_t$ présente une connexion de selles rétrograde (cf. F).
On montre ci-dessous comment appliquer le lemme d'élimination 2.14 près des surfaces $F_t$, $t \in \Sigma$, pour obtenir le résultat voulu.

Si $t_0 \in \Sigma_1^+$, le feuilletage $\xi F_{t_0}$ n'a que des singularités non dégénérées et aucune connexion de selles mais présente une orbite fermée dégénérée, laquelle indique la naissance de deux orbites fermées non dégénérées (lemme 2.11). Par suite, dans un intervalle $]t_0, t_0+\varepsilon_0]$, les surfaces $F_t$ sont toutes convexes et les feuilletages $\xi F_t$ sont structurellement stables (de type Morse-Smale). L'étude des surfaces convexes et de leurs feuilletages dans les structures tendues [El4, Gi2] montre alors que, si $\xi F_t$ compte $2n + |\chi(F)|$ singularités pour $t$ dans $]t_0, t_0 + \varepsilon_0]$, on peut grouper $2n$ d'entre elles par paires disjointes en position d'élimination (cf. G). De plus, chaque paire, avec son arc d'élimination, se laisse suivre par continuité quand $t$ passe en-deçà de $t_0$, jusqu'en $t_0 - \varepsilon_1$. Par ailleurs, on peut trouver des réels positifs $\delta < \varepsilon \leq \min(\varepsilon_0, \varepsilon_1)$ tels que, dans l'intervalle $[t_0 - \varepsilon, t_0 - \delta]$, les surfaces $F_t$ soient toutes convexes et portent des feuilletages $\xi F_t$ structurellement stables. Le lemme 2.14 (avec la remarque 2.15) fournit alors une isotopie à support dans $F \times [t_0 - \varepsilon, t_0 + \varepsilon]$ au terme de laquelle les surfaces $F_t$ non convexes, pour $t \in [t_0 - \varepsilon, t_0 + \varepsilon]$, ont exactement $|\chi(F)|$ points complexes.

Si $t_0 \in \Sigma_2$, le feuilletage $\xi F_{t_0}$ a des singularités et des orbites fermées toutes non dégénérées mais présente une et une seule connexion de selles, laquelle est rétrograde et indique le croisement de deux séparatrices (lemme 2.13). Les autres séparatrices de $\xi F_{t_0}$ (celles qui ne participent pas à la connexion) sont stablement accrochées à des foyers ou à des orbites fermées. Par suite, dans un intervalle entrecoupé $[t_0 - \varepsilon, t_0 + \varepsilon] \setminus \{t_0\}$, les surfaces $F_t$ sont toutes convexes et les feuilletages $\xi F_t$ sont structurellement stables.

Afin de décrire précisément le processus d'élimination, on introduit les objets et les notations suivantes pour tout $t \in [t_0 - \varepsilon, t_0 + \varepsilon]$ :
  – $\Gamma_t^+$ (resp. $\Gamma_t^-$) est le « graphe » formé, sur $F_t$, par les orbites fermées répulsives (resp. attractives) de $\xi F_t$, les singularités positives (resp. négatives) et leurs variétés stables (resp. instables) ;
  – $b_t^+$ et $b_t^-$ sont les selles de $\xi F_t$, respectivement positive et négative, qui sont connectées à l'instant $t = t_0$ ;
  – $S_t^+$ et $S_t^-$ sont les séparatrices respectives de $b_t^+$ et $b_t^-$ qui se confondent à l'instant $t = t_0$ ;
  – $\underline{a}_t^{\pm}$ (resp. $\overline{a}_t^{\pm}$) est le foyer ou l'orbite fermée de $\xi F_t$ — variant continûment avec $t$ — sur lequel, quand $t < t_0$ (resp. quand $t > t_0$), s'accumule le bout de $S_t^{\pm}$ opposé à $b_t^{\pm}$.
Pour $t \neq t_0$, les graphes $\Gamma_t^+$ et $\Gamma_t^-$ sont fermés, disjoints, et l'inégalité de Bennequin [El4, Gi2] affirme qu'aucune de leurs composantes connexes n'est un arbre, sauf si $F$ est la sphère $\mathbf{S}^2$, auquel cas $\Gamma_t^+$ et $\Gamma_t^-$ sont deux arbres connexes. Les graphes $\Gamma_{t_0}^{\pm}$, en revanche, ont en commun la séparatrice $S = S_{t_0}^{\pm}$ et ne sont pas fermés puisque $b_{t_0}^{\mp} \notin \Gamma_{t_0}^{\pm}$. Dans la suite, on omet l'indice $t$ dans les notations ci-dessus quand $t = t_0$ et on note $A^{\pm}$ (resp. $B^{\pm}$) la composante connexe de $\Gamma^\pm$ qui contient $\underline{a}^\pm$ (resp. $b^\pm$). L'assertion ci-dessous donne la clé pour mener l'élimination à son terme :

**Assertion.** *Si $F \neq \mathbf{S}^2$ et si $A^+$ (resp. $A^-$) est un arbre fermé, alors $B^-$ (resp. $B^+$) est un arbre — non fermé.*

*Preuve.* On suppose que $A^+$ est un arbre fermé et que $B^-$ contient des cycles. On observe alors que :
- comme $A^+$ est un arbre, $\underline{a}^+$ est un foyer ;
- comme $A^+$ est fermé, il ne contient pas $S$ et diffère donc de $B^+$ ;
- comme, pour $t \neq t_0$, aucune composante connexe de $\Gamma_t^+$ n'est un arbre, $B^+$ contient au moins un cycle et $A^+$ contient $\overline{a}^+$, lequel est donc aussi un foyer.

Ainsi, les deux bouts de la variété stable $R$ de $b^-$ s'appuient sur $A^+$ et le cycle $A^+ \cup R$ n'est pas contractile, sans quoi il enfermerait dans un disque l'un des cycles de $B^+$ ou de $B^-$ (ce qu'interdit l'inégalité de Bennequin). On explique ci-dessous comment bouger $F_{t_0}$ par une isotopie pour faire apparaître dans son feuilletage une nouvelle configuration interdite.

On choisit sur $B^-$ un cycle $C$ et un arc simple $B$ qui joint $b^-$ à $C$. On plie ensuite $F_{t_0}$ deux fois (en sens inverses) le long de $C$ du côté opposé à $B$. Cette opération change les graphes $\Gamma^\pm$ comme suit : le cycle $C$ se dédouble en deux cycles négatifs parallèles, $C$ et $C'''$, qui sont séparés par un cycle positif isolé $C'$ ; l'arc $B$ joint toujours $b^-$ à $C$. Maintenant, si $C$ est une orbite fermée, on y crée une paire de singularités négatives en position d'élimination *cyclique*, au sens où toute la variété stable de la selle aboutit au foyer et couvre $C$. L'extrémité de $B$ située sur $C$ est alors un foyer $a^-$. L'ultime opération consiste à éliminer toutes les singularités qui se trouvent sur $B$ (y compris les extrémités $a^-$ et $b^-$). Le lemme 2.14-b permet de s'arranger pour qu'à la fin, les séparatrices qui aboutissaient en $a^-$ sur $C$ aillent l'une vers $\underline{a}^-$ et l'autre vers $\overline{a}^-$. Dans ces conditions, la séparatrice $S$ de $b^+$ qui venait de $b^-$ vient du cycle positif $C'$ et l'arbre $A^+$ se trouve isolé des autres parties du graphe positif. Les techniques de [El4, Gi2] s'appliquent alors pour montrer que la structure n'est pas tendue. $\square$

On aborde maintenant l'élimination. Si $2n_0$ des $2n + \chi(F)$ singularités de $\xi F_{t_0}$ forment $n_0$ paires disjointes en position d'élimination et ne contenant pas les selles $b^\pm$, le lemme 2.14 et la remarque 2.15 donnent une isotopie à support dans $F \times [t_0 - \varepsilon, t_0 + \varepsilon]$ au terme de laquelle les surfaces $F_t$ non convexes, pour $t \in [t_0 - \varepsilon, t_0 + \varepsilon]$, ont exactement $2(n - n_0) + |\chi(F)|$ points complexes. Le nombre maximal $n_0$ de ces paires « inoffensives » vaut
- $n$ si $B^+$ et $B^-$ contiennent tous deux des cycles (l'assertion assure qu'alors, ni $A^+$ ni $A^-$ ne sont des arbres) ;
- $n - 2$ si $B^+$ et $B^-$ sont tous deux des arbres ;
- $n - 1$ dans tout autre cas.

Après cette élimination, on choisit un nombre positif $\delta < \varepsilon$ satisfaisant aux conditions suivantes :
- $F_t$ a exactement $2(n - n_0) + \chi(F)$ points complexes si $|t - t_0| < \delta$ ;
- $F_t$ est convexe si $\delta \leq |t - t_0| \leq \varepsilon$.

Pour fixer les idées, on suppose désormais que $B^-$ est un arbre et on raisonne cas par cas.

**1)** Si $B^+$ est aussi un arbre, l'élimination des $n - 2$ paires inoffensives laisse sur $F_t$, pour tout $t \in [t_0 - \delta, t_0 + \delta]$, deux points elliptiques $a_t^+$ et $a_t^-$ qui sont en position d'élimination avec $b_t^+$ et $b_t^-$ respectivement. On se donne alors un nombre positif $\gamma < \delta$ assez proche de $\delta$ pour que $F_t$ soit convexe si $\gamma \leq |t - t_0| \leq \delta$. On applique ensuite le lemme 2.14 simultanément aux deux paires $(a^\pm, b^\pm)$ en prenant pour $\Lambda$ la réunion des intervalles $[t_0 - \delta, t_0 - \gamma]$ et $[t_0 + \gamma, t_0 + \delta]$.

**2)** Si ni $B^+$ ni $A^+$ ne sont des arbres, on élimine comme en 1 le seul point elliptique qui reste après la suppression des $n-1$ paires inoffensives.

**3)** Si $A^+$ est un arbre (fermé), toute la variété stable de $b^-$ vient de $A^+$ et forme avec $A^+$ un cycle contractile qui enferme dans un disque l'arbre $B^- \setminus S$ (voir la preuve de l'assertion). Après avoir supprimé les $n-1$ paires inoffensives, on élimine $b^-$, comme en 2, avec l'unique point elliptique $a^-$ qui reste sur $B^-$. Ceci fait, la séparatrice $S$ de $b^+$ qui venait de $b^-$ vient du point elliptique $a^+$ sur lequel s'est effondré $A^+$ et il ne reste plus alors qu'à éliminer la paire $(a^+, b^+)$. □

## I  Lemmes de bifurcation et théorème de Bennequin

La preuve du lemme 2.16 invoque abusivement les lemmes 2.11 et 2.13. En effet, l'important n'est pas de connaître le sens des bifurcations mais simplement de savoir qu'elles ont lieu. Or cela fait partie du comportement générique du film. Pour clore cette partie, on esquisse une démonstration du théorème de D. BENNEQUIN [Be] dans laquelle les deux lemmes de bifurcation 2.11 et 2.13 jouent un rôle crucial. Tout d'abord, on observe que la structure de contact ordinaire sur $\mathbf{R}^3$, d'équation $dz + x\, dy - y\, dx = 0$, est invariante par le flot des transformations $(x, y, z) \mapsto (e^t x, e^t y, e^{2t} z)$. Toutes les sphères euclidiennes centrées à l'origine sont donc convexes et leur feuilletage caractéristique est scindé par une courbe connexe. Le théorème de BENNEQUIN — selon lequel cette structure de contact est tendue — est donc équivalent au résultat suivant :

**Théorème 2.18** [Be]. *Soit $\xi$ une structure de contact sur $V = \mathbf{S}^2 \times [-1, 1]$ pour laquelle toutes les sphères $S_t = \mathbf{S}^2 \times \{t\}$, $t \in [0, 1]$, sont convexes. Si chaque feuilletage $\xi S_t$ est scindé par une courbe connexe, la structure $\xi$ est tendue.*

*Esquisse de démonstration.* On suppose que $\xi$ est vrillée. Il existe alors une $\partial$-isotopie $\phi_s$ de $V$, $s \in [0, 1]$, qui amène $S_0$ sur une sphère convexe $\phi_1(S_0)$ contenant un disque partout tangent à $\xi$ le long de son bord. Le feuilletage $\xi\phi_1(S_0)$ est alors scindé par une courbe non connexe. Dans $[0, 1] \times [-1, 1]$, on note $\Omega$ l'ouvert dense des points $(s, t)$ pour lesquels $\phi_s(S'_t)$ est une sphère convexe. Cet ouvert est l'union des ouverts $\Omega_t$ et $\Omega_v$ où $\Omega_t$ (resp. $\Omega_v$) est formé des points $(s, t)$ pour lesquels $\xi\phi_s(S_t)$ est scindé par une courbe connexe (resp. non connexe). On pose alors $s_0 = \inf \pi(\Omega_v)$ où $\pi$ est la projection $[0, 1] \times [-1, 1] \to [0, 1]$ et on choisit, au-dessus de $s_0$, un point $(s_0, t_0)$ dans l'adhérence de $\Omega_v$.

Par hypothèse, $s_0$ est strictement positif et, par construction, $(s_0, t_0)$ appartient au complémentaire $\Sigma$ de $\Omega$. L'ensemble $\Sigma$ est un fermé naturellement stratifié en strates de codimensions 1 et 2. Comme tout feuilletage de $\mathbf{S}^2$ vérifie la propriété de Poincaré-Bendixson par le théorème du même nom, les strates de codimension 1 correspondent aux feuilletages $\xi\phi_s(S_t)$ qui présentent soit une orbite fermée dégénérée de naissance ou de mort, soit une connexion rétrograde entre deux selles. On note $\Sigma^1$ la réunion de ces strates. Les strates de codimension 2 sont de deux types :
  – les points d'intersection entre deux strates de codimension 1 — on note $\Sigma^{11}$ leur ensemble ;
  – les points $(s, t)$ correspondant aux feuilletages $\xi\phi_s(S_t)$ qui présentent soit une orbite fermée plus dégénérée (mais stable), soit une connexion rétrograde entre une selle et un nœud-selle — on note $\Sigma^2$ l'ensemble de ces points.

La première observation clé est que le point $(s_0, t_0)$ n'est pas dans $\Sigma^1$. En effet, les lemmes 2.11 et 2.13 assurent que la restriction de $\pi$ à ces strates est toujours un homéomorphisme local. Un argument facile montre que $(s_0, t_0)$ n'est pas non plus dans $\Sigma^2$ : un point $(s, t) \in \Sigma^2$ adhère à une seule strate de codimension 1 et les composantes de $\Omega$ voisines de $(s, t)$ se projettent sur tout un voisinage de $s$. La dernière possibilité est donc que $(s_0, t_0)$ soit à l'intersection (génériquement transversale) de deux strates $\Sigma_0$ et $\Sigma_1$ de codimension 1. Aucune de ces strates ne correspond à une orbite fermée puisque sa projection sur $[0, 1]$ se prolonge en-deçà de $s_0$. Ainsi, $\Sigma_0$ et $\Sigma_1$ correspondent à des connexions de selles. Par ailleurs, elles partagent localement $\Omega$ en quatre régions : $G$ (gauche), $D$ (droite), $H$ (haut) et $B$ (bas). Par hypothèse, $D$ est dans $\Omega_v$ tandis que $B$ et $H$ sont dans $\Omega_t$. On laisse au lecteur le plaisir de découvrir comment le lemme 2.13 implique alors que $G$ est dans $\Omega_v$. $\square$

# 3  Formes normales des structures de contact tendues

Dans cette partie, on cherche des formes normales pour les structures de contact tendues sur les variétés fibrées en tores au-dessus de l'intervalle ou du cercle. Ces formes normales sont données par des films simples sur le tore. On désigne donc désormais par $F$ un tore orienté de dimension deux.

## A  Structures rotatives

Il y a deux familles typiques de structures de contact tendues sur $F \times [0, 1]$ : les structures pour lesquelles tous les tores $F_t$ sont convexes et les structures rotatives.

**Définition 3.1.** Soit $V$ une variété fibrée en tores au-dessus d'un intervalle ou d'un cercle. Une structure de contact sur $V$ est une *structure rotative* si le feuilletage caractéristique de chaque fibre est une suspension.

Les cycles asymptotiques [Sc] d'une suspension non orientée $\sigma$ sur $F \simeq \mathbf{T}^2$ engendrent, dans le plan $H_1(F; \mathbf{R})$, une droite vectorielle qu'on note $D(\sigma)$. Si $I$ est un intervalle et $\xi$ une structure rotative sur $F \times I$, la condition de contact assure que la fonction

$$I \longrightarrow \mathbf{P}\bigl(H_1(F; \mathbf{R})\bigr) \simeq \mathbf{R}/\pi\mathbf{Z} \quad t \longmapsto D_t = D(\xi F_t),$$

est une fonction monotone, décroissante au sens large si on oriente le plan $H_1(F; \mathbf{R})$ au moyen de la forme d'intersection. On appelle *amplitude* de $\xi$ la variation totale de cette fonction.

Les structures rotatives d'amplitude non nulle ont une classification simple qui sera établie plus loin et qui repose sur le fait suivant :

**Lemme de redressement 3.2.** *Soit $\xi$ une structure rotative sur $V = F \times [0, 1]$ dont l'amplitude est non nulle et strictement supérieure à $-\pi$. Si les feuilletages $\xi F_0$ et $\xi F_1$ ont des projections sur $F$ transversales l'une à l'autre, $\xi$ est $\partial$-isotope, parmi les structures rotatives, à une structure de contact qui contient le champ de vecteurs $\partial_t$.*

*Démonstration.* Il s'agit de trouver dans $(V, \xi)$ un feuilletage legendrien $\partial$-isotope au feuilletage par les segments $\{*\} \times [0, 1]$.

L'amplitude de $\xi$ étant supérieure à $-\pi$, on peut feuilleter $V$ par des anneaux $A_s$, $s \in \mathbf{R}/\mathbf{Z}$, qui coupent transversalement chaque tore $F_t$ selon une courbe $C_{st} = A_s \cap F_t$ elle-même transversale à $\xi F_t$. La réunion des divers feuilletages $\xi A_s$ forme alors un feuilletage legendrien de $(V, \xi)$ transversal aux tores $F_t$ et dont on veut corriger l'holonomie $F_0 \to F_1$ en utilisant la différence de pentes entre $\xi F_0$ et $\xi F_1$.

Pour contrôler la propriété d'isotopie requise, on passe dans l'espace $\widetilde{V} = \widetilde{F} \times [0, 1]$ où $\widetilde{F}$ est le revêtement universel de $F$. Les relèvements respectifs de $\xi$, $F_t$, $A_s$ et $C_{st}$ sont notés $\widetilde{\xi}$, $\widetilde{F}_t$, $\widetilde{A}_s$ et $\widetilde{C}_{st}$, $s \in \mathbf{R}$. On désigne par $\lambda$ le feuilletage legendrien formé des feuilletages $\widetilde{\xi}\widetilde{A}_s$, $s \in \mathbf{R}$, on note $\phi$ son holonomie $\widetilde{F}_0 \to \widetilde{F}_1$ et, pour $i \in \{0, 1\}$, on pose $\tau_i = \widetilde{\xi}\widetilde{F}_i$.

La condition de contact et le fait que $\phi$ envoie chaque courbe $\widetilde{C}_{s0}$ sur $\widetilde{C}_{s1}$ assurent que les feuilletages $\tau_0$ et $\phi^*\tau_1$ sont transversaux l'un à l'autre. La non-nullité de l'amplitude garantit de plus que chaque feuille de $\tau_0$ coupe chaque feuille de $\phi^*\tau_1$ en un point (forcément unique). Pour tout $(x, 0) \in \widetilde{F}_0$, on note alors $\psi(x, 0)$ le point d'intersection des feuilles de $\tau_0$ et $\phi^*\tau_1$ passant respectivement par $(x, 0)$ et $\phi^{-1}(x, 1)$. l'application $\psi$ ainsi définie est un difféomorphisme de $\widetilde{F}_0$. En effet, comme $\tau_0$ et $\tau_1$ ont des projections sur $\widetilde{F}$ transversales l'une à l'autre, $\psi$ est injective : si deux points $(x, 0)$ et $(x', 0)$ appartiennent à une même feuille de $\tau_0$, les points $\phi^{-1}(x, 1)$ et $\phi^{-1}(x', 1)$ sont sur des feuilles distinctes de $\phi^*\tau_1$.

Le difféomorphisme $\psi$ offre une méthode pour relier tout point $(x, 0) \in \widetilde{F}_0$ au point $(x, 1) \in \widetilde{F}_1$ par un arc legendrien : on suit la feuille de $\tau_0$ entre $(x, 0)$ et $\psi(x, 0)$, puis la feuille de $\lambda$ jusqu'en $\phi \circ \psi(x, 0)$ et enfin la feuille de $\tau_1$ jusqu'en $(x, 1)$.

Ces arcs sont isotopes aux segments verticaux relativement à leurs bouts mais ne sont pas lisses et ne feuillettent pas $\widetilde{V}$ — ils se chevauchent sur le bord. Ce sont en fait les caractéristiques de bandes anguleuses obtenues en collant aux bandes $\widetilde{A}_s$ des morceaux de $\partial\widetilde{V}$ et bordées par des courbes qui feuillettent $\partial\widetilde{V}$. De plus, comme la construction est invariante par l'action de $\pi_1(F)$, ces bandes se projettent sur $V$ en des anneaux. Il est facile de lisser et de bouger légèrement ces anneaux, sans toucher leurs bords, pour obtenir un feuilletage de $V$ par des anneaux $A'_s$ transversaux aux tores $F_t$. Si la perturbation est assez petite, l'holonomie du feuilletage legendrien $\bigcup_s \xi A'_s$ est $\mathcal{C}^1$-proche de l'« identité », i.e. de l'application $(x, 0) \in F_0 \mapsto (x, 1) \in F_1$. Il existe donc sur $V$ une structure de contact $\xi'$ qui est $\mathcal{C}^1$-proche de $\xi$, coïncide avec $\xi$ près du bord et pour laquelle l'holonomie du feuilletage legendrien $\bigcup \xi' A'_s$ est l'identité. Vu la stabilité des structures de contact, $\xi$ est $\partial$-isotope à $\xi'$. □

## B Classification des structures rotatives

**Théorème 3.3.** *Sur $V = F \times [0, 1]$, deux structures rotatives qui coïncident au bord et qui ont des amplitudes égales et non nulles sont $\partial$-isotopes.*

Avant de démontrer ce théorème, il est utile de noter que les structures rotatives d'amplitude nulle se comportent autrement.

**Exemple 3.4.** Toute structure rotative sur $F \times [0, 1]$ pour laquelle chaque tore $F_t$ est convexe a une amplitude nulle. Une telle structure $\xi_0$ étant choisie, on paramètre $F$ par $\mathbf{T}^2$ de telle sorte que le découpage des tores $F_t$ soit parallèle aux courbes $\{*\} \times \mathbf{S}^1 \times \{t\}$ et on se donne un difféomorphisme $\phi$ de $F \times [0, 1]$ qui induit un twist de Dehn sur chaque anneau $\mathbf{S}^1 \times \{*\} \times [0, 1]$ (et coïncide donc avec l'identité près du bord).

**Assertion.** *Les structures de contact $\xi_0$ et $\xi_1 = \phi_*\xi_0$ ne sont pas $\partial$-isotopes.*

Le lemme 2.6 montre qu'il suffit, pour établir ce fait, de traiter le cas où $\xi_0$ a pour

équation

$$\cos(nx_1)\,dx_2 - \sin(nx_1)\,dt = 0\;.$$

En recollant $F \times \{1\}$ sur $F \times \{0\}$ par l'identité, on fabrique un tore de dimension trois sur lequel $\xi_0$ et $\phi$ induisent respectivement une structure de contact $\overline{\xi}_0$ et un difféomorphisme $\overline{\phi}$. Si $\xi_0$ et $\xi_1$ étaient $\partial$-isotopes, les structures $\overline{\xi}_0$ et $\overline{\xi}_1 = \overline{\phi}_*\overline{\xi}_0$ seraient isotopes. Or il n'en est rien car leurs tores prélagrangiens incompressibles ne sont pas isotopes [Gi4, lemme 10]. □

Pour établir le théorème 3.3, on se ramène au cas de deux structures ayant une amplitude supérieure à $-\pi$.

**Lemme 3.5.** *Soit $\xi_0$ et $\xi_1$ deux structures rotatives sur $V = F \times [0,1]$ qui coïncident au bord et ont des amplitudes égales et non nulles. Quitte à déformer $\xi_0$ et $\xi_1$ parmi les structures rotatives et relativement au bord, on peut subdiviser $[0,1]$ en sous-intervalles $[t_{i-1}, t_i]$, $1 \leq i \leq p$, de telle sorte que les restrictions de $\xi_0$ et $\xi_1$ à chaque produit $F \times [t_{i-1}, t_i]$ coïncident au bord et aient des amplitudes égales, non nulles et strictement supérieures à $-\pi$.*

*Démonstration.* Quitte à reparamétrer $[0,1]$ dans $(V, \xi_1)$ par exemple, on prend une subdivision $0 = t_0 < t_1 < \cdots < t_p = 1$ vérifiant les conditions suivantes :
  – pour $1 \leq i \leq p-1$, les feuilletages $\xi_0 F_{t_i}$ et $\xi_1 F_{t_i}$ ont la même direction asymptotique qui est une droite irrationnelle ;
  – pour $1 \leq i \leq p$, les restrictions de $\xi_0$ et $\xi_1$ à $F \times [t_{i-1}, t_i]$ ont des amplitudes égales, non nulles et strictement supérieures à $-\pi$.

On perturbe alors $\xi_0$ et $\xi_1$ par des isotopies $\mathcal{C}^1$-petites pour rendre convexe chaque tore $F_{t_i}$, $1 \leq i \leq p-1$, à la fois dans $\xi_0$ et $\xi_1$, tout en s'assurant que les feuilletages $\xi_0 F_{t_i}$ et $\xi_1 F_{t_i}$ gardent une direction asymptotique commune (maintenant rationnelle) et possèdent exactement deux orbites fermées. Les lemmes 2.3 et 2.8 fournissent les ultimes isotopies qui font coïncider $\xi_0 F_{t_i}$ et $\xi_1 F_{t_i}$. □

**Remarque 3.6.** Le théorème de M. HERMAN sur la conjugaison des difféomorphismes du cercle [He] assure en fait qu'il suffit de reparamétrer $[0,1]$ (dans $(V, \xi_1)$ par exemple) pour trouver la subdivision voulue.

*Démonstration du théorème 3.3.* Vu le lemme 3.5, on prend deux structures rotatives d'amplitude supérieure à $-\pi$. Par une même isotopie préservant chaque tore $F_t$, on les déforme pour que les feuilletages caractéristiques de $F_0$ et $F_1$ aient des projections sur $F$ transversales l'une à l'autre. Le redressement du lemme 3.2 fournit alors, après choix d'un paramétrage de $F$ par $\mathbf{T}^2$, deux structures $\xi_0$ et $\xi_1$ ayant des équations du type

$$\cos\theta_i(x,t)\,dx_1 + \sin\theta_i(x,t)\,dx_2 = 0, \qquad i = 0, 1,$$

où, d'après la condition de contact, les fonctions $\theta_i \colon \mathbf{T}^2 \times [0,1] \to \mathbf{R}$ vérifient $\partial_t \theta_i < 0$.

Comme les amplitudes de $\xi_0$ et $\xi_1$ sont égales, il en est de même de toutes les variations ponctuelles

$$\begin{cases} \theta_0(x,1) - \theta_0(x,0) \\ \theta_1(x,1) - \theta_1(x,0) \end{cases} \quad x \in \mathbf{T}^2.$$

En effet, pour tout $t \in [0,1]$, la demi-orbite de $\xi_i F_t$ partant du point $x \in F_t \cong \mathbf{T}^2$ se relève dans le plan — par la projection $u \in \mathbf{R}^2 \mapsto x + u \in \mathbf{T}^2 = \mathbf{R}^2/\mathbf{Z}^2$ — en une courbe $C_i^t$ issue de l'origine 0 et paramétrée par $r \in [0, \infty[$. La variation angulaire ponctuelle et l'amplitude mesurent de combien tournent respectivement la tangente de $C_i^t$ à l'origine et la direction du vecteur limite $\lim_{r \to \infty} C_i^t(r)/r$. On voit donc comment passer continûment de la variation ponctuelle à l'amplitude par un arc $w_i(r)$, où $r$ varie de 0 à $\infty$. Comme les structures $\xi_0$ et $\xi_1$ coïncident au bord, les courbes $C_0^0$ et $C_1^0$ (resp. $C_0^1$ et $C_1^1$) sont identiques. Par suite, la fonction $w_1 - w_0$ prend ses valeurs dans $2\pi \mathbf{Z}$ et, comme elle est nulle à l'infini, elle l'est partout.

On peut ainsi prendre les fonctions $\theta_i$ égales sur $\partial V$, ce qui permet de les relier par un chemin de fonctions $\theta_s \colon V \to \mathbf{R}$, $s \in [0, 1]$, strictement décroissantes en $t$ et toutes identiques sur le bord $\partial V$. Les structures de contact $\xi_s$ d'équations

$$\cos \theta_s(x, t) \, dx_1 + \sin \theta_s(x, t) \, dx_2 = 0$$

forment, entre $\xi_0$ et $\xi_1$, une homotopie relative au bord que le théorème de Gray convertit en isotopie. $\square$

Grâce à un résultat de [Gi4], on peut aussi classifier les structures rotatives sur les variétés $\mathbf{T}_A^3$, $A \in \mathrm{SL}_2(\mathbf{Z})$. On rappelle que $\mathbf{T}_A^3$ est le quotient de $\mathbf{T}^2 \times \mathbf{R}$ par la transformation $(x, t) \mapsto (Ax, t+1)$, orienté par la forme $dx_1 \wedge dx_2 \wedge dt$. Étant donné une fonction $\theta \colon \mathbf{R} \to \mathbf{R}$, l'équation de Pfaff

$$\cos \theta(t) \, dx_1 + \sin \theta(t) \, dx_2 = 0, \qquad (x, t) \in \mathbf{T}^2 \times \mathbf{R},$$

définit une structure de contact sur $\mathbf{T}_A^3$ si et seulement si, pour tout $t \in \mathbf{R}$,

$$\theta'(t) < 0 \quad \text{et} \quad A^* \big( \cos \theta(t+1) \, dx_1 + \sin \theta(t+1) \, dx_2 \big) \wedge \big( \cos \theta(t) \, dx_1 + \sin \theta(t) \, dx_2 \big) = 0 \, .$$

On note $\zeta(\theta)$ cette structure.

**Théorème 3.7** [Gi4, proposition 2 et théorème 6].
**a)** *Chaque structure $\zeta(\theta)$ est tendue et sa classe d'isotopie ne dépend que du nombre $n \in \mathbf{N}$ qui est le plus grand entier strictement inférieur à*

$$\frac{1}{2\pi} \sup_{t \in \mathbf{R}} \big( \theta(t) - \theta(t+1) \big) \, .$$

*Dans la suite, $\zeta_n$ désigne un représentant $\zeta(\theta)$ quelconque de la classe d'isotopie associée à $n$.*
**b)** *Deux structures $\zeta_m$ et $\zeta_n$ sont toujours homotopes comme champs de plans mais ne sont conjuguées que si $m = n$.*

**Corollaire 3.8.** *Toute structure rotative sur $\mathbf{T}_A^3$ qui se relève à $\mathbf{T}^2 \times \mathbf{R}$ en une structure d'amplitude non nulle est isotope à l'une des structures $\zeta_n$.*

**Remarque 3.9.** Par un abus de langage commode, on dira parfois que deux structures rotatives sur $\mathbf{T}_A^3$ ont *la même amplitude* si elles sont isotopes à la même structure $\zeta_n$.

*Démonstration.* Soit $\xi$ la structure considérée sur $\mathbf{T}_A^3$ et $\widetilde{\xi}$ son rappel sur $\mathbf{T}^2 \times \mathbf{R}$. On se donne une fonction strictement décroissante $\theta \colon \mathbf{R} \to \mathbf{R}$ ayant les propriétés suivantes :
- $A^* \big( \cos \theta(t+1) \, dx_1 + \sin \theta(t+1) \, dx_2 \big) \wedge \big( \cos \theta(t) \, dx_1 + \sin \theta(t) \, dx_2 \big) = 0$ pour tout $t \in \mathbf{R}$ ;
- $\cos \theta(0) \, x_1 + \sin \theta(0) \, x_2 = 0$ définit la direction asymptotique de $\widetilde{\xi}(\mathbf{T}^2 \times \{0\})$ dans $H_1(\mathbf{T}^2; \mathbf{R})$ ;

— $\theta(1) - \theta(0)$ est l'amplitude de $\widetilde{\xi}$ restreinte à $\mathbf{T}^2 \times [0, 1]$.

La structure $\zeta(\theta)$ définie sur $\mathbf{T}_A^3$ par l'équation

$$\cos\theta(t)\, dx_1 + \sin\theta(t)\, dx_2 = 0, \qquad (x, t) \in \mathbf{T}^2 \times \mathbf{R},$$

est isotope à l'une des structures $\zeta_n$ (par définition de celles-ci). Comme dans la preuve du lemme 10 de [Gi4], on peut de plus déformer $\zeta(\theta)$, parmi les structures rotatives induisant un feuilletage de direction fixe sur $\mathbf{T}^2 \times \{0\}$, en une structure $\xi'$ qui imprime sur $\mathbf{T}^2 \times \{0\}$ le même feuilletage que $\xi$. Le théorème 3.3 s'applique alors aux restrictions de $\widetilde{\xi}$ et $\widetilde{\xi}'$ à $\mathbf{T}^2 \times [0, 1]$ — où $\widetilde{\xi}'$ est le rappel de $\xi'$ sur $\mathbf{T}^2 \times \mathbf{R}$ — et donne le résultat. □

**Corollaire 3.10.** *Sur $V = F \times [0, 1]$, deux structures rotatives qui coïncident au bord et ont des amplitudes différentes ne sont pas $\partial$-isotopes.*

La démonstration de ce corollaire utilise l'observation simple suivante :

**Lemme 3.11.** *Soit $\sigma_0$ et $\sigma_1$ des suspensions sur $F$. Pour tout entier $n \geq 0$, il existe sur $F \times [0, 1]$ une structure rotative d'amplitude comprise entre $-n\pi$ et $-(n+1)\pi$ qui imprime $\sigma_0$ sur $F_0$ et $\sigma_1$ sur $F_1$.*

*Démonstration.* Soit $\eta$ la structure rotative définie sur $\mathbf{T}^2 \times \mathbf{R}$ par

$$\cos\theta\, dx_1 - \sin\theta\, dx_2 = 0, \qquad (x, \theta) \in \mathbf{T}^2 \times \mathbf{R}.$$

Tout feuilletage non singulier $\sigma$ sur $\mathbf{T}^2$ admet une équation du type

$$\cos\theta(x)\, dx_1 - \sin\theta(x)\, dx_2 = 0, \qquad x \in \mathbf{T}^2,$$

et est donc le feuilletage caractéristique du graphe de la fonction $\theta$ dans $(\mathbf{T}^2 \times \mathbf{R}, \eta)$. En outre, toute suspension est isotope à un feuilletage $\sigma$ pour lequel la fonction $\theta$ oscille arbitrairement peu. Il existe donc un plongement incompressible $\phi \colon \mathbf{T}^2 \times [0, 1] \to \mathbf{T}^2 \times \mathbf{R}$ tel que la structure de contact $\phi^*\eta$ soit une structure rotative et imprime $\sigma_0$ sur $F_0$ et $\sigma_1$ sur $F_1$. On peut de plus choisir arbitrairement son amplitude parmi les valeurs compatibles avec les directions de $\sigma_0$ et $\sigma_1$. □

*Démonstration du corollaire 3.10.* On note $\xi_0$, $\xi_1$ les structures et $\sigma_0$, $\sigma_1$ les feuilletages qu'elles impriment sur $F_0$ et $F_1$ respectivement. Le lemme 3.11 fournit une structure rotative $\xi$ sur $F \times [1, 2]$ qui imprime $\sigma_1$ sur $F_1$ et $\sigma_0$ sur $F_2$. En empilant $\xi$ sur $\xi_0$ et $\xi_1$ puis en collant $F_2$ sur $F_0$ par l'identité, on obtient sur le tore $\mathbf{T}^3$ des structures rotatives $\overline{\xi}_0$ et $\overline{\xi}_1$ qui, vu le corollaire 3.8, sont respectivement isotopes à $\zeta_{n_0}$ et $\zeta_{n_1}$ où $n_i$, $i \in \{0, 1\}$, est le plus grand entier strictement inférieur à

$$\frac{1}{2\pi}\big|\mathrm{ampl}(\xi_i) + \mathrm{ampl}(\xi)\big|.$$

Comme les amplitudes de $\xi_0$ et $\xi_1$ diffèrent par un multiple non nul de $2\pi$, les entiers $n_0$ et $n_1$ sont distincts et le théorème 3.7-b) affirme alors que $\overline{\xi}_0$ et $\overline{\xi}_1$ ne sont pas isotopes. Il en résulte que $\xi_0$ et $\xi_1$ ne sont pas $\partial$-isotopes. □

## C  Structures élémentaires

On aborde ici la recherche de formes normales. Pour les structures de contact sur $V = F \times [0, 1]$, $F \simeq \mathbf{T}^2$, on impose une contrainte minime sur le bord qui simplifie les énoncés :

**Définition 3.12.** Dans toute la suite, on dit qu'un feuilletage de $F$ est *admissible* s'il est scindé par une multi-courbe essentielle (dont on indique souvent le nombre de composantes) ou s'il est topologiquement linéarisable (auquel cas on convient qu'il est scindé par 0 courbes essentielles). On dit aussi qu'une structure de contact sur $V$ est *admissible* (resp. est *admise*) si elle est directe, orientable et si le feuilletage caractéristique de chaque composante de $\partial V$ est admissible (resp. est une suspension admissible).

Étant donné des feuilletages admissibles $\sigma_0$ et $\sigma_1$ de $F$, on désigne désormais par $\mathcal{SCT}(V; \sigma_0, \sigma_1)$ l'espace des structures de contact tendues sur $V$ qui impriment $\sigma_0$ sur $F_0$ et $\sigma_1$ sur $F_1$. Le lemme ci-dessous montre que, pour comprendre la topologie de cet espace, on peut toujours supposer que $\sigma_0$ et $\sigma_1$ sont des suspensions. Autrement dit, la classification des structures de contact admissibles se ramène à celle des structures admises.

**Lemme 3.13.** *Soit $\sigma_0$, $\sigma_1$ et $\sigma$ des feuilletages admissibles de $F$. Si $\sigma_0$ et $\sigma$ sont scindés par une même multi-courbe non vide, les espaces $\mathcal{SCT}(V; \sigma_0, \sigma_1)$ et $\mathcal{SCT}(V; \sigma, \sigma_1)$ sont homéomorphes.*

*Démonstration.* Le lemme de réalisation 2.3 fournit sur $F \times [-1, 0]$ des structures $\xi_-$ et $\xi^-$ pour lesquelles les tores $F_t$, $t \in [-1, 0]$, sont tous convexes et qui vérifient les conditions suivantes :
$$\begin{aligned} \xi_- F_0 &= \sigma_0 \\ \xi_- F_{-1} &= \sigma \end{aligned} \quad \text{et} \quad \begin{aligned} \xi^- F_0 &= \sigma \\ \xi^- F_{-1} &= \sigma_0 \,. \end{aligned}$$

On dispose alors d'une application continue de $\mathcal{SCT}(V; \sigma_0, \sigma_1)$ dans $\mathcal{SCT}(V; \sigma, \sigma_1)$ qui consiste à coller $\xi_-$ sur tout élément $\xi$ de la source et à reparamétrer $[-1, 1]$ par $[0, 1]$. L'application dans l'autre sens, qui consiste à coller $\xi^-$, n'est pas exactement l'inverse mais la composée des deux est isotope à l'identité en vertu du lemme 2.6. □

Pour s'épargner des répétitions inutiles, on se place dorénavant dans un cadre un peu plus large. $V$ désigne une variété orientée fibrée en tores au-dessus de $K$ où $K$ est le cercle ou l'intervalle $[0, 1]$. On note $f$ la projection $V \to K$ et $F_t$, pour $t \in K$, sa fibre $f^{-1}(t)$.

**Définition 3.14.** Soit $\xi$ une structure de contact sur $V$, tendue ou non. On dit que $\xi$ est *élémentaire* — relativement à la fibration $f$ — si $\xi$ est admise et si, pour tout $t \in K$, les conditions suivantes sont remplies :
- l'union des variétés stables (resp. instables) des singularités positives (resp. négatives) de $\xi F_t$ forme un nombre fini de cercles lisses ;
- chacun de ces cercles est la trace sur $F_t$ d'un anneau transversal lisse qui coupe toute fibre voisine $F_s$ selon une courbe invariante de $\xi F_s$ — un cercle lisse formé de singularités et de variétés (in)stables ou une orbite fermée.

La seconde condition dit simplement que les cercles décrits dans la première (qui peuvent être entièrement constitués de singularités) changent par isotopie avec $t$, se transformant à l'occasion en orbites fermées.

**Proposition 3.15.** *Toute structure de contact tendue et admise sur $V$ est $\partial$-isotope à une structure élémentaire.*

*Démonstration.* Le lemme de préparation 2.16 permet de supposer que les tores $F_t$ dont le feuilletage caractéristique est singulier sont convexes. Il suffit alors de traiter le cas où $K$ est un intervalle $[t_0, t_1]$ et où les tores $F_t$, $t \in [t_0, t_1]$, sont tous convexes. Or, dans ce cas, le lemme est une conséquence directe des lemmes 2.3 et 2.6. En effet, pour un

paramétrage bien choisi

$$\phi\colon \mathbf{T}^2 \times K \longrightarrow V = F \times K, \quad (x,t) \longmapsto \phi_t(x) \in F_t$$

la multi-courbe qui scinde chaque feuilletage $\xi F_t$, $t \in K$, est l'image $\phi_t(\Gamma)$ d'une multi-courbe fixe $\Gamma \subset \mathbf{T}^2$. Pour chaque composante $S$ de $\mathbf{T}^2 \setminus \Gamma$, l'anneau $\phi_{t_i}(S)$, $i \in \{0,1\}$, contient une seule orbite fermée $C_i$ de $\xi F_{t_i}$ (le tore $f_{t_i} \subset \partial V$ est totalement réel) et on peut comparer les orientations de $C_0$ et $C_1$. Le lemme 2.3 permet de construire sur $\mathbf{T}^2 \times K$ une structure de contact $\eta$ ayant les propriétés suivantes :
- $\phi_*\eta$ coïncide avec $\xi$ près du bord de $V$ ;
- les tores $\mathbf{T}^2 \times \{t\}$ sont tous convexes et leur feuilletage $\eta(\mathbf{T}^2 \times \{t\})$ est scindé par $\Gamma \times \{t\}$ ;
- $\eta(\mathbf{T}^2 \times \{t\})$ est non singulier pour tout $t \neq t_{1/2} = (t_0 + t_1)/2$ ;
- $\eta(\mathbf{T}^2 \times \{t_{1/2}\})$ a exactement un cercle de singularités dans $S$ si et seulement si $C_0$ et $C_1$ ont des orientations opposées.

La structure $\phi_*\eta$ est alors $\partial$-isotope à $\xi$. □

Cette mise sous forme élémentaire fait apparaître dans $V$ des surfaces pré-lagrangiennes qu'on va exploiter pour identifier la structure de contact.

**Définition 3.16.** Soit $\xi$ une structure de contact élémentaire sur $V$. En réunissant les orbites fermées, les variétés stables des singularités positives et les variétés instables des singularités négatives de tous les feuilletages $\xi F_t$, $t \in K$, on obtient un sous-ensemble $\mathfrak{F}(\xi)$ de $V$ qu'on appelle le *feuillage* de $\xi$. Ses composantes connexes — les *feuilles* de $\xi$ — sont des surfaces (remarque 2.12) et plus précisément des tores, des bouteilles de Klein et éventuellement aussi des anneaux si $V$ a un bord. En fait, le feuillage $\mathfrak{F}(\xi)$ se décompose en quatre parties :
- l'ensemble $\mathfrak{F}_\parallel(\xi)$ des feuilles isotopes aux fibres ;
- l'ensemble $\mathfrak{F}_\perp(\xi)$ des feuilles $T$ — tores, bouteilles de Klein ou anneaux — pour lesquelles l'application $T/f \to K$ induite par $f$ est de degré non nul ;
- l'ensemble $\mathfrak{F}_o(\xi)$ des tores compressibles du feuillage ;
- l'ensemble $\mathfrak{F}_\partial(\xi)$ des anneaux $\partial$-compressibles (ou parallèles au bord).

Il est clair que $\mathfrak{F}_\partial$ est vide si $K = \mathbf{S}^1$ et que $\mathfrak{F}_\parallel$ et $\mathfrak{F}_\perp$ ne peuvent être simultanément non vides. En fait, si $\mathfrak{F}_\perp$ est non vide et si $K = \mathbf{S}^1$, la monodromie du fibré $V \to K$ préserve une courbe fermée simple essentielle et $V$ fibre donc en cercles au-dessus du tore ou de la bouteille de Klein.

**Lemme 3.17.** *Chaque feuille $T$ d'une structure de contact élémentaire est une surface pré-lagrangienne.*

*Démonstration.* La seule chose à voir est que $T$ est transversale à la structure de contact $\xi$, son feuilletage $\xi T$ étant alors formé des courbes legendriennes $T \cap F_t$, $t \in K$. On va vérifier cette propriété le long d'une composante connexe quelconque $C$ de l'intersection de $T$ avec un tore $F_{t_0}$, $t_0 \in K$. Si $C$ est une orbite fermée dégénérée de $\xi F_{t_0}$, la transversalité est manifeste car $T$ est tangente à $F_{t_0}$, lequel est transversal à $\xi$. On suppose donc maintenant que $C$ est une orbite fermée non dégénérée, ou une réunion de variétés stables (resp. instables) de singularités positives (resp. négatives) de $\xi F_{t_0}$. On paramètre un voisinage $U$ de $C$ par $\mathbf{S}^1 \times [-\varepsilon, \varepsilon] \times [t_0 - \varepsilon, t_0 + \varepsilon]$ de telle sorte que :
- $f|_U$ soit la projection sur $[t_0 - \varepsilon, t_0 + \varepsilon]$ ;
- $T \cap U$ soit l'anneau $\mathbf{S}^1 \times \{0\} \times [t_0 - \varepsilon, t_0 + \varepsilon]$.

Comme $C$ n'est pas une orbite fermée dégénérée de $\xi F_{t_0}$, les feuilletages $\xi(F_t \cap U)$ sont définis par des 1-formes $\beta_t$ sur $\mathbf{S}^1 \times [-\varepsilon, \varepsilon]$ dont la différentielle ne s'annule pas, et est par exemple positive. La structure $\xi$ a alors une équation du type $\beta_t + u_t\, dt = 0$, $t \in [t_0 - \varepsilon, t_0 + \varepsilon]$, et la condition de contact s'écrit

$$u_t\, d\beta_t + \beta_t \wedge \left(du_t - \dot\beta_t\right) > 0 \,.$$

Or, le long de $T \cap U = \mathbf{S}^1 \times \{0\} \times [t_0 - \varepsilon, t_0 + \varepsilon]$, le produit $\beta_t \wedge \dot\beta_t$ s'annule de sorte que l'inégalité ci-dessus y devient

$$u_t\, d\beta_t + \beta_t \wedge du_t > 0 \,.$$

Cette inégalité assure que le minimum de la fonction $u_{t_0}$ sur le cercle $\mathbf{S}^1 \times \{0\}$ est strictement positif, ce qui montre que la feuille $T$ est transversale à $\xi$ le long de $C$. □

**Définition 3.18.** Dans la suite, on dit qu'un ouvert connexe $J$ de $K$ est une *séquence de rotation* — pour une structure de contact donnée sur $V$ — si la restriction de $\xi$ à $f^{-1}(J)$ est une structure rotative d'amplitude non nulle. Si $J = \mathbf{S}^1$, cette dernière condition signifie que $\xi$ se relève à $F \times \mathbf{R}$ en une structure d'amplitude non nulle.

**Lemme 3.19.** *Soit $\xi$ une structure de contact élémentaire sur $V$.*
**a)** *Toute feuille $T$ de $\xi$ dont la projection sur $K$ est contenue dans une séquence de rotation est isotope aux fibres $F_t$.*
**b)** *Les feuilles de $\xi$ dont la projection sur $K$ n'est pas contenue dans une séquence de rotation sont en nombre fini.*

*Démonstration.*
**a)** C'est une conséquence directe du lemme 3.17 : $T$ est, comme les fibres $F_t$, transversale à $\xi$.
**b)** Il suffit d'observer que, si $J$ est un intervalle assez petit autour d'un point $t_0 \in K$, ou bien $J$ est une séquence de rotation, ou bien $f^{-1}(J)$ ne rencontre qu'un nombre fini de feuilles. Or cela est évident si $f_{t_0}$ est convexe et découle de la remarque 2.12 sinon. □

**Proposition 3.20.** *Au cours d'une homotopie de structures de contact élémentaires, homotopie relative au bord éventuel, les feuilles non isotopes aux fibres se déforment par isotopie.*

*Démonstration.* On note $\xi_s$, $s \in [0,1]$, l'homotopie des structures élémentaires et on choisit un instant quelconque $s_0 \in [0,1]$. La proposition est une conséquence directe des deux assertions suivantes qu'on justifie ci-après :
1) toute feuille $T_{s_0}$ de $\xi_{s_0}$ se déforme par isotopie en une feuille $T_s$ de $\xi_s$ pour $s$ dans un voisinage de $s_0$ ;
2) chaque fibre $F_{t_0}$ possède un voisinage dans lequel les seules feuilles de $\xi_s$ qui entrent, pour $s$ assez proche de $s_0$, sont isotopes aux fibres ou obtenues par déformation de feuilles de $\xi_{s_0}$ qui touchent $F_{t_0}$.

L'assertion 1) découle de la remarque 2.12 et du fait que chaque courbe $T_{s_0} \cap F_t$ qui n'est pas un cercle critique de $f \,|_{T_{s_0}}$ est une courbe invariante stable de $\xi_{s_0} F_t$. L'assertion 2) repose aussi sur le lemme 2.11 et la remarque 2.12.

Pour démontrer l'assertion 2), on note $S_{s_0}$ la réunion des feuilles de $\xi_{s_0}$ qui touchent $F_{t_0}$ et $S_s \subset \mathfrak{F}(\xi_s)$ sa déformation par isotopie pour $s$ dans un voisinage $I$ de $s_0$. Le lemme 2.11 et la remarque 2.12 montrent que, dans un intervalle $J$ autour de $t_0$, le complémentaire $J \setminus f(S_s)$ de $f(S_s)$ est une séquence de rotation (de $\xi_s$) pour tout $s$ assez proche de $s_0$.

Par suite, si $J' \subset J$ est un voisinage assez petit de $t_0$, les feuilles de $\xi_s$ dont la projection rencontre $J' \setminus f(S_s)$ sont entièrement incluses dans $F \times (J \setminus f(S_s))$ et sont donc isotopes aux fibres. D'autre part, l'ensemble

$$\bigcup_{\substack{s \in I \\ t \in f(S_s)}} (F_t \cap \mathfrak{F}(\xi_s)) \times \{s\}$$

est fermé dans $F \times J \times I$. L'assertion 1) assure alors que, pour $t \in J \cap f(S_s)$ et $s$ assez proche de $s_0$, chaque feuille de $\xi_s$ qui touche $F_t$ est contenue dans $S_s$, d'où l'assertion 2) et la proposition. □

La proposition 3.20 ci-dessus justifie la terminologie suivante :

**Définition 3.21.** On appelle *feuillage persistant* d'une structure élémentaire $\xi$ la réunion des ensembles $\mathfrak{F}_\perp(\xi)$, $\mathfrak{F}_o(\xi)$ et $\mathfrak{F}_\partial(\xi)$. On le note $\mathfrak{F}'(\xi)$. D'après le lemme 3.19, les feuilles persistantes sont en nombre fini.

## D  Formes normales des structures élémentaires

On décrit ci-dessous des isotopies qui permettent de mettre chaque structure de contact élémentaire sous une forme simple au sens où la projection $f$ a, sur le feuillage, un lieu critique le plus réduit possible. $V$ désigne toujours une variété orientée fibrée en tores au-dessus de $K$, où $K$ est le cercle ou l'intervalle $[0, 1]$.

**Proposition 3.22.** *Toute structure de contact élémentaire $\xi$ sur $V$ est $\partial$-isotope, parmi les structures élémentaires, à une structure $\xi'$ ayant les propriétés suivantes :*
- *sur chaque feuille de $\mathfrak{F}_o(\xi')$ (resp. de $\mathfrak{F}_\partial(\xi')$, resp. de $\mathfrak{F}_\perp(\xi')$), la fibration $f$ a exactement deux cercles critiques (resp. un seul, resp. aucun) ;*
- *chaque feuille de $\mathfrak{F}_\parallel(\xi')$ a une projection sur $K$ contenue dans une séquence de rotation ;*
- *les projections $f(S)$, $f(T)$ de deux feuilles distinctes quelconques $S, T$ ne s'intersectent que si c'est topologiquement inévitable, c'est-à-dire s'il n'y a aucune surface $S'$ isotope à $S$ relativement à $T \cup \partial V$ dont l'image $f(S')$ soit disjointe de $f(T)$.*

**Remarque 3.23.** Les configurations de feuilles $S, T$ pour lesquelles on ne peut disjoindre $f(S)$ et $f(T)$ sont clairement les suivantes :
- $S$ ou $T$ est dans $\mathfrak{F}_\perp$ ;
- $S$ et $T$ sont dans $\mathfrak{F}_\partial$ et s'appuient sur la même composante du bord de $V$ ;
- $S$ et $T$ sont dans $\mathfrak{F}_o \cup \mathfrak{F}_\partial$ et bordent — au besoin avec des anneaux de $\partial V$ — des tores pleins qui sont inclus l'un dans l'autre.

Pour toute feuille $T$ de $\mathfrak{F}(\xi)$, le lieu critique de $f|_T$ est une réunion de cercles. Ces cercles sont les composantes des intersections $T \cap F_t$ qui sont des orbites fermées dégénérées de $\xi F_t$. On obtient la forme normale promise dans la proposition ci-dessus en simplifiant au maximum ce lieu critique.

Soit $T$ une feuille et $R \subset T$ un anneau saturé par le feuilletage $\xi T$, on appellera *hauteur* de $R$ la différence $\sup \widetilde{f} - \inf \widetilde{f}$ où $\widetilde{f} \colon R \to \mathbf{R}$ désigne, selon que $K$ est un intervalle ou un cercle, la restriction ou un relèvement de $f$. On dira d'autre part que $R$ est *simple* s'il existe, le long de $R$, un champ de vecteurs $\nu$ qui est transversal à la fois aux tores $F_t$ et à $T$.

**Lemme 3.24.** *Soit $\xi$ une structure élémentaire, $T$ une de ses feuilles, $R \subset T \cap \operatorname{Int} V$ un anneau simple et $U \subset \operatorname{Int} V$ un voisinage de $R$ disjoint de $\mathfrak{F}(\xi) \setminus T$.*

*Il existe une isotopie $\phi_s$ de $V$, $s \in [0,1]$, ayant les propriétés suivantes :*
- *l'isotopie $\phi_s$ est arbitrairement $\mathcal{C}^0$-petite, à support dans $U$ et stationnaire sur $T$ ainsi que sur tout un voisinage du lieu critique de $f \mid_T$ ;*
- *pour tout $s \in [0,1]$, la structure $\phi_s^* \xi$ est élémentaire et a le même feuillage que $\xi$ ;*
- *$\phi_1^* \xi$ est transversale à $F_t$ le long de chaque composante de $R \cap F_t$.*

*Si l'anneau $R$ ne contient aucun cercle critique de $f \mid_T$, on peut troquer la dernière propriété ci-dessus contre la suivante :*
- *$\phi_1^* \xi$ est tangente à $F_t$ le long de chaque composante de $R \cap F_t$.*

*Démonstration.* Soit $\nu$ un champ de vecteurs legendrien transversal à $R$. La fonction $\nu \cdot f \colon R \to \mathbf{R}$ ne peut s'annuler sur le lieu critique de $f \mid_T$ et, comme $R$ est un anneau simple, elle est de signe constant sur ce lieu, positive pour fixer les idées. On construit ci-dessous l'isotopie $\phi_s$ de façon que la fonction $\nu \cdot (f \circ \phi_1^{-1})$ soit partout positive sur $R$. Pour cela, quitte à découper $R$ en sous-anneaux simples, on suppose que $R$ ne contient aucun cercle critique de $f \mid_T$ et se projette par $f$ sur un intervalle $[a,b] \subset K$. On met alors sur un voisinage de $R$ contenu dans $U$ des coordonnées

$$(x, y, z) \in \mathbf{S}^1 \times [-\varepsilon, \varepsilon] \times [a-\varepsilon, b+\varepsilon]$$

satisfaisant aux conditions suivantes :

1) $R_0 = \mathbf{S}^1 \times \{0\} \times [a-\varepsilon, b+\varepsilon]$ est un anneau simple de $T$ qui contient $R$ dans son intérieur et dont l'intersection avec $F_t$, $t \in [a-\varepsilon, b+\varepsilon]$, est la courbe d'équation $z = t$ ;

2) pour tout $y \in [-\varepsilon, \varepsilon] \setminus \{0\}$ et tout $t \in [a-\varepsilon, b+\varepsilon]$, l'anneau $R_y = \mathbf{S}^1 \times \{y\} \times [a-\varepsilon, b+\varepsilon]$ coupe transversalement $F_t$ suivant une courbe elle-même transversale à $\xi F_t$ ;

3) le champ de vecteurs $\partial_y$ est legendrien et coïncide avec $\nu$ le long de $R$.

Dans ces coordonnées, $\xi$ a une équation du type $dz - u(x,y,z)\,dx = 0$ où la fonction $u$ est nulle sur $R_0$. D'autre part, chaque tore $F_t$ a, près de $R_0$, une équation du type $z = h_t(x,y)$ où $h_t(x,0) = t$. Dans cette région, le feuilletage $\xi F_t$ est défini par

$$dh_t(x,y) - u\big(x,y,h_t(x,y)\big)\,dx = 0$$

et la condition 2) assure que la quantité

$$\partial_x h_t(x,y) - u\big(x,y,h_t(x,y)\big)$$

ne s'annule pas si $y \neq 0$.

On choisit alors un nombre positif $\delta < \varepsilon$ assez petit pour que

$$\sup\Big\{\big|h_t(x,y) - t\big|,\ (x,y,t) \in \mathbf{S}^1 \times [-\delta, \delta] \times [a-\delta, b+\delta]\Big\} < \varepsilon\,.$$

On se donne ensuite une fonction $\rho \colon \mathbf{R} \to \mathbf{R}$, à support dans $[-\delta, \delta]$, vérifiant

$$\rho(0) = 0,$$
$$\rho'(0) < -\sup\{\partial_y h_t(x,0),\ (x,t) \in \mathbf{S}^1 \times [a-\delta, b+\delta]\},$$
$$|\rho(y)| < \frac{\delta}{2} \inf\{\dot{h}_t(x,y),\ (x,t) \in \mathbf{S}^1 \times [a-\delta, b+\delta]\}\,.$$

On définit maintenant l'isotopie $\phi_s$ sur la partie des surfaces $F_t$ comprise entre $R_{-\delta}$ et $R_\delta$ par

$$\phi_s(x, y, h_t(x,y)) = (x, y, h_t(x,y) + s\lambda(t)\rho(y)), \qquad (x,y,t) \in \mathbf{S}^1 \times [-\delta, \delta] \times [a-\delta, b+\delta],$$

où $\lambda$ est une fonction lisse à support dans $[a-\delta, b+\delta]$, valant 1 sur $[a,b]$ et vérifiant $|\lambda'| \leq 2/\delta$.

Compte tenu de la majoration de $|\rho(y)|$, ce choix de $\lambda$ assure que $\phi_s$ est une isotopie. La propriété 2) garantit que les feuilletages des surfaces $\phi_s(F_t)$, $t \in [a-\delta, b+\delta]$, restent non singuliers dans $U \setminus R_0$ et transversaux aux courbes $R_y \cap \phi_s(F_t)$, $y \neq 0$. En outre, la majoration de $\rho'(0)$ implique que $\nu \cdot (f \circ \phi_1^{-1})$ est bien positive sur $R$.

Des arguments analogues permettent de construire une autre isotopie $\phi_s$ qui amène $\xi$ sur une structure $\phi_1^*\xi$ tangente à $F_t$ le long de chaque composante de $R \cap F_t$. $\square$

Avant de démontrer la proposition, on fait encore une observation simple :

**Lemme 3.25.** *Soit $\xi$ une structure élémentaire et $\phi_s$, $s \in [0,1]$, une isotopie de $V$ à support dans $U = f^{-1}(J)$, où $J$ est un intervalle de $K$. On suppose qu'il existe sur $U$ des coordonnées $(x,y,t) \in \mathbf{T}^2 \times J$ dans lesquelles les propriétés suivantes sont satisfaites :*
- *la fibration $f$ est la projection de $\mathbf{T}^2 \times J$ sur $J$ ;*
- *chaque cercle $\mathbf{S}^1 \times \{(y,t)\}$, $(y,t) \in \mathbf{S}^1 \times J$, est partout transversal à $\xi$ ou legendrien ;*
- *chaque surface $\phi_s(F_t)$, $(t,s) \in J \times [0,1]$, est le produit par $\mathbf{S}^1$ d'une courbe dans l'anneau des coordonnées $(y,t)$.*

*Les structures $\phi_s^*\xi$, $s \in [0,1]$, sont alors toutes élémentaires.* $\square$

*Démonstration de la proposition 3.22.* Pour fixer les idées, on suppose que $K$ est un cercle. En outre, on perturbe $\xi$ par une petite homotopie élémentaire pour que chaque feuilletage $\xi F_t$ possède au plus une orbite fermée dégénérée et que, sur chaque surface du feuillage, $f$ n'ait qu'un nombre fini de cercles critiques.

On se donne une feuille $T$ et un anneau $T' \subset T$ de plus petite hauteur parmi tous les anneaux bordés par deux cercles critiques consécutifs de $f|_T$. Des arguments classiques de topologie des courbes planes montrent que l'anneau $T'$ est simple, sauf si $T$ est un tore compressible sur lequel $f$ n'a que deux cercles critiques. On note $C$ (resp. $C'$) la composante de $\partial T'$ formée de maxima (resp. de minima) et $T''$ l'anneau de $T \setminus \text{Int } T'$ bordé par $C$ et le cercle critique suivant $C'''$, formé de minima. On relève les restrictions de $f$ à $T'$ et $T''$ en des fonctions $f' \colon T' \to \mathbf{R}$ et $f'' \colon T'' \to \mathbf{R}$ qui coïncident sur $C$ et on note $d$ la différence $f'(C') - f''(C''')$. Comme $T'$ est de hauteur minimale, $d$ est un réel positif ou nul. Par ailleurs, $d$ est un entier si et seulement si $C' = C'''$ car chaque feuilletage $\xi F_t$ compte au plus une orbite fermée dégénérée. L'entier $d$ est, dans ce cas, le degré de l'application $T/f \to \mathbf{S}^1$ induite par $f$.

Si $d = 0$, la feuille $T$ est un tore qui est soit compressible, soit isotope à une fibre $F_t$. Dans les deux cas, $f|_T$ n'a que deux cercles critiques. Si $T$ est compressible, on est dans la situation voulue. Si $T$ est isotope à une fibre et si sa projection $f(T)$ est disjointe de celle des autres feuilles, l'isotopie du lemme 3.24 transforment un voisinage de $f(T)$ en une séquence de rotation. Sinon, les isotopies construites ci-dessous permettent de disjoindre préalablement $f(T)$ des projections des autres feuilles.

On procède par récurrence sur la partie entière de $c = f'(C) - f'(C')$ et, pour commencer, on suppose $c < 1$. On se donne un réel positif $\varepsilon$ tel que $f(C')$ soit le seul niveau critique de $f$ sur le feuillage dans l'intervalle $[f(C') - \varepsilon, f(C') + \varepsilon]$ et on note $R$ le plus petit anneau de $T$ contenant $C \cup C'$ et dont les composantes de bord sont, l'une sur $F_{f(C')-\varepsilon}$, l'autre sur $F_{f(C')+\varepsilon}$. Comme l'anneau $R$ est simple (car $T'$ l'est), le lemme 3.24 permet de supposer que $\xi$ est transversale à $F_t$ le long de chaque courbe de $R \cap F_t$. On choisit alors un champ de vecteurs legendrien $\nu$ normal à $R$ et vérifiant $\nu \cdot f = 1$ et on munit un voisinage $U$ de $R$ de coordonnées $(x, y, z) \in \mathbf{S}^1 \times [-\delta, \delta] \times [-1, 1]$ ayant les caractéristiques suivantes :

- $R = \mathbf{S}^1 \times \{0\} \times [-1, 1]$ et $f|_R$ est une fonction $\widehat{f}(z)$ de la seule variable $z$ ;
- le point $(x, y, z)$ est l'image du point $(x, 0, z)$ par le flot de $\nu$ au temps $y$.

En particulier, chaque courbe $\mathbf{S}^1 \times \{y\} \times \{z\}$ est située sur $F_t$, où $t = \widehat{f}(z) + y$, et est transversale à $\xi F_t$ ; ainsi, $T$ est la seule feuille de $\xi$ qui rencontre $U$.

On prend alors pour $\phi_s$ une isotopie de $V$ ayant les propriétés suivantes :

4) dans $U$, chaque surface $\phi_s(F_t)$ est formée d'un ou deux anneaux $\mathbf{S}^1$-invariants et transversaux à $\nu$ (i.e. ayant une équation du type $y = g(z)$) ;

5) hors de $U$, le feuilletage par les surfaces $\phi_s(F_t)$, $t \in K = \mathbf{S}^1$, est indépendant de $s \in [0, 1]$ ;

6) chaque surface $\phi_1(F_t) \cap U$ est transversale à $R$.

Les propriétés 4) et 5) assurent que $\xi$ est élémentaire relativement à chaque fibration $f \circ \phi_s^{-1}$, $s \in [0, 1]$, avec toujours le même feuillage. La propriété 5) garantit aussi que, sur toute feuille autre que $T$, la fonction $f \circ \phi_s^{-1}$ a autant de cercles critiques que $f$. Par contre, grâce à la propriété 6), $f \circ \phi_1^{-1}$ a sur $T$ deux cercles critiques de moins que $f$.

On suppose maintenant $c > 1$. On note $\overline{C}'$ (resp. $\overline{C}''$) la première courbe d'intersection de $T'$ (resp. $T''$) avec $F_{f(C')}$ qu'on rencontre en partant de $C$. Une manipulation très voisine de celle effectuée plus haut, à présent autour de l'anneau $\overline{R} \subset T$ contenant $C$ et bordé par $\overline{C}' \cup \overline{C}''$, amène $f'(C)$ en-dessous de $f'(\overline{C}')$ et abaisse ainsi d'une unité la partie entière de $c = f'(C) - f'(C')$. À terme, on fabrique une structure élémentaire dont toute feuille $T$ non incluse dans l'image inverse d'une séquence de rotation a deux cercles critiques ou aucun selon que l'application $T/f \to K = \mathbf{S}^1$ induite par $f$ est ou non homotope à une constante[3]. Si deux feuilles ont des projections sur $K$ qui se rencontrent sans être dans l'une des configurations de la remarque 3.23, les isotopies du lemme 3.25 permettent facilement de modifier les hauteurs et de disjoindre les images. On obtient ainsi une structure qui vérifie les propriétés 1) et 3) de l'énoncé et une ultime application du lemme 3.24 donne 2). □

**Remarque 3.26.** Soit $\xi_0$ et $\xi$ deux structures de contact élémentaires sur $V$ qui coïncident au bord et sont sous forme normale. Si leurs feuillages persistants sont $\partial$-isotopes, les isotopies du lemme 3.25 permettent, comme ci-dessus, de déformer $\xi$ (relativement au bord et parmi les structures sous forme normale) en une structure $\xi_1$ dont le feuillage persistant coïncide avec celui de $\xi_0$.

**Corollaire 3.27.** *Soit $\xi_0$ et $\xi_1$ des structures de contact élémentaires sur $V$ qui coïncident au bord. Si $\mathfrak{F}_\perp(\xi_0)$ n'est pas vide et si $\mathfrak{F}(\xi_1)$ est $\partial$-isotope à $\mathfrak{F}(\xi_0)$, alors $\xi_1$ est $\partial$-isotope à $\xi_0$.*

*Démonstration.* Les propositions 3.20, 3.22 et la remarque 3.26 permettent de supposer que $\xi_0$ et $\xi_1$ sont sous forme normale et ont le même feuillage. Le reste de la démonstration utilise des arguments très semblables à ceux développés dans la partie 2 et on se borne à indiquer la marche à suivre. Les isotopies qu'on effectue à chaque étape laissent les structures sous forme normale et préservent leur feuillage.

---

[3]Lorsque $K = [0, 1]$, les anneaux parallèles au bord gardent un seul cercle critique.

Par une première modification, on fait coïncider $\xi_0$ et $\xi_1$ sur un voisinage du lieu critique de $f\,|\,_{\mathfrak{F}(\xi_0)}$ (il suffit de prendre un modèle local près d'une courbe legendrienne). Par une seconde modification, on fait coïncider $\xi_0$ et $\xi_1$ sur tout un voisinage des tores $F_t$ qui contiennent un cercle critique de $f\,|\,_{\mathfrak{F}(\xi_0)}$ (on invoque ici une version du lemme 2.8 pour des surfaces à bord, les structures en jeu étant égales sur le produit du bord par $\mathbf{R}$). On termine enfin à l'aide du lemme 2.6. □

Le corollaire ci-dessus est faux si $\mathfrak{F}_\perp(\xi_0)$ est vide. La raison est que le feuillage ne détermine pas l'amplitude de $\xi_0$ au-dessus des séquences de rotation (qui existent dès que $\xi_0$ est par exemple sous forme normale). On doit donc ajouter cette information et, en incorporant le théorème 3.3 à la démonstration du corollaire 3.27, on obtient le résultat suivant :

**Corollaire 3.28.** *Soit $\xi_0$ et $\xi_1$ des structures de contact élémentaires sur $V$ qui coïncident au bord et sont sous forme normale. On suppose que :*
  *– l'ensemble $\mathfrak{F}_\perp(\xi_0)$ est vide ;*
  *– les feuillages persistants $\mathfrak{F}'(\xi_0)$ et $\mathfrak{F}'(\xi_1)$ coïncident ;*
  *– pour toute composante connexe $J$ de l'ensemble*

$$K \setminus \mathrm{Int}\, f\bigl(\mathfrak{F}'(\xi_0)\bigr),$$

  *les restrictions de $\xi_0$ et $\xi_1$ à $f^{-1}(J)$ ont la même amplitude.*
*Alors $\xi_1$ est $\partial$-isotope à $\xi_0$.*

## E   Formes normales des structures universellement tendues

Un problème désormais important est d'identifier les structures de contact tendues parmi toutes les structures élémentaires. La proposition ci-dessous, qui est un corollaire immédiat des lemmes 3.33 et 3.35, caractérise les structures universellement tendues.

**Proposition 3.29.** *Une structure de contact élémentaire $\xi$ sur $V$ est universellement tendue si et seulement si son feuillage ne contient aucun tore compressible, autrement dit si et seulement si $\mathfrak{F}_o(\xi)$ est vide.*

**Exemple 3.30.** Pour tout entier $n > 0$, on note $\eta_n$ la structure de contact définie sur $\mathbf{T}^2 \times \mathbf{R}$ par l'équation de Pfaff

$$\cos(2n\pi y)\,dx + \sin(2n\pi y)\,dz = 0, \qquad (x,y,z) \in \mathbf{T}^2 \times \mathbf{R}.$$

C'est une structure universellement tendue car l'équation ci-dessus, vue sur $\mathbf{R}^3$, définit la structure de contact ordinaire. De plus, $\eta_n$ est élémentaire et son feuillage est $\mathbf{S}^1 \times \Delta_n$, où $\Delta_n \subset \mathbf{S}^1 \times \mathbf{R}$ est la réunion des $2n$ droites d'équation $2ny = 0 \pmod 1$.

Soit $n > 0$ un entier et $\phi \colon \mathbf{T}^2 \times [0,1] \to \mathbf{T}^2 \times \mathbf{R}$ un plongement incompressible du type

$$(x_1, x_2, t) \longmapsto \bigl(x_1, \hat\phi(x_2, t)\bigr) \in \mathbf{S}^1 \times \bigl(\mathbf{S}^1 \times \mathbf{R}\bigr).$$

La structure $\phi^*\eta_n$ est admise si et seulement si *chaque courbe paramétrée* $\hat\phi(\mathbf{S}^1 \times \{i\})$, $i \in \{0,1\}$, *coupe transversalement $\Delta_n$ avec, de plus, un vecteur tangent dont la composante sur $\partial_z$ est non nulle et de signe constant*. En outre, $\phi^*\eta_n$ est élémentaire dès qu'elle est admise et son feuillage est $\mathbf{S}^1 \times \Sigma$ où $\Sigma = \hat\phi^{-1}(\Delta_n)$.

On se donne maintenant, dans l'anneau $\mathbf{S}^1 \times [0,1]$, une famille $\Sigma$ de $2n + n_0 + n_1$ arcs proprement plongés et disjoints dont $2n$ joignent un bord à l'autre et $n_i$, pour $i \in \{0,1\}$, vont de $\mathbf{S}^1 \times \{i\}$ à lui-même. On fabrique facilement un plongement incompressible $\hat{\phi}$ de $\mathbf{S}^1 \times [0,1]$ dans $\mathbf{S}^1 \times \mathbf{R}$ qui remplisse la condition d'admission vue plus haut et vérifie $\hat{\phi}^{-1}(\Delta_n) = \Sigma$. La structure $\phi^*\eta_n$ est alors une structure élémentaire universellement tendue ayant pour feuillage $\mathbf{S}^1 \times \Sigma$.

L'enseignement de cet exemple est le suivant :

**Lemme 3.31.** *Soit $\mathfrak{A}$ une famille finie d'anneaux proprement plongés et disjoints dans $F \times [0,1]$, parmi lesquels $2n$, $n > 0$, vont de $F_0$ à $F_1$.*
**a)** *Pour $i \in \{0,1\}$, on se donne une suspension $\sigma_i$ sur $F_i$ dont les orbites fermées sont non dégénérées et sont les composantes de $F_i \cap \mathfrak{A}$. Il existe alors une structure de contact élémentaire $\xi$ sur $F \times [0,1]$ qui imprime $\sigma_i$ sur $F_i$, $i \in \{0,1\}$, et dont le feuillage est $\partial$-isotope à $\mathfrak{A}$.*
**b)** *On suppose que $\mathfrak{A}$ est le feuillage d'une structure élémentaire $\xi$ et ne contient aucun anneau allant de $F_0$ à lui-même. Tout plongement $\phi_0$ de $F_0$ dans une variété de contact $(M, \zeta)$ qui envoie $\xi F_0$ sur $\zeta \phi_0(F_0)$ se prolonge en un plongement de contact de $(F \times [0,1], \xi)$ dans un voisinage arbitrairement petit de $\phi_0(F_0)$.*

*Démonstration.*
**a)** Soit $\varepsilon > 0$ un nombre assez petit pour que $f|_{\mathfrak{A}}$ n'ait aucune valeur critique dans $[0, \varepsilon] \cup [1-\varepsilon, 1]$. L'exemple 3.30 fournit sur $F \times [\varepsilon, 1-\varepsilon]$ une structure de contact élémentaire $\xi$ dont le feuillage est isotope à $\mathfrak{A} \cap \big(F \times [\varepsilon, 1-\varepsilon]\big)$. On prolonge ensuite $\xi$ à $F \times [0,1]$ en une structure qui imprime la suspension $\sigma_i$ sur $F_i$, $i \in \{0,1\}$, et pour laquelle chaque tore $f_t$, $t \in [0, \varepsilon] \cup [1-\varepsilon, 1]$, est convexe.
**b)** La surface $\phi(F_0)$ est convexe et possède donc un système fondamental de voisinages

$$U \simeq \phi_0(F_0) \times \mathbf{R} \supset \phi_0(F_0) \times \{0\} = \phi_0(F_0)$$

dans lesquels la structure de contact $\eta$ est $\mathbf{R}$-invariante. Comme aucun anneau de $\mathfrak{A}$ ne va de $F_0$ à lui-même, l'exemple 3.30 et les lemmes de convexité — comme en a) — permettent de prolonger $\phi_0$ en un plongement $\phi \colon F \times [0,1] \to U$ pour lequel la structure élémentaire $\phi^*\zeta$ coïncide avec $\xi$ près du bord et a un feuillage $\partial$-isotope à $\mathfrak{A}$. Le corollaire 3.27 assure alors que $\phi^*\zeta$ est $\partial$-isotope à $\xi$. □

La partie b) du lemme ci-dessus a en particulier la conséquence suivante :

**Lemme 3.32.** *Soit $\xi$ une structure de contact sous forme normale sur $F \times [0,1]$ et $J$ l'intervalle fermé*

$$J = [0,1] \setminus \operatorname{Int} f\big(\mathfrak{F}_\partial(\xi)\big).$$

*Si la restriction de $\xi$ à un voisinage de $f^{-1}(J)$ est tendue (resp. universellement tendue), alors $\xi$ est tendue (resp. universellement tendue).* □

On démontre maintenant la proposition 3.29.

**Lemme 3.33.** *Soit $\xi$ une structure de contact élémentaire sur $V$. Si $\mathfrak{F}_o(\xi)$ est vide, $\xi$ est universellement tendue.*

*Démonstration.* Compte tenu des propositions 3.22 et 3.20, on suppose que $\xi$ est sous forme normale. Si $\mathfrak{F}_\partial(\xi)$ est vide (ce qui est en particulier le cas si $K = \mathbf{S}^1$), la situation est limpide : ou bien $\mathfrak{F}_\perp(\xi)$ n'est pas vide et tous les tores $F_t$ sont convexes, ou bien $\mathfrak{F}_\perp(\xi)$ est vide et $\xi$ est une structure rotative. Le lemme résulte alors du lemme 2.6 et du théorème 3.3. On suppose désormais que $\mathfrak{F}_\partial(\xi)$ est non vide, donc que $V = F \times [0,1]$.

Si $\mathfrak{F}_\perp(\xi)$ est non vide, on choisit un point $t \in [0,1] \setminus f\big(\mathfrak{F}_\partial(\xi)\big)$. Le tore $F_t$ est convexe et le lemme 3.31 montre que $(V, \xi)$ se plonge dans un voisinage arbitrairement petit de $F_t$. Par suite, la structure $\xi$ est universellement tendue.

Si $\mathfrak{F}_\perp(\xi)$ est vide, on pose $J = [0,1] \setminus \mathrm{Int}\, f\big(\mathfrak{F}_\partial(\xi)\big)$. La restriction de $\xi$ à un voisinage de $f^{-1}(J)$ est une structure rotative et est donc universellement tendue. Le lemme 3.32 assure alors que $\xi$ est aussi universellement tendue. $\square$

Avant la réciproque du lemme 3.33, on montre que certaines configurations de tores compressibles sont exclues dans le feuillage d'une structure élémentaire tendue.

**Lemme 3.34.** *Soit $\xi$ une structure de contact élémentaire sur $V$. Si le feuillage de $\xi$ contient un tore compressible dont la projection sur $K$ intersecte celle d'autres feuilles, $\xi$ est vrillée.*

*Démonstration.* Parmi les tores compressibles du feuillage dont la projection sur $K$ rencontre celle d'autres feuilles, l'un — au moins — borde un tore plein qui n'enferme aucune feuille ; on le note $T$. D'autre part, quitte à déformer $\xi$ comme dans les lemmes 3.24, 3.25 et dans la démonstration de la proposition 3.22, on peut supposer que les conditions suivantes sont remplies :
- $f|_T$ n'a que deux cercles critiques ;
- si $S$ est une feuille dont la projection $f(S)$ rencontre $f(T)$, alors en fait $f(S)$ recouvre un voisinage ouvert connexe $J$ de $f(T)$ et $f|_S$ n'a pas de valeurs critiques dans $J$ ;
- pour tout $t \in J$ et toute feuille $S \neq T$, la structure $\xi$ est tangente à $F_t$ le long de $F_t \cap S$.

On note alors $U$ l'intersection de $f^{-1}(J)$ avec la composante connexe de $\big(V \setminus \mathfrak{F}(\xi)\big) \cup T$ qui contient $T$. Par construction, $U$ est un tore plein ouvert dont le bord est la réunion de quatre anneaux dont deux sont inclus dans le feuillage et deux dans des fibres $F_t$. De plus, $T$ est la seule feuille qui rencontre $U$. On trouve alors sans peine un disque orienté $D$, méridien de $\mathrm{Adh}\, U$, ayant les propriétés suivantes :
- $D$ rencontre le tore plein bordé par $T$ suivant un disque $D'$ ;
- le feuilletage $\xi D$ sort (transversalement) le long de $\partial D$ et entre le long de $\partial D'$ ;
- les singularités de $\xi D$ sont négatives dans $D'$ et positives dans $D \setminus D'$.

Le bord orienté de $D$ est ainsi une courbe transversale à $\xi$ qui viole l'inégalité de Bennequin, de sorte que $\xi$ est vrillée. $\square$

**Lemme 3.35.** *Soit $\xi$ une structure élémentaire sur $V$. Si $\xi$ est tendue et si son feuillage contient un tore compressible, $\xi$ est virtuellement vrillée.*

*Démonstration.* Vu les propositions 3.20 et 3.22, on suppose $\xi$ sous forme normale. La multiplication par un entier $n > 1$ dans $\pi_1(F_t) \simeq \mathbf{Z}^2$ induit un morphisme injectif $\pi_1(V) \to \pi_1(V)$. Sur le revêtement $\widetilde{V}$ de $V$ associé à ce morphisme, $\xi$ induit une structure de contact élémentaire $\widetilde{\xi}$ dont le feuillage $\mathfrak{F}(\widetilde{\xi})$ est l'image inverse de $\mathfrak{F}(\xi)$. Or, l'image inverse d'un tore compressible de $\mathfrak{F}(\xi)$ est formée de $n$ tores compressibles qui ont tous la même projection composée sur $K$ (en particulier, $\widetilde{\xi}$ n'est pas sous forme normale). D'après le lemme 3.35, la structure $\widetilde{\xi}$ est vrillée. $\square$

## F    Formes normales des structures virtuellement vrillées

On montre ici qu'une structure de contact virtuellement vrillée a presque toujours plusieurs formes normales (aux feuillages persistants non isotopes) et que celles-ci ont des combinatoires très particulières.

Soit $\sigma_0$ et $\sigma_1$ des suspensions admissibles sur $F$ scindées respectivement par $2n_0$ et $2n_1$ courbes essentielles. On note :
- $D_i$, $i \in \{0, 1\}$, la direction asymptotique de $\xi F_i$ dans $H_1(F; \mathbf{R})$ ;
- $C = C(\sigma_0, \sigma_1) \subset H_1(F; \mathbf{R}) \setminus \{0\}$ l'un des deux cônes convexes bordés à gauche[4] par $D_0$, à droite par $D_1$ et *ouvert* (resp. *fermé*) du côté de $D_i$ si $n_i = 0$ (resp. si $n_i \geq 1$) ;
- $E = E(\sigma_0, \sigma_1)$ l'enveloppe convexe de $C \cap H_1(F; \mathbf{Z})$ ;
- $B = B(\sigma_0, \sigma_1)$ l'ensemble des points entiers primitifs de $\partial E$, ordonné par la relation $w \preceq w'$ si $w \wedge w' \leq 0$.

De plus, pour toute partie finie $Q$ de $B$, on pose

$$\gamma(Q) = \sum_{j=1}^{k} (-1)^{j-1} w_j \quad \text{où} \quad Q = \{w_1, w_2, \ldots, w_k\} \text{ avec } w_1 \prec w_2 \prec \cdots \prec w_k.$$

On prend maintenant dans $\mathcal{SCT}_v(V; \sigma_0, \sigma_1)$ une structure de contact $\xi$ sous forme normale. D'après la proposition 3.29 et le lemme 3.34, $\mathfrak{F}_o(\xi)$ est non vide mais $\mathfrak{F}_\perp(\xi)$ est vide. De plus, tous les feuilletages $\xi F_t$ ont une direction asymptotique $D_t$ bien définie et localement constante sur la projection $f(\mathfrak{F}'(\xi))$ du feuillage persistant $\mathfrak{F}'(\xi) = \mathfrak{F}_o(\xi) \cup \mathfrak{F}_\partial(\xi)$.
— La direction asymptotique du feuilletage d'un tore convexe est la droite contenant la classe des courbes qui scindent le feuilletage. —

**Proposition 3.36.** *Soit $\xi \in \mathcal{SCT}_v(V; \sigma_0, \sigma_1)$ une structure de contact sous forme normale. Quand $t$ varie de $0$ à $1$, la direction asymptotique $D_t$ de $\xi F_t$ fait au plus une fois le tour de la droite projective de $H_1(F; \mathbf{R})$. De plus, pour $t \in f(\mathfrak{F}_o(\xi))$, le vecteur directeur primitif de $D_t$ situé dans $\operatorname{Int} C$ est sur le bord de l'enveloppe convexe $E$.*

En fait, il s'ensuit des propositions 1.7 et 3.40 que les propriétés ci-dessus caractérisent les formes normales des structures de contact virtuellement vrillées. D'autre part, ces propriétés sont très restrictives ; elles entraînent par exemple directement que toute structure de contact dans $\mathcal{SCT}(V; \sigma_0, \sigma_1)$ a une torsion nulle. Elle montrent aussi accessoirement que, vu les conditions au bord mises dans la définition du cône $C$, les droites qui rencontrent $C$ sont exactement les droites $D_t$, $t \in \,]0, 1[$.

Pour prouver la seconde affirmation de la proposition, on a besoin de caractériser les points entiers de $\operatorname{Int} E$ :

**Lemme 3.37.** *Un point primitif $w \in \operatorname{Int} C$ est dans $\operatorname{Int} E$ si et seulement s'il existe un point entier $v \in C$ tel que $w - v$ soit aussi dans $C$.*

*Démonstration.* Chaque arête de $\partial E$ non contenue dans $\partial C$ a pour direction une droite vectorielle qui évite $C$. Par suite, si $w \in \partial E$ et si $v \in C$ est un point entier quelconque, $w - v$ est hors de $C$. Inversement, si $w \in \operatorname{Int} E$, le segment $[0, w]$ coupe $\partial E$ en un point $x$ et on note $v$ (resp. $u$) le point entier le plus proche de $x$ sur $\partial E$ tel que $w \wedge v > 0$ (resp. $w \wedge u < 0$). Comme $v$ et $u$ sont des points entiers consécutifs sur $\partial E$, le triangle $[0, u, v]$ ne contient aucun point entier autre que ses sommets et, ainsi, $u \wedge v = 1$. De plus, $w$ est aussi un point entier donc $u \wedge w \geq 1$ et, du coup, $u \wedge (w - v) \geq 0$. D'autre part, $v \wedge (w - v) = v \wedge w < 0$ de sorte que $w - v$ est dans $C$ puisque $u$ et $v$ y sont. □

*Démonstration de la proposition 3.36.* Soit $T_1, \ldots, T_k$ les composantes de $\mathfrak{F}_o(\xi)$ indexées de telle sorte que les intervalles $f(T_j)$ soient rangés en ordre croissant.

---
[4]Comme le tore $F$ est orienté, le plan $H_1(F; \mathbf{R})$ est orienté par la forme d'intersection.

La fonction
$$\delta : [0,1] \longrightarrow \mathbf{P}^1 = \mathbf{P}\big(H_1(F, \mathbf{R})\big), \quad t \longmapsto D_t,$$
est continue, décroissante au sens large et localement constante sur $f\big(\mathfrak{F}'(\xi)\big)$. Si $D_t$ fait strictement plus d'un tour sur $\mathbf{P}^1$, on peut trouver un intervalle $\Delta$ non ponctuel dans $\mathbf{P}^1$ dont l'image inverse par $\delta$ ne soit pas connexe et qui, de plus, ne contienne aucune des droites $D_t$, $t \in f\big(\mathfrak{F}_0(\xi)\big)$. On choisit alors deux composantes connexes consécutives $I'$, $I''$ de $\delta^{-1}(\Delta) \subset [0,1]$ et on note $J$ l'ensemble des entiers $j \in [1, k]$ pour lesquels $f(T_j)$ est situé entre $I'$ et $I''$. Pour tout $j \in J$, la droite rationnelle $D_t$, $t \in f(T_j)$, a donc un vecteur directeur primitif $w_j$ dans $\operatorname{Int} C$. On se donne ensuite un vecteur primitif $v$ tel que la droite $D = \mathbf{R}v$ soit dans $\Delta$ et que le déterminant
$$v \wedge \sum_{j \in J} (-1)^{j-1} w_j$$
soit non nul. On prend enfin des points $t' \in \delta^{-1}(D) \cap I'$, $t'' \in \delta^{-1}(D) \cap I''$ et on regarde dans $F \times [t', t'']$ un anneau $A$ transversal aux tores $F_t$, $t \in [t', t'']$, dont les composantes de bord sont des orbites fermées respectives de $\xi F_{t'}$ et $\xi F_{t''}$. Pour cet anneau à bord legendrien, l'inégalité de Bennequin dit que les indices des singularités positives de $\xi A$ ont une somme égale à celle des indices des singularités négatives. Or la différence entre ces deux sommes vaut ici $\pm 2 v \wedge \sum_{j \in J} (-1)^{j-1} w_j$ et est donc non nulle. Par conséquent, $\xi$ est vrillée.

Pour tout $j \in \{1, \ldots, k\}$, on sait désormais que la droite $D_t$, $t \in f(T_j)$, a un vecteur directeur primitif dans $\operatorname{Int} C$ et on note $w_j$ ce vecteur. On suppose que les points $w_j$ ne sont pas tous sur $\partial E$ et on montre que $\xi$ est vrillée. Pour chaque $w_j \in \operatorname{Int} E$, le lemme 3.37 fournit un point entier $v_j \in C$ tel que $w_j - v_j$ soit aussi dans $C$. Au besoin, on permute $v_j$ avec $w_j - v_j$ et on remplace $v_j$ par un autre point entier du triangle $[0, w_j, v_j]$ pour imposer $w_j \wedge v_j = 1$, condition qui fige le choix de $v_j$. On note ensuite $P_j$ le parallélogramme $[0, v_j, w_j, w_j - v_j] \subset C$ et on observe que $\operatorname{Int} P_j$ ne contient aucun point entier, donc *a fortiori* aucun point $w_j$.

**Assertion.** *Il existe un point $w_i \in \operatorname{Int} E$ pour lequel $\mathbf{R} w_i$ est, parmi les droites $\mathbf{R} w_j$, la seule qui rencontre $\operatorname{Int} P_i$.*

*Preuve.* Quel que soit $w_i \in \operatorname{Int} E$, les entiers $j$ pour lesquels la droite $\mathbf{R} w_j$ rencontre $\operatorname{Int} P_i$ forment un intervalle $[i_0, i_1] \cap \mathbf{Z} \ni \{i\}$. On choisit pour $w_i$ un point tel que le nombre $i_1 - i_0 + 1$ de ces droites soit minimal et on montre que ce minimum vaut 1.

Si $i_1 - i_0 \geq 1$, l'un des nombres $i - i_0$, $i_1 - i$ est non nul et, pour fixer les idées, on suppose que c'est $i - i_0$. On distingue alors deux cas, selon qu'il existe ou non un entier $j \in [i_0, i-1]$ tel que $w_i \wedge w_j > 1$. Si un tel $w_j$ existe, $w_j - v_i$ appartient à $C$ et le triangle $[0, w_j, v_i]$ contient donc le point $v_j$. Du coup, $\operatorname{Int} P_j$ rencontre au plus $i - i_0$ des droites $\mathbf{R} w_j$ et $i_1 - i_0 + 1$ n'est pas minimal. Si $w_i \wedge w_j = 1$ pour tout entier $j \in [i_0, i-1]$, alors $\mathbf{R} w_{i-1}$ est la seule des droites $\mathbf{R} w_j$ à rencontrer $\operatorname{Int} P_{i-1}$ et $i_1 - i_0 + 1$ est encore moins minimal. □

Dans la suite de la preuve de la proposition 3.40, $w_i$ est un point de $\operatorname{Int} E$ pour lequel $\mathbf{R} w_i$ est, parmi les droites $\mathbf{R} w_j$, la seule qui rencontre $\operatorname{Int} P_i$. D'autre part, on oriente $\xi$ de telle sorte que $w_i$ soit la classe d'homologie des orbites fermées de $\xi F_{r_i}$ où $r_i = \inf f(T_i)$. On peut alors trouver un paramétrage
$$\mathbf{T}^2 \times [0,1] \longrightarrow V = F \times [0,1]$$
qui respecte la projection sur $[0,1]$ et donne les propriétés suivantes :

1) $w_i = (1,0)$ et $v_i = (0,1)$ ;
2) il existe des points $t' < t'' \in [0,1]$ situés de part et d'autre de $f(T_i)$ tels que les feuilletages $\xi F_{t'}$ et $\xi F_{t''}$ aient pour directions asymptotiques respectives les droites $\mathbf{R}v_i$ et $\mathbf{R}(w_i - v_i)$ ;
3) chaque anneau $A_x = \{x\} \times \mathbf{S}^1 \times [t', t'']$, $x \in \mathbf{S}^1$, muni de l'orientation induite par $\partial_x$, a deux points complexes — un point elliptique positif et un point hyperbolique négatif ;
4) chaque courbe $A_x \cap F_t$ est transversale à $\xi$ pour $t \in [s_i, t'']$ et, pour $t \in [t', r_i]$, est soit transversale, soit legendrienne.

En outre, quitte à déformer $\xi$ par une $\partial$-isotopie $\mathcal{C}^\infty$-petite — parmi les structures sous forme normale —, on suppose que les tores $F_{t'}$ et $F_{t''}$ sont convexes.

Pour tout $x \in \mathbf{S}^1$, l'une des séparatrices instables du point hyperbolique de $A_x$ coupe $F_{t''}$ en un point $g(x)$ et les points $g(x)$, $x \in \mathbf{S}^1$, décrivent une courbe orientée $G$ sur $F_{t''}$ dont la classe d'homologie est $(1,0)$. Par suite, $G$ coupe chaque orbite fermée de $\xi F_{t''}$. On choisit alors un point $g(x_0)$ situé sur une orbite fermée répulsive $L$ de $\xi F_{t''}$. Pour $x_1 > x_0$ et $x_1 - x_0$ assez petit, la demi-orbite de $\xi F_{t''}$ qui part de $g(x_1)$ revient taper $G$ en un point $g(x_2)$ avec $x_2 > x_1$. On note $I$ l'intervalle du cercle $[x_2, x_1] \ni x_0$ et on observe que l'ensemble polyédral

$$F_{t'} \cup A_{x_1} \cup A_{x_2} \cup \left( I \times \mathbf{S}^1 \times \{t''\} \right)$$

contient un anneau polyédral plongé $A \supset I \times \mathbf{S}^1 \times \{t''\}$ dont le feuilletage $\xi A$ a les propriétés suivantes :
- $\xi A$ a exactement deux orbites fermées — les composantes du bord — dont l'une est attractive et l'autre répulsive ;
- $\xi A$ a quatre singularités, un foyer et une selle de chaque signe, qu'on note respectivement $a^+$, $a^-$, $b^+$ et $b^-$ ;
- les deux séparatrices stables (resp. instables) de $b^-$ (resp. de $b^+$) viennent de $a^+$ (resp. vont vers $a^-$) et forment un lacet non contractile ;
- l'une des séparatrices instables (resp. stables) de $b^-$ (resp. de $b^+$) va vers le cycle attractif (resp. vient du cycle répulsif) tandis que l'autre est en connexion rétrograde avec l'autre selle.

En créant deux singularités positives sur le bord puis en éliminant $b^+$ avec le nouveau foyer, on fait apparaître un disque vrillé (voir la démonstration du lemme de préparation). □

On discute maintenant la multiplicité et l'équivalence des formes normales. Pour des raisons techniques, on ne traite que les cas où $n_0$, $n_1$ valent 0 ou 1. On commence par évaluer la classe d'Euler relative d'une structure de contact $\xi \in \mathcal{SCT}(V; \sigma_0, \sigma_1)$ qui est sous forme normale. À cette fin, on oriente $\sigma_0$ — de telle sorte que ses cycles asymptotiques soient sur le bord gauche du cône $C$ — et $\xi$ de façon compatible. D'après la proposition 3.36, toutes les droites $D_t$, $t \in f(\mathfrak{F}_o(\xi))$, intersectent $\partial E$ en un point entier ; on désigne par $R_0(\xi) \subset B$ l'ensemble (fini) de ces points.

Lorsque $n_0 = n_1 = 0$, la classe d'Euler relative de $\xi$ vaut

$$\chi_\partial(\xi) = 2\gamma\big(R_0(\xi)\big).$$

En effet, il existe une section de $\xi$ qui, sur chaque tore $F_t$, engendre le feuilletage $\xi F_t$. Or, quitte à déformer $\xi$ parmi les structures de contact sous forme normale, on peut supposer que, pour toute composante $T$ de $\mathfrak{F}_o(\xi)$, l'intervalle $f(T) = [r, s]$ contient un seul instant $t$ où $\xi F_t$ est singulier : cela résulte du lemme 2.6 et du fait que les orbites fermées de $\xi F_r$ et $\xi F_s$ ont des orientations opposées (cf. lemme 3.17). Les singularités de $\xi F_t$ forment alors deux courbes où la section de $\xi$ s'annule transversalement et qui, comme courbes de zéros, ont la même classe d'homologie que l'orbite fermée de $\xi F_r$.

Lorsque $n_i = 1$, soit $T$ la composante de $\mathfrak{F}_\partial(\xi)$ qui s'appuie sur $F_i$ et $[r,s] = f(T)$ sa projection. Le raisonnement ci-dessus montre que, si les orbites fermées de $\xi F_r$ et $\xi F_s$ ont des orientations opposées, $T$ compte comme un tore compressible dans le calcul de la classe d'Euler relative. On dira dans ce cas que l'anneau $T$ est *impair*. Cette définition et l'observation qui la fonde s'appliquent aussi aux anneaux de $\mathfrak{F}_\perp(\xi)$ lorsque $n_0 = n_1 = 1$ et que $\mathfrak{F}_\partial(\xi)$ est vide.

**Définition 3.38.** Soit $\sigma_0$, $\sigma_1$ des suspensions admissibles scindées par 0 ou 2 courbes essentielles et soit $\xi \in \mathcal{SCT}(V; \sigma_0, \sigma_1)$ une structure de contact sous forme normale. On appelle *lieu de retournement* l'ensemble des points de $B$ situés sur une droite $D_t$ où $t$ appartient à la projection soit d'un tore de $\mathfrak{F}_o(\xi)$, soit d'un anneau impair de $\mathfrak{F}_\partial(\xi)$ ou de $\mathfrak{F}_\perp(\xi)$.

**Lemme 3.39.** *Soit $\sigma_0$, $\sigma_1$ des suspensions admissibles scindées par 0 ou 2 courbes essentielles. La classe d'Euler relative de toute structure de contact $\xi$ sous forme normale dans $\mathcal{SCT}(V; \sigma_0, \sigma_1)$ vaut*
$$\chi_\partial(\xi) = 2\gamma\big(R(\xi)\big)$$
*où $R(\xi)$ est le lieu de retournement de $\xi$.*

On observe au passage que le lieu de retournement $R(\xi)$ est de cardinal pair ou impair selon que les cycles asymptotiques du feuilletage orienté $\xi F_1$ sont ou non sur le bord de $C$.

**Proposition 3.40.** *Soit $\sigma_0$, $\sigma_1$ des suspensions admissibles scindées par 0 ou 2 courbes essentielles et $\xi$ une structure de contact dans $\mathcal{SCT}_v(V; \sigma_0, \sigma_1)$. Pour tout ensemble fini $Q \subset B$ tel que $2\gamma(Q) = \chi_\partial(\xi)$, il existe une structure de contact sous forme normale et $\partial$-isotope à $\xi$ dont le lieu de retournement est $Q$.*

La preuve de cette proposition repose sur l'étude détaillée d'un exemple suggéré par une construction de V. COLIN [Co2].

**Exemple 3.41.** Soit $\zeta$ la structure de contact définie sur $\mathbf{R}^2 \times \mathbf{S}^1$ par
$$\cos(2\pi z)\, dx - \sin(2\pi z)\, dy = 0, \qquad (x, y, z) \in \mathbf{R}^2 \times \mathbf{R}/\mathbf{Z},$$
et $\widehat{\phi}_0$ le plongement
$$\begin{cases} \mathbf{T}^2 \longrightarrow \mathbf{R}^2 \times \mathbf{S}^1 = \mathbf{C} \times \mathbf{S}^1, \\ (x_1, x_2) \longmapsto \widehat{\phi}_0(x_1, x_2) = \Big(e^{-2i(n+1)\pi x_2}\big(1 - e^{2i\pi x_1}\big), x_2\Big), \end{cases} \qquad n \geq 1.$$

Un calcul direct montre que le feuilletage caractéristique de $\widehat{\phi}_0(\mathbf{T}^2)$ est dirigé par le champ de vecteurs
$$-2(n+1)\sin(\pi x_1)\cos(2n\pi x_2 - \pi x_1)\,\partial_{x_1} + \sin(2n\pi x_2 - 2\pi x_1)\,\partial_{x_2}.$$

Ce champ de vecteurs a $4n$ singularités, $n$ selles et $n$ foyers de chaque signe. Les selles négatives (resp. positives) sont les points de coordonnées $(0, m/n)$ (resp. $(0, (2m+1)/2n)$) et les foyers négatifs (resp. positifs) sont les points $(1/2, (2m+1)/2n)$ (resp. $(1/2, m/n)$), $0 \leq m \leq n - 1$. De plus, la courbe d'équation $x_1 = 0$ est formée de $2n$ connexions de selles qui sont rétrogrades, *i.e.* vont des selles négatives aux selles positives. Dans la suite, on dira qu'une quelconque de ces connexions est de signe $+$ ou $-$ selon qu'elle relie un point $(0, m/n)$ au point $(0, (2m+1)/2n)$ ou au point $(0, (2m-1)/2n)$. Enfin, le feuilletage $\zeta \, \phi_0(\mathbf{T}^2)$ n'a aucune orbite fermée et chaque selle $(0, m/2n)$ est reliée au foyer $(1/2, m/2n)$ par deux séparatrices.

Afin d'enlever des symétries, on prend deux entiers $p, q \in \{1, \ldots, n\}$ et on déforme $\widehat{\phi}_0$ par une petite isotopie en un plongement $\phi_0$ dont l'image ne garde, dans son feuilletage caractéristique, que $q$ (resp. $p$) des $n$ connexions positives (resp. négatives) de $\zeta \, \widehat{\phi}_0(\mathbf{T}^2)$. Pour cela, on choisit arbitrairement $n - q$ (resp. $n - p$) des $n$ points de coordonnées $(0, (4m + 1)/4n)$ (resp. $(0, (4m - 1)/4n)$) et on pousse $\widehat{\phi}_0(\mathbf{T}^2)$ vers le bas près des $n - q$ points situés sur des connexions positives et vers le haut près des $n - p$ autres.

Soit alors $\phi \colon \mathbf{T}^2 \times \mathbf{R} \to \mathbf{R}^2 \times \mathbf{S}^1$ un plongement dont la restriction à $\mathbf{T}^2 \times \{0\}$ est $\phi_0$. Le lemme 2.13 montre que, pour $\varepsilon > 0$ assez petit, tous les feuilletages

$$\zeta\Big(\phi\big(\mathbf{T}^2 \times \{t\}\big)\Big), \qquad t \in [-\varepsilon, \varepsilon] \setminus \{0\},$$

sont de type Morse-Smale, de sorte que les tores $\phi(\mathbf{T}^2 \times \{t\})$ sont convexes. En outre, il est facile de voir que $\phi^* \zeta$ est virtuellement vrillée sur $\mathbf{T}^2 \times [-\varepsilon, \varepsilon]$ (voir [Co2]).

On revient maintenant au cadre habituel en posant $F = \mathbf{T}^2$ et en identifiant $[0, 1]$ à $[-\varepsilon, \varepsilon]$ par le difféomorphisme $t \mapsto \varepsilon(2t - 1)$. On note $\xi$ la structure de contact ainsi obtenue sur $V = F \times [0, 1]$ et on l'oriente pour que le feuilletage

$$\xi F_{1/2} = \phi_0^*\big(\zeta \, \phi_0(\mathbf{T}^2)\big)$$

ait les caractéristiques décrites plus haut. Les feuilletages $\sigma_i = \xi F_i, i \in \{0, 1\}$, sont singuliers mais ils sont scindés par deux courbes essentielles et leurs directions asymptotiques $D_i$ sont respectivement engendrées par $(-p, 1)$ et $(q, 1)$ (par construction de $\phi_0$). On prend alors pour $C$ le cône

$$C = \big\{(x_1, x_2) \in \mathbf{R} \times \,]0, \infty[ \, \subset H_1(\mathbf{T}^2; \mathbf{R}) = \mathbf{R}^2 \mid -px_2 \leq x_1 \leq qx_2\big\}.$$

Dans la suite, on confectionne explicitement diverses structures de contact sous forme normale[5] $\partial$-isotopes à $\xi$.

La remarque initiale est que le film $\xi F_t$, $t \in [0, 1]$, n'est pas du tout générique : dans un film générique, les connexions de selles arrivent une par une. Une approximation générique du film $\xi F_t$ est donc codée (du moins topologiquement) par une énumération des connexions de $\xi F_{1/2}$ (en fait, par le mot de $p + q$ lettres en $+$ et $-$ qui affiche les signes des connexions successives). On choisit alors une telle énumération et on note $\xi'$ une déformation de $\xi$ qui la réalise, *i.e.* dont le film fait apparaître les connexions dans l'ordre prescrit, à des instants $t_1 < \cdots < t_{p+q}$.

Pour mettre $\xi'$ sous forme normale, la proposition 3.15 montre qu'il faut essentiellement éliminer les singularités des feuilletages $\xi' F_{t_i}$. On choisit pour cela un nombre $\delta > 0$ tel que les intervalles $[t_i - 2\delta, t_i + 2\delta]$, $1 \leq i \leq p + q$, soient disjoints et inclus dans $[0, 1]$. Comme chaque feuilletage $\xi' F_{t_i}$ présente une et une seule connexion de selles, on peut grouper ses singularités par paires en position d'élimination et ce, de façon unique.

---
[5]et donc élémentaires, à ceci près qu'elles ne sont pas admises puisque les composantes du bord gardent leur feuilletage caractéristique singulier.

Précisément, selon que la connexion est de signe $+$ ou $-$, on apparie $(0, m/2n)$ avec $(1/2, (2m-1)/2n)$ ou avec $(1/2, (2m+1)/2n)$. Grâce au lemme 2.14 appliqué sur ces paires, on déforme $\xi'$ en une structure $\xi''$ dans laquelle chaque tore $F_t$ est totalement réel ou convexe selon que $t$ appartient ou non à l'un des intervalles $[t_i - \delta, t_i + \delta]$. En particulier, la structure $\xi''$ est élémentaire et il reste à comprendre son feuillage.

Le signe de la connexion de $\xi' F_{t_i}$ détermine la direction orientée des feuilletages non singuliers $\xi'' F_t$, $t \in [t_i - \delta, t_i + \delta]$. La classe d'homologie des orbites fermées de $\xi'' F_{t_i - \delta}$ (resp. de $\xi'' F_{t_i + \delta}$) est ainsi $(-p + i - 1, 1)$ (resp. $(-p + i, 1)$) si le signe est $+$ et l'opposé sinon. Or le feuillage de $\xi''$ a un tore compressible dans $F \times [t_i, t_{i+1}]$ si et seulement si $\xi'' F_{t_i + \delta}$ et $\xi'' F_{t_{i+1} - \delta}$ ont des directions opposées (cf. lemme 3.17). Le feuillage de $\xi''$ contient donc autant de tores compressibles qu'il y a de changements de signes dans l'énumération choisie des connexions de $\xi F_{1/2}$.

En outre, le feuillage de $\xi''$ contient deux anneaux qui s'appuient l'un sur $F_0$ et l'autre sur $F_1$. Leurs classes d'isotopie sont déterminées par les signes respectifs de la première et de la dernière connexion du film $\xi' F_t$ (lemme 2.11).

En résumé, on obtient autant de formes normales pour $\xi$ (à $\partial$-isotopie des feuillages persistants près) qu'il y a de mots de $q+p$ lettres en $+$ et $-$, à savoir $\binom{p+q}{q}$.

Dans la démonstration qui suit, on appelle *support* d'une partie $Q$ de $B$ la plus petite réunion connexe d'arêtes de $\partial E$ qui contient $Q$. Le lemme 1.9-b) et le corollaire 1.10 montrent que si $Q$ et $R$ sont deux parties de $B$ telles que $\gamma(Q) = \gamma(R)$, alors elles ont le même support.

*Démonstration de la proposition* 3.40. Compte tenu des propositions 3.15 et 3.22, on prend $\xi$ sous forme normale et, pour fixer les idées, on suppose que son lieu de retournement $R(\xi)$ a un nombre pair de points (autrement dit, les cycles asymptotiques de $\xi F_1$ sont sur $\partial C$). Comme $\gamma(Q) = \gamma\bigl(R(\xi)\bigr)$, l'ensemble $Q$ compte aussi un nombre pair de points et a le même support que $R(\xi)$. On note $v_0 \prec v_1 \prec \cdots \prec v_m$ les sommets de $\partial E$ situés dans ce support et $J_i$, $0 \leq i \leq m$, l'intervalle des $t$ pour lesquels $D_t$ contient $v_i$. Si $t_i$ est un point pris dans $J_i$ assez près du bord, $\xi F_{t_i}$ est une suspension scindée par deux courbes essentielles. Par suite, $F_{t_i}$ rencontre une seule feuille $S_i \subset \mathfrak{F}(\xi)$ et la partage en deux anneaux $S_i^+$ et $S_i^-$ situés respectivement au-dessus et en-dessous de $F_{t_i}$. De plus, quitte à déformer $\xi$ parmi les structures sous forme normale, on peut s'arranger pour que la classe d'homologie des orbites fermées de $\xi F_{t_i}$ soit dans $C$ (lemmes 2.3 et 2.6). On pose, pour tout $1 \leq i \leq m$,

$$V_i = F \times [t_{i-1}, t_i] \quad \text{et} \quad \xi_i = \xi \,|\, V_i \,.$$

Chaque structure $\xi_i$ est ainsi sous forme normale sur $V_i$, son lieu de retournement $R(\xi_i)$ a un nombre pair de points et $\chi_\partial(\xi) = \sum_1^m \chi_\partial(\xi_i)$. En outre, si $S_i$ est un tore incompressible (resp. compressible), $S_i^- \subset V_i$ et $S_i^+ \subset V_{i+1}$ sont des anneaux de même parité (resp. de parités différentes).

**Assertion.** *Il existe des ensembles $Q_i$, $1 \leq i \leq m$, ayant tous un cardinal pair et les propriétés suivantes :*
- *$Q_i \subset [v_{i-1}, v_i]$ ;*
- *$\bigcup_{i=1}^m Q_i \subset Q \cup \{v_0, \ldots, v_m\}$ ;*
- *$2\gamma(Q_i) = \chi_\partial(\xi_i)$.*

*Preuve.* On construit d'abord $Q_1$. Si $Q_1' = Q \cap [v_0, v_1]$ contient un nombre pair de points,

on prend $Q_1 = Q_1'$. Sinon, on pose

$$Q_1 = \begin{cases} Q_1' \setminus \{v_1\} & \text{si } v_1 \in Q, \\ Q_1' \cup \{v_1\} & \text{si } v_1 \notin Q. \end{cases}$$

Il existe alors un seul ensemble $P$ contenant $Q \setminus Q_1$ et inclus dans $(Q \setminus Q_1) \cup \{v_1\}$ tel que $\gamma(Q) = \gamma(Q_1 \cup P)$. On pose $Q_2 = P_1$ où $P_1$ est obtenu à partir de $P$ comme $Q_1$ à partir de $Q$. On construit ainsi des ensembles $Q_i$ ayant chacun un cardinal pair et vérifiant les deux premières propriétés. En outre,

$$\sum_{i=1}^{m} \gamma(Q_i) = \gamma(Q) = \tfrac{1}{2}\chi_\partial(\xi) = \tfrac{1}{2}\sum_{i=1}^{m}\chi_\partial(\xi_i)$$

et il en résulte, d'après le lemme 1.9, que $2\gamma(Q_i) = \chi_\partial(\xi_i)$. □

La discussion de l'exemple 3.41 permet maintenant de finir la preuve de la proposition 3.40. En effet, pour tout $i$, on peut trouver un paramétrage (fibré)

$$\mathbf{T}^2 \times [t_{i-1}, t_i] \longrightarrow V_i = F \times [t_{i-1}, t_i]$$

qui identifie l'arête $[v_{i-1}, v_i] \subset H_1(F; \mathbf{R}) \cong \mathbf{R}^2$ au segment entre des points du type $(-p, 1)$ et $(q, 1)$ où $q = \tfrac{1}{2}\chi_\partial(\xi_i)$. De plus, chaque composante du bord de $(V_i, \xi_i)$ est convexe et a un feuilletage scindé par deux courbes. Ainsi, modulo une isotopie parmi les structures à bord convexe, $\xi_i$ est $\partial$-isotope à la structure étudiée dans l'exemple 3.41. Par suite, $\xi_i$ est $\partial$-isotope à une structure $\xi_i'$ sous forme normale dont le lieu de retournement $R(\xi_i')$ est égal à $Q_i$. En recollant les morceaux, on obtient la structure cherchée. □

## G  Torsion des structures sous forme normale

On calcule ici la torsion des structures de contact sous forme normale.

**Proposition 3.42.** *Soit $\xi$ une structure de contact tendue sous forme normale sur $V$ et $J$ l'intervalle*

$$J = [0, 1] \setminus \operatorname{Int} f\bigl(\mathfrak{F}_\partial(\xi)\bigr).$$

**a)** *Si $\mathfrak{F}_\perp(\xi)$ ou $\mathfrak{F}_o(\xi)$ est non vide, la torsion de $\xi$ est nulle.*
**b)** *Si $\mathfrak{F}_\perp(\xi)$ et $\mathfrak{F}_o(\xi)$ sont vides, la $\pi$-torsion de $\xi$ est le plus grand entier strictement inférieur à $-\theta/\pi$, où $\theta$ désigne l'amplitude — non nulle — de la restriction de $\xi$ à $F \times J$ — qui est bien une structure rotative.*

*Démonstration.* si $\mathfrak{F}_\perp(\xi)$ n'est pas vide, $\mathfrak{F}_o(\xi)$ est vide car $\xi$ est tendue (lemme 3.34). Le lemme 3.31 assure alors que $(V, \xi)$ se plonge dans $(F \times J, \xi)$. Or la proposition 8 de [Gi4] montre que $\mathbf{T}^2 \times [0, 1]$ muni d'une structure rotative d'amplitude non nulle n'admet aucun plongement incompressible dans $(F \times J, \xi)$. Par suite, la torsion de $\xi$ est nulle.

Si $\mathfrak{F}_o(\xi)$ est vide, $\xi$ est virtuellement vrillée et sa torsion est nulle en vertu de la proposition 3.36 et de l'argument développé dans le b) ci-dessous.
**b)** Dans un premier temps, on suppose que $\mathfrak{F}_\partial(\xi)$ est vide (tout comme $\mathfrak{F}_\perp(\xi)$ et $\mathfrak{F}_o(\xi)$), et donc que $J = [0, 1]$. Le théorème 3.3 et le lemme 3.11 (démonstration incluse) montrent alors que la torsion de $\xi$ est au moins égale au plus grand entier strictement inférieur à $-\theta/\pi$. On se donne donc maintenant un plongement

$$\Bigl(\mathbf{T}^2 \times [0, 1], \ker\bigl(\cos(n\pi z)\, dx - \sin(n\pi z)\, dy\bigr)\Bigr) \longrightarrow (V, \xi), \qquad (x, y, z) \in \mathbf{T}^2 \times [0, 1],$$

et on note $W$ son image. D'après un théorème de J. STALLINGS, $V \setminus \operatorname{Int} W$ a deux composantes connexes, $V_0 \supset F_0$ et $V_1 \supset F_1$, toutes deux difféomorphes à $\mathbf{T}^2 \times [0,1]$. D'autre part, la restriction $\xi_i$ de $\xi$ à $V_i$, $i \in \{0,1\}$, est une structure de contact tendue et admise. La proposition 3.15 dit par suite que, pour un paramétrage bien choisi de $V_i$ par $\mathbf{T}^2 \times [0,1]$, la structure $\xi_i$ est élémentaire. Comme elle est universellement tendue et imprime sur le bord de $V_i$ des feuilletages linéarisables (car $\mathfrak{F}_\partial(\xi)$ est vide), les ensembles $\mathfrak{F}_o(\xi_i)$, $\mathfrak{F}_\perp(\xi_i)$ et $\mathfrak{F}_\partial(\xi_i)$ sont vides et $\xi_i$ est donc une structure rotative. Il découle alors immédiatement du corollaire 3.10 que $n\pi < -\theta$, ce qu'on voulait démontrer. Enfin, si $\mathfrak{F}_\partial(\xi)$ n'est pas vide, on conclut en appliquant l'argument qui précède à la structure $\xi'$ fournie par le lemme ci-dessous. $\square$

**Lemme 3.43.** *On reprend les notations de la proposition 3.42 et on suppose que $\mathfrak{F}_\perp(\xi)$ et $\mathfrak{F}_o(\xi)$ sont vides. Pour tout $\varepsilon > 0$, il existe sur $V' = \mathbf{T}^2 \times [0,1]$ une structure rotative $\xi'$ d'amplitude comprise entre $\theta - \varepsilon$ et $\theta$ qui a les propriétés suivantes :*
  – *$(V, \xi)$ se plonge incompressiblement dans $(V', \xi')$ ;*
  – *$\mathfrak{F}_\partial(\xi')$ est vide.*

*Démonstration.* On raisonne par récurrence sur la somme $k = k_0 + k_1$ où $k_i$ est le nombre de composantes connexes de $\mathfrak{F}_\partial(\xi)$ qui s'appuient sur $F_i$. Si $k = 0$, il est clair que $\xi' = \xi$ convient.

Si $k_0 > 1$, on prend dans $\mathfrak{F}_\partial(\xi)$ une composante $A$ de plus petite hauteur $\sup f \mid_A$ parmi tous les anneaux qui s'appuient sur $F_0$. Le lemme 3.31 donne une structure élémentaire $\xi_-$ sur $F \times [-1, 0]$ qui imprime $\xi F_0$ sur $F_0$ et dont le feuillage est formé d'un anneau reliant les orbites $A$ à un autre anneau de $\mathfrak{F}_\partial(\xi)$ et de $2k_0 - 2$ anneaux allant de $F_0$ à $F_{-1}$. En collant $\xi_-$ avec $\xi$, on obtient sur $F \times [-1, 1]$ une structure élémentaire dont le feuillage a un anneau de moins parallèle au bord inférieur. De plus, le fait que la hauteur de $A$ soit minimale permet de mettre la structure obtenue sous forme normale par des isotopies dont le support est disjoint de $F \times J$. On peut donc répéter l'argument et procéder de même sur le bord supérieur pour se ramener au cas où $k_0$ et $k_1$ valent au plus 1.

Si $k_0 = 1$, on note $A$ la composante de $\mathfrak{F}_\partial(\xi)$ qui s'appuie sur $F_0$. On construit alors sans peine sur $F \times [-1, 0]$ une structure élémentaire $\xi_-$ ayant les propriétés suivantes :
  – $\xi_-$ imprime $\xi F_0$ sur $F_0$ ;
  – le feuillage persistant $\mathfrak{F}'(\xi_-)$ est constitué d'un seul anneau $A_-$ dont la réunion avec $A$ dans $F \times [-1, 1]$ forme un tore isotope aux fibres $F_t$ ;
  – la restriction de $\xi_-$ à $F \times J_-$, où $J_- = [\inf f \mid_{A_-}, 0]$, est une structure rotative d'amplitude supérieure à $-\varepsilon/2$.

Pour trouver la structure rotative voulue, il suffit alors de coller $\xi_-$ sur $\xi$, de mettre la structure obtenue sous forme normale puis de procéder de même sur le bord supérieur si $k_1 = 1$. $\square$

# 4 Démonstrations des résultats de classification

## A Les espaces lenticulaires

*Démonstration du théorème 1.1.* Soit $\pi \colon \mathbf{T}^2 \times [0,1] \to \mathbf{L}_{p,q}$ la projection dont les fibres non ponctuelles sont les cercles de direction $(0,1)$ sur $\mathbf{T}^2 \times \{0\}$ et ceux de direction $(p,q)$ sur $\mathbf{T}^2 \times \{1\}$. Pour tout $t \in [0,1]$, on note $F_t$ l'image de $\mathbf{T}^2 \times \{t\}$ dans $\mathbf{L}_{p,q}$ ; c'est un tore, sauf pour $t \in \{0,1\}$ où c'est une courbe qu'on oriente comme l'image de $\mathbf{S}^1 \times \{0\} \times \{t\}$.

On dira qu'une structure de contact orientée $\xi$ sur $\mathbf{L}_{p,q}$ est sous forme normale si :

1) $F_0$ et $F_1$ sont des transversales positives de $\xi$ ;

2) la structure de contact $\widetilde{\xi}$ induite par $\xi$ sur $\mathbf{T}^2 \times [0,1]$ est sous forme normale.

La condition 1) implique que les caractéristiques de $\widetilde{\xi}$ sur $\mathbf{T}^2 \times \{0\}$ et $\mathbf{T}^2 \times \{1\}$ sont les cercles orientés de directions respectives $(0,-1)$ et $(-p,-q)$. Par suite, le lieu de retournement $R(\xi) = R(\widetilde{\xi})$ est de cardinal pair.

On pose alors
$$C = \{(x,y) \in \mathbf{R}^2 \mid x < 0,\ py < qx\}$$
et on note $E$ l'enveloppe convexe des points entiers de $C$.

Le théorème 1.1 découle directement de la proposition 1.7 et du lemme ci-dessous.

**Lemme 4.1.** *Soit $\widehat{B}$ l'ensemble des points entiers primitifs de $\partial E$ situés sur une arête de longueur finie.*
**a)** *Toute structure de contact tendue sur $\mathbf{L}_{p,q}$ est isotope à une structure sous forme normale dont le lieu de retournement est inclus dans $\widehat{B}$.*
**b)** *Soit $\xi$ une structure de contact sous forme normale dont le lieu de retournement $R(\xi)$ est inclus dans $\widehat{B}$. La classe d'isotopie de $\xi$ est déterminée par $\gamma(R(\xi))$. De plus, la structure $\xi$ est universellement tendue si et seulement si $R(\xi)$ est inclus dans $\partial \widehat{B}$.*
**c)** *Sur l'ensemble des parties de cardinal pair dans $\widehat{B}$, la fonction $\gamma$ prend un nombre de valeurs égal à $\prod_{i=0}^{n}(a_i - 1)$.*

*Démonstration.*
**a)** Toute structure de contact tendue sur $\mathbf{L}_{p,q}$ est isotope à une structure $\xi$ sous forme normale (propositions 3.15 et 3.22). Soit maintenant $A_0$ l'arête extrême de $\partial E$ parallèle à $D_0$. La proposition 3.40 montre que $\xi$ est isotope, relativement à $F_0 \cup F_1$, à une structure $\xi'$ dont le lieu de retournement a un nombre pair de points dans l'intérieur de $A_0$ (si $R(\xi) \cap A_0$ compte un nombre impair de points, $R(\xi')$ contient l'extrémité de $A_0$). On choisit alors un point $t_0 \in [0,1]$ pour lequel la direction asymptotique de $\widetilde{\xi}' F_{t_0}$ est une droite de pente comprise entre 1 et 2 et on pose $V_0 = \pi(\mathbf{T}^2 \times [0,t_0]) \subset \mathbf{L}_{p,q}$. La restriction de $\xi'$ à $V_0$ est une structure de contact tendue dont la classe d'Euler relative est le générateur de $H_1(V_0; \mathbf{Z})$. Dès lors, un résultat de S. MAKAR-LIMANOV [Ma] (directement issu de la classification des structures de contact sur la boule $\mathbf{D}^3$ [El4]) assure que $\xi'$ est $\partial$-isotope, sur $V_0$, à « la structure de contact ordinaire » qui est transversale à $F_0$ et imprime sur tous les tores $F_t$, $0 < t \leq t_0$, une suspension. En remplaçant $\xi'$ dans $V_0$ par cette structure ordinaire, on obtient sur $\mathbf{L}_{p,q}$ une structure de contact $\xi''$ sous forme normale dont le lieu de retournement évite l'intérieur de $A_0$. Le même argument avec l'autre arête extrême de $\partial E$ donne la structure voulue.

**b)** Si $R(\xi)$ est non vide, la première assertion vient de la proposition 3.40. Si $R(\xi)$ est vide, $\widetilde{\xi}$ est une structure rotative. De plus, la direction asymptotique $D_t$ de $\widetilde{\xi} F_t$ ne peut coïncider avec $D_0$ ou $D_1$ pour $t \in \,]0,1[$, sans quoi toute orbite fermée de $\xi F_t$ borderait un disque vrillé dans $\mathbf{L}_{p,q}$. Par suite, l'amplitude de $\widetilde{\xi}$ est supérieure à $-\pi$ et la première assertion découle du théorème 3.3.

D'autre par, le revêtement
$$\rho \colon \begin{cases} \mathbf{T}^2 \times [0,1] \longrightarrow \mathbf{T}^2 \times [0,1] \\ (x,y,t) \longmapsto (px, qx+y, t) \end{cases}$$
induit un revêtement de $\mathbf{S}^3 = \mathbf{L}_{1,0}$ sur $\mathbf{L}_{p,q}$ qui envoie les cercles de Hopf de direction $(1,-1)$ sur ceux de direction $(p, q-1)$. Si $R(\xi)$ est vide, $\widetilde{\xi}$ est transversale à ces derniers cercles ; $\xi$ se relève alors à $\mathbf{S}^3$ en une structure transversale à la fibration de Hopf, donc isotope à la structure tendue ordinaire (théorème 3.3). Si $R(\xi)$ se réduit aux extrémités de $\widehat{B}$, le résultat de S. MAKAR-LIMANOV invoqué plus haut montre que $-\xi$ est isotope à une structure $\xi'$ sous forme normale dont le lieu de retournement est vide. Enfin, pour toute structure $\xi$ sous forme normale dont le lieu de retournement n'est ni vide, ni égal à $\partial\widehat{B}$, le relèvement à $\mathbf{S}^3 = \mathbf{L}_{1,0}$ est vrillé d'après la proposition 3.36.

**c)** Cette assertion provient directement du corollaire 1.10. $\square$

## B  Le tore plein

*Démonstration du théorème 1.6.* On pose $V = \mathbf{S}^1 \times \mathbf{D}^2$ et $F_t = \mathbf{S}^1 \times t\mathbf{S}^1$ où $t\mathbf{S}^1$, pour $t \in [0, 1]$, est le cercle de rayon $t$. On note aussi $i$ l'inclusion de $\partial V = \mathbf{T}^2$ dans $V = \mathbf{S}^1 \times \mathbf{D}^2$ et on observe que l'application $i_* \colon H_1(\partial V; \mathbf{Z}) \to H_1(V; \mathbf{Z})$ s'identifie à la projection de $\mathbf{Z}^2$ sur son premier facteur. On rappelle d'autre part que $C = C(\sigma_0, \sigma) \subset \mathbf{R}^2 \setminus \{0\}$ est ici le cône bordé à gauche par la demi-droite $\mathbf{R}_+(0, -1)$ et ouvert de ce côté. Enfin, on oriente $\sigma$ de telle sorte que ses cycles asymptotiques soient sur le bord droit de $C$ et ce choix fixe une orientation sur toute structure de contact dans $\mathcal{SCT}(V; \sigma)$.

On dira qu'une structure de contact $\xi \in \mathcal{SCT}(V; \sigma)$ est sous forme normale si :

1) $F_0$ est une transversale positive de $\xi$ ;
2) la structure de contact $\widetilde{\xi}$ induite par $\xi$ sur $\mathbf{T}^2 \times [0, 1]$ est sous forme normale.

La condition 1) implique que les caractéristiques de $\widetilde{\xi}$ sur $\mathbf{T}^2 \times \{0\}$ sont les cercles orientés de direction $(0, -1)$. Par suite, le lieu de retournement $R(\xi) = R(\widetilde{\xi})$ est de cardinal pair.

Le théorème 1.6 est une conséquence directe du lemme ci-dessous.

**Lemme 4.2.** *Soit $\widehat{B}$ l'ensemble des points entiers primitifs de $\partial E$ situés sur une arête non verticale.*
**a)** *Soit $\xi$ une structure de contact sous forme normale sur $V$ qui imprime $\sigma$ sur $F_1 = \partial V$. Si le lieu de retournement $R(\xi)$ est inclus dans $\partial E$ et si la direction asymptotique de $\widetilde{\xi}F_t$ rencontre $C$ pour tout $t \in\, ]0, 1[$, la structure $\xi$ est tendue.*
**b)** *Toute structure de contact tendue dans $\mathcal{SCT}(V; \sigma)$ est $\partial$-isotope à une structure sous forme normale dont le lieu de retournement est inclus dans $\widehat{B}$.*
**c)** *Soit $\xi$ une structure de contact sous forme normale dont le lieu de retournement $R(\xi)$ est inclus dans $\widehat{B}$. La classe d'Euler relative de $\xi$ vaut $1 + 2i_*\gamma\bigl(R(\xi)\bigr)$ et détermine entièrement la classe de $\partial$-isotopie. De plus, la structure $\xi$ est universellement tendue si et seulement si $R(\xi)$ est inclus dans $\partial\widehat{B}$.*

*Démonstration.*
**a)** Soit $[u, v]$, $u \wedge v < 0$, l'arête de $\partial E$ qui contient l'ultime point de $R(\xi)$ du côté de $D = D(\sigma)$ et soit $v'$ le premier point entier situé au-delà de $v$ sur le prolongement de cette arête. On note $C'$ le cône ouvert bordé à gauche par $D_0$, à droite par $D' = \mathbf{R}v'$ et qui contient $C$. Par construction, l'ensemble $R(\xi)$ est inclus dans le bord de l'enveloppe convexe $E'$ des points entiers de $C'$. On écrit alors $v' = -(p, mp + q)$, où $p$ et $q$ sont des entiers premiers entre eux tels que $0 < q < p$. D'après la proposition 1.7 et le lemme 4.1, il existe sur $\mathbf{L}_{p,q}$ une structure de contact tendue $\xi'$ sous forme normale dont le lieu de retournement $R(\xi')$ est l'image de $R(\xi)$ par la transformation $(x, y) \mapsto (x, -mx + y)$. Le fait que la direction asymptotique du feuilletage $\widetilde{\xi}F_t$ reste dans $C$ pour tout $t \in\, ]0, 1[$ garantit alors que $(V, \xi)$ se plonge dans $(\mathbf{L}_{p,q}, \xi')$, de sorte que $\xi$ est tendue.

La démonstration des parties b) et c) est identique à celle des parties a) et b) du lemme 4.1. □

## C  Le tore épais

*Démonstration du théorème 1.5.* On pose $V = T^2 \times [0,1]$ et $F_t = \mathbf{T}^2 \times \{t\}$, $t \in [0,1]$.

Soit $\xi \in \mathcal{SCT}(V; \sigma_0, \sigma_1)$. Pour calculer la $\pi$-torsion et la classe d'Euler relative, on suppose $\xi$ sous forme normale (propositions 3.15 et 3.22). Si $\xi$ est virtuellement vrillée, $\mathfrak{F}_o(\xi)$ est non vide (proposition 3.29). Par suite, $R(\xi)$ rencontre Int $B$ et $\chi_\partial(\xi)$ appartient à $X_v$ (lemme 3.39). En outre, la $\pi$-torsion de $\xi$ est nulle (proposition 3.42). Si $\xi$ est universellement tendue, $\mathfrak{F}_o(\xi)$ est vide. Le lieu $R(\xi)$ est alors inclus dans $\partial B$ car les seuls retournements possibles viennent des éventuels anneaux impairs dans $\mathfrak{F}_\partial(\xi)$ ou $\mathfrak{F}_\perp(\xi)$. ainsi, $\chi_\partial(\xi)$ appartient à $X_u$.

Pour tout $w \in X_v$, la fibre $\chi_\partial^{-1}(w)$ est connexe d'après la proposition 3.40. Elle est d'autre part non vide en vertu du lemme suivant qui découle de la proposition 1.7 exactement comme le lemme 4.2-a) :

**Lemme 4.3.** *Soit $\xi$ une structure de contact sous forme normale sur $V$ qui imprime $\sigma_0$ sur $F_0$ et $\sigma_1$ sur $F_1$. Si le lieu de retournement $R(\xi)$ est inclus dans $\partial E$ et si la direction asymptotique de chaque feuilletage $\xi F_t$ rencontre $C$, la structure $\xi$ est tendue.* □

Pour tout $w \in X_u$ et tout $n \in \mathbf{N}$, la proposition 3.42 et le lemme 3.11 fournissent une structure de contact $\xi$ sous forme normale dont la classe d'Euler relative vaut $w$, la $\pi$-torsion $n$ et dont le feuillage ne contient aucun anneau allant d'un bord à l'autre ($\mathfrak{F}_\perp(\xi) = \varnothing$). Le corollaire 3.28 montre de plus que cette structure est unique à $\partial$-isotopie près. Par suite, la fibre $(\chi_\partial \times \tau_\pi)^{-1}(w, n)$ est connexe sauf éventuellement si elle contient une structure de contact $\xi$ sous forme normale pour laquelle $\mathfrak{F}_\perp(\xi)$ est non vide. Cela suppose que les conditions suivantes soient remplies :
  – $\sigma_0$ et $\sigma_1$ sont tous deux de type 1 et ont la même direction asymptotique ;
  – la $\pi$-torsion $n$ est nulle (proposition 3.42) ;
  – $w = \chi_\partial(\xi)$ est égal à 0 ou au double du vecteur directeur primitif de $D_0$, selon que les anneaux de $\mathfrak{F}_\perp(\xi)$ sont pairs ou impairs. (Comme $\sigma_0$ et $\sigma_1$ sont de type 1, *i.e.* scindés par deux courbes, la parité des anneaux de $\mathfrak{F}_\perp(\xi)$ détermine $\mathfrak{F}_\perp(\xi)$ à l'action près des difféomorphismes relatifs au bord.)

Dans ces conditions, le corollaire 3.27 dit que les structures $\xi \in (\chi_\partial \times \tau_\pi)^{-1}(w, n)$ telles que $\mathfrak{F}_\perp(\xi)$ soit non vide forment une seule orbite sous l'action des difféomorphismes relatifs au bord. Le fait qu'elles forment une infinité de composantes connexes résulte de la proposition 4.6. □

Avant de compléter la démonstration ci-dessus, on donne un énoncé plus général pour les structures de contact universellement tendues sur le tore épais. Ce théorème montre que, pour des suspensions $\sigma_0$, $\sigma_1$ scindées par plus de deux courbes, la classe d'Euler relative ne suffit pas à distinguer les structures de contact de torsion fixée : il faut prendre en compte la combinatoire des systèmes d'anneaux du feuilletage.

**Théorème 4.4.** *Soit $\sigma_0$ et $\sigma_1$ des suspensions admissibles de $\mathbf{T}^2$ scindées respectivement par $2n_0$ et $2n_1$ courbes essentielles. La torsion définit une application*

$$\tau_\pi \colon \mathcal{SCT}_u(V; \sigma_0, \sigma_1) \longrightarrow \mathbf{N}$$

*qui est surjective et dont les fibres ont toutes exactement $\binom{2n_0}{n_0}\binom{2n_1}{n_1}$ composantes connexes, sauf dans un cas : lorsque les directions asymptotiques de $\sigma_0$, $\sigma_1$ coïncident et que les* entiers $n_0$, $n_1$ sont tous deux non nuls, la fibre $\tau_\pi^{-1}(0)$ contient une infinité de composantes connexes supplémentaires qui, sous l'action des difféomorphismes relatifs au bord, se rangent en un nombre fini d'orbites égal à

$$\sum_{n=1}^{\inf(n_0,n_1)} 2n \binom{2n_0}{n_0-n}\binom{2n_1}{n_1-n}.$$

On élucide d'abord la présence des coefficients binomiaux.

**Lemme 4.5.**
**a)** *On se donne $2n$ points sur le bord du disque $\mathbf{D}^2$. À $\partial$-isotopie près, il existe $\binom{2n}{n}$ familles d'arcs disjoints dans $\mathbf{D}^2 \setminus \{0\}$ qui relient ces points deux par deux.*
**b)** *On se donne $2(n_0+n_1)$ points sur le bord de l'anneau $\mathbf{S}^1 \times [0,1]$, à savoir $2n_i$ points sur $\mathbf{S}^1 \times \{i\}$, $i \in \{0,1\}$. À l'action près des difféomorphismes égaux à l'identité sur le bord, les familles d'arcs disjoints dans $\mathbf{S}^1 \times [0,1]$ qui relient ces points deux par deux et contiennent au moins un arc allant d'un bord à l'autre sont en nombre*

$$\sum_{n=1}^{\inf(n_0,n_1)} 2n \binom{2n_0}{n_0-n}\binom{2n_1}{n_1-n}.$$

*Démonstration.* **a)** Tout arc dans $\mathbf{D}^2 \setminus \{0\}$ qui joint deux points de $\partial \mathbf{D}^2$ a une *origine* définie comme l'extrémité de plus petit argument. Quand on choisit $n$ points parmi les $2n$ donnés, il y a une unique possibilité (à $\partial$-isotopie près) pour tracer dans $\mathbf{D}^2 \setminus \{0\}$ des arcs disjoints qui partent de ces points et aboutissent aux $n$ autres points. En fait, si un point choisi est suivi d'un point non choisi sur le cercle orienté trigonométriquement, on est forcé de les relier ; on efface alors cette paire et on recommence.

**b)** Soit $n \in \{1, \ldots, \inf(n_0, n_1)\}$. Quand on choisit, parmi les points donnés, $n_0 - n$ points sur $\mathbf{S}^1 \times \{0\}$ (resp. $n_1 - n$ points sur $\mathbf{S}^1 \times \{1\}$), il y a une unique possibilité pour tracer dans $\mathbf{S}^1 \times [0,1]$ des arcs disjoints qui partent de ces points, aboutissent à d'autres points donnés sur $\mathbf{S}^1 \times \{0\}$ (resp. sur $\mathbf{S}^1 \times \{1\}$) et n'emprisonnent pas les $2n$ points restants. On a alors $2n$ possibilités (à l'action près des difféomorphismes relatifs au bord) pour relier un à un, par des arcs disjoints, ces $2n$ points restants aux $2n$ points libres sur le bord d'en face. □

*Démonstration du théorème 4.4.* Chaque structure de contact dans $\mathcal{SCT}_u(V;\sigma_0,\sigma_1)$ est $\partial$-isotope à une structure de contact $\xi$ sous forme normale (propositions 3.15 et 3.22). En outre, la classe de $\partial$-isotopie de $\xi$ est déterminée par sa $\pi$-torsion et la combinatoire du système d'anneaux $\mathfrak{F}_\perp(\xi) \cup \mathfrak{F}_\partial(\xi)$ (corollaires 3.27, 3.28 et proposition 3.42). Le lemme 4.5 dit que le nombre de combinatoires possibles pour les anneaux est exactement le nombre de composantes connexes annoncé pour $\tau_\pi^{-1}\big(\tau_\pi(\xi)\big)$. Comme toutes ces combinatoires sont réalisables (lemme 3.31), la démonstration du théorème se réduit à celle de la proposition 4.6 ci-dessous. □

**Proposition 4.6.** *Soit $\xi_0$, $\xi_1 \in \mathcal{SCT}_u(V;\sigma_0,\sigma_1)$ des structures de contact sous forme normale. Si $\mathfrak{F}_\perp(\xi_0) \cup \mathfrak{F}_\partial(\xi_0)$ et $\mathfrak{F}_\perp(\xi_1) \cup \mathfrak{F}_\partial(\xi_1)$ ne sont pas $\partial$-isotopes, les structures $\xi_0$ et $\xi_1$ ne le sont pas non plus.*

La preuve repose en partie sur l'inégalité de Bennequin par le biais du lemme suivant.

**Lemme 4.7.** *Soit $\eta$ une structure de contact sur $\mathbf{T}^2 \times \mathbf{R}$ pour laquelle tous les tores $\mathbf{T}^2 \times \{t\}$ sont convexes et ont un feuilletage caractéristique scindé par $2n$ courbes isotopes* à $\{x\} \times \mathbf{S}^1 \times \{t\}$, où $n \geq 1$. Si $L$ est une courbe legendrienne fermée simple tracée sur un tore $T \subset \mathbf{T}^2 \times \mathbf{R}$ isotope à $\mathbf{T}^2 \times \{0\}$, le nombre de tours que fait la structure $\eta$ par rapport au plan tangent à $T$ le long de $L$ vérifie

$$\deg(\eta, T; L) \leq -n \Big| \int_L dx \Big|.$$

*Démonstration.* On oriente $L$ de telle sorte que l'entier $p = \int_L dx$ soit positif ou nul et on pose $q = \int_L dy$. Comme $L$ est tracée sur un tore $T$ isotope à $\mathbf{T}^2 \times \{0\}$, les entiers $p$ et $q$ sont premiers entre eux. On distingue alors deux cas.

Si $p = 0$, on note $\zeta$ la structure de contact ordinaire sur $\mathbf{S}^3$ — champ des tangentes complexes à la sphère unité de $\mathbf{C}^2$ — et $S_0$ le tore de $\mathbf{S}^3$ d'équation $|z| = |w|$, $(z, w) \in \mathbf{C}^2$. On voit sans peine que $S_0$ est un tore pré-lagrangien, que ses caractéristiques ne sont pas nouées et que leur invariant de Thurston-Bennequin vaut $-1$. Par une isotopie $\mathcal{C}^\infty$-petite, on perturbe $S_0$ en un tore $S$ qui intersecte $S_0$ transversalement le long d'exactement $2n$ caractéristiques de $S_0$. Ce tore $S$ est alors convexe et son feuilletage $\zeta S$ est scindé par $2n$ courbes essentielles. Vu le lemme 2.8, il existe un difféomorphisme $\psi$ de $\mathbf{T}^2 \times \mathbf{R}$ sur un voisinage tubulaire $U$ de $S$ qui pousse $\eta$ sur $\zeta$ et les courbes $\{x\} \times \mathbf{S}^1 \times \{t\}$ sur des courbes isotopes dans $U$ aux composantes de $S \cap S_0$. La courbe legendrienne $\psi(L)$ est ainsi non nouée dans $\mathbf{S}^3$ et son invariant de Thurston-Bennequin vaut $-1 + \deg(\eta, T; L)$. L'inégalité de Bennequin assure donc que $\deg(\eta, T; L) \leq 0$.

Si $p > 0$, on considère le revêtement

$$\rho \colon \mathbf{T}^2 \times \mathbf{R} \longrightarrow \mathbf{T}^2 \times \mathbf{R}, \quad (x, y, t) \longmapsto (px, y + qx, t).$$

Tous les tores $\mathbf{T}^2 \times \{t\}$ sont convexes pour $\widehat{\eta} = \rho^*\eta$ et ont un feuilletage caractéristique scindé par $2np$ courbes isotopes à $\{x\} \times \mathbf{S}^1 \times \{t\}$. En outre, le tore $\widehat{T} = \rho^{-1}(T)$ est isotope à $\mathbf{T}^2 \times \{0\}$ et chaque composante $\widehat{L}$ de $\rho^{-1}(L)$ vérifie

$$\begin{aligned} \int_{\widehat{L}} dx &= 1 \\ \int_{\widehat{L}} dy &= 0 \end{aligned} \quad \text{et} \quad \deg(\widehat{\eta}, \widehat{T}; \widehat{L}) = \deg(\eta, T; L).$$

Il suffit ainsi de traiter le cas où $p = 1$ et $q = 0$. En outre, d'après le lemme 2.8, on peut supposer que $\eta$ a pour équation

$$\cos(2n\pi x)\, dy - \sin(2n\pi x)\, dt = 0.$$

Soit alors $\widetilde{\eta}$ (resp. $\widetilde{L}$) le rappel de $\eta$ (resp. une préimage de $L$) par le revêtement évident $\mathbf{S}^1 \times \mathbf{R} \times \mathbf{R} \to \mathbf{T}^2 \times \mathbf{R}$. Puisque $\widetilde{\eta}$ est invariante par toutes les transformations $(x, y, t) \mapsto (x, \lambda y, \lambda t)$ et que la courbe $\mathbf{S}^1 \times \{(0, 0)\}$ est legendrienne, on peut plonger $(\mathbf{S}^1 \times \mathbf{R}^2, \widetilde{\eta})$ dans toute variété de contact comme un voisinage tubulaire d'une courbe legendrienne quelconque. Si on effectue ce plongement dans $\mathbf{S}^3$ autour d'une courbe legendrienne non nouée dont l'invariant de Thurston-Bennequin vaut $-1$, l'inégalité de Bennequin fournit immédiatement l'estimation voulue. $\square$

*Démonstration de la proposition* 4.6. On note d'abord que, si $\mathfrak{F}_\partial(\xi_0)$ et $\mathfrak{F}_\partial(\xi_1)$ sont vides, le lemme 2.6 et l'exemple 3.4 prouvent la proposition. On observe ensuite que, si $\mathfrak{F}_\perp(\xi_0)$ a moins de composantes que $\mathfrak{F}_\perp(\xi_1)$, il n'y a aucun plongement incompressible de $(V, \xi_0)$ dans $(V, \xi_1)$. En effet, si $\mathfrak{F}_\perp(\xi_i)$ a $2m_i$ composantes, le lemme 3.31 fournit un plongement incompressible de $(V, \xi_1)$ dans $\mathbf{T}^2 \times \mathbf{R}$ équipé d'une structure de contact $\eta$ pour laquelle les tores $\mathbf{T}^2 \times \{*\}$ sont convexes et ont un feuilletage scindé par $2m_1$ courbes. D'autre part, selon que $m_0$ est nul ou non, l'intervalle $J_0 = [0,1] \setminus \operatorname{Int} f(\mathfrak{F}_\partial(\xi_0))$ est une séquence rotative ou formé de paramètres $t$ pour lesquels $F_t$ est un tore convexe muni d'un feuilletage scindé par $2m_0$ courbes. Or le lemme 4.7 empêche de plonger incompressiblement dans $(\mathbf{T}^2 \times \mathbf{R}, \eta)$ un tel tore avec $m_0 < m_1$ ou une structure rotative d'amplitude non nulle.

On procède maintenant par récurrence sur l'entier $n_0 + n_1$. L'idée est de construire des structures de contact élémentaires $\xi'_0$ et $\xi'_1$ sur $F \times [-1, 2]$ qui aient les propriétés suivantes :

1) $\xi'_0$ et $\xi'_1$ induisent respectivement $\xi_0$ et $\xi_1$ sur $F \times [0, 1]$ ;
2) $\xi'_0$ et $\xi'_1$ coïncident sur $F \times ([-1, 0] \cup [1, 2])$ ;
3) $\xi'_0$ et $\xi'_1$ ne sont pas $\partial$-isotopes — soit par hypothèse de récurrence, soit parce que l'une est vrillée ou virtuellement vrillée mais pas l'autre.

Pour construire ces structures, on va distinguer plusieurs cas. Auparavant, on observe que les familles d'anneaux $\mathfrak{F}_\perp(\xi_0) \cup \mathfrak{F}_\partial(\xi_0)$ et $\mathfrak{F}_\perp(\xi_1) \cup \mathfrak{F}_\partial(\xi_1)$ sont $\partial$-isotopes si et seulement si chaque anneau de la première famille est individuellement $\partial$-isotope à un anneau de la seconde.

**1)** On suppose d'abord que $n_0 = 1$ et que les ensembles $\mathfrak{F}_\perp(\xi_i)$ sont vides. Dans chacun des ensembles $\mathfrak{F}_\partial(\xi_i)$, $i \in \{0, 1\}$, il y a donc une seule composante connexe $A_i$ qui s'appuie sur $F_0$. Il est facile de construire sur $F \times [-1, 0]$ une structure élémentaire $\xi$ sous forme normale qui imprime $\sigma_0$ sur $F_0$ et dont le feuillage $\mathfrak{F}(\xi) \setminus \mathfrak{F}_\|(\xi)$ se réduise à un anneau $A$ — s'appuyant sur $F_0$ — tel que $A_0 \cup A$ soit un tore incompressible dans $F \times [-1, 1]$. On note alors $\xi'_0$ et $\xi'_1$ les structures obtenues sur $F \times [-1, 1]$ en collant $\xi$ sur $\xi_0$ et $\xi_1$ respectivement. Si $A_0$ et $A_1$ sont isotopes, $\xi'_0$ et $\xi'_1$ sont universellement tendues mais les ensembles $\mathfrak{F}_\partial(\xi_i)$ ont une composante connexe de moins. Si $A_0$ et $A_1$ ne sont pas isotopes, le tore $A_1 \cup A$ est compressible. Du coup, $\xi'_1$ est virtuellement vrillée tandis que $\xi'_0$ demeure universellement tendue.

**2)** On suppose maintenant que l'une des composantes de $\mathfrak{F}_\perp(\xi_0) \cup \mathfrak{F}_\partial(\xi_0)$ est un anneau $A_0$ isotope, relativement à son bord, à un anneau $A_1$ de $\mathfrak{F}_\perp(\xi_1) \cup \mathfrak{F}_\partial(\xi_1)$. On note $C$ l'une des composantes de $\partial A_0 = \partial A_1$ et, pour fixer les idées, on imagine $C$ contenue dans $F_0$. Si la situation n'est pas celle envisagée en 1), l'une des orbites fermées $C'$ de $\sigma_0$ voisines de $C$ ne fait pas partie de $\partial A_0$. On se donne alors sur $F \times [-1, 0]$ une structure élémentaire $\xi$ sous forme normale qui imprime $\sigma_0$ sur $F_0$ et dont le feuillage $\mathfrak{F}_\partial(\xi)$ se compose d'un seul anneau $A$ de bord $\partial A = C \cup C'$. Les structures $\xi'_0$ et $\xi'_1$ obtenues en collant $\xi$ respectivement sur $\xi_0$ et $\xi_1$ restent universellement tendues. En outre, les ensembles $\mathfrak{F}_\perp(\xi'_0) \cup \mathfrak{F}_\partial(\xi'_0)$ et $\mathfrak{F}_\perp(\xi'_1) \cup \mathfrak{F}_\partial(\xi'_1)$ demeurent non isotopes mais ont chacun une composante connexe de moins.

**3)** On suppose désormais, compte tenu de 2), qu'aucune composante de $\mathfrak{F}_\perp(\xi_0) \cup \mathfrak{F}_\partial(\xi_0)$ n'est $\partial$-isotope, à une composante de $\mathfrak{F}_\perp(\xi_1) \cup \mathfrak{F}_\partial(\xi_1)$. On choisit, dans $\mathfrak{F}_\partial(\xi_1)$, un anneau $A_1$ de « profondeur minimale » au sens où il borde, avec un anneau de $F \times \{0, 1\}$, un tore plein disjoint du feuillage $\mathfrak{F}(\xi_0)$. Pour fixer les idées, on imagine toujours $\partial A_1$ contenu dans $F_0$. On prend alors sur $F \times [-1, 0]$ une structure élémentaire $\xi$ sous forme normale qui imprime $\sigma_0$ sur $F_0$ et dont le feuillage $\mathfrak{F}_\partial(\xi)$ se compose d'un seul anneau $A$ tel que $A_1 \cup A$ soit un tore compressible dans $F \times [-1, 1]$. En collant $\xi$ sur $\xi_1$, on obtient une structure élémentaire $\xi'_1$ qui est virtuellement vrillée — et même vrillée dès que $n_0 > 1$. En revanche, la structure $\xi'_0$ obtenue en collant $\xi$ sur $\xi_0$ est universellement tendue. □

# D   Les fibrés en tores au-dessus du cercle

*Démonstration du théorème* 1.3. On note $\pi$ la projection canonique $\mathbf{T}^2 \times \mathbf{R} \to \mathbf{T}^3_A$ et $F_t$, $t \in \mathbf{R}$, le tore $\mathbf{T}^2 \times \{t\}$. Pour toute structure de contact $\xi$ sur $\mathbf{T}^3_A$, on pose $\widetilde{\xi} = \pi^*\xi$. L'idée de la démonstration qui suit est simplement que, si $\xi$ est sous forme normale et tendue, alors $\widetilde{\xi}$ est sous forme normale et sa restriction à chaque produit $\mathbf{T}^2 \times [t, t+1]$, $t \in \mathbf{R}$, est tendue. D'après les propositions 3.15 et 3.22, toute structure de contact tendue sur $\mathbf{T}^3_A$ est isotope à une structure sous forme normale.

*Démonstration du* a). Soit $\xi$ une structure de contact universellement tendue et sous forme normale sur $\mathbf{T}^3_A$. On distingue deux cas.

Si $\mathfrak{F}_\perp(\xi)$ est vide, $\xi$ est une structure rotative et $\widetilde{\xi}$ a une amplitude non nulle. Le corollaire 3.8 montre alors que $\xi$ est isotope à l'une des structures $\zeta_n$ du théorème 3.7. D'autre part, la proposition 3.42 montre que la $2\pi$-torsion de $\zeta_n$ vaut $n$.

Si $\mathfrak{F}_\perp(\xi)$ est non vide, tous les tores $F_t$ sont convexes. Les feuilletages $\widetilde{\xi} F_t$ sont donc scindés par des courbes essentielles et $A$ préserve la classe de ces courbes, au signe près. Par suite, $A$ est conjuguée à une matrice du type $\varepsilon \left( \begin{smallmatrix} 1 & 0 \\ k & 1 \end{smallmatrix} \right)$, où $\varepsilon = \pm 1$. Si $\varepsilon$ vaut 1 (resp. $-1$), il existe une fibration de $\mathbf{T}^3_A$ au-dessus du tore $\mathbf{T}^2$ (resp. de la bouteille de Klein $\mathbf{K}^2$) telle que $\mathfrak{F}_\perp(\xi)$ soit isotope à l'image inverse d'un système de courbes fermées simples. Si ces courbes sont en nombre $2m$, il résulte du lemme 2.6 que $\xi$ est conjuguée à la structure de contact $\eta_m$ d'équation

$$\sin(2m\pi x)\,(dy - kt\,dx) + |k+1|\cos(2m\pi x)\,dt = 0, \qquad (x, y, t) \in \mathbf{T}^2 \times \mathbf{R}.$$

Si $\varepsilon = 1$, cette structure est elle-même conjuguée à une des structures $\zeta_n$ mais ceci est faux si $\varepsilon = -1$ (voir [Gi4] où ces structures sont classifiées à conjugaison et isotopie près).

Ainsi, si $\operatorname{tr}(A) \neq -2$, les structures de contact universellement tendues et de $2\pi$-torsion $n$ sur $\mathbf{T}^3_A$ forment une unique classe de conjugaison représentée par $\zeta_n$. Lorsque $\operatorname{tr}(A) \neq \pm 2$, cette classe de conjugaison est une classe d'isotopie. En effet, tout difféomorphisme de $\mathbf{T}^3_A$ est isotope à un difféomorphisme fibré, lequel envoie $\zeta_n$ sur une structure isotope à $\zeta_n$ [Gi4]. □

*Démonstration du* b). Soit $\xi$ une structure de contact virtuellement vrillée, orientée et sous forme normale sur $\mathbf{T}^3_A$. D'après la proposition 3.42, la torsion de $\xi$ est nulle. Pour tout $t \in \mathbf{R}$, on note $D_t \subset H_1(\mathbf{T}^2; \mathbf{R})$ la direction asymptotique du feuilletage $\widetilde{\xi} F_t$.

**Lemme 4.8.** *Dans* $\mathbf{P}^1 = \mathbf{P}\bigl(H_1(\mathbf{T}^2; \mathbf{R})\bigr)$, *l'image* $\Delta$ *de l'application* $t \in \mathbf{R} \mapsto D_t$ *est* :
  - $\mathbf{P}^1$ *si $A$ est d'ordre fini ou conjuguée à une matrice du type* $\left( \begin{smallmatrix} 1 & 0 \\ k & 1 \end{smallmatrix} \right)$ *avec* $k \geq 0$ ;
  - $\mathbf{P}^1$ *privé du point fixe de $A$ si $A$ est conjuguée à une matrice du type* $\left( \begin{smallmatrix} 1 & 0 \\ k & 1 \end{smallmatrix} \right)$ *avec* $k < 0$ ;
  - *l'intervalle de* $\mathbf{P}^1$ *limité à gauche (resp. à droite) par la direction stable (resp. instable) de $A$ si $A$ est hyperbolique.*

*En particulier, $\Delta$ ne dépend pas de $\xi$.*

*Démonstration.* L'ensemble $\Delta$ est un intervalle de $\mathbf{P}^1$ invariant par $A$ qui contient l'arc décroissant $D_t$, $t \in [0, 1]$. Si $A$ est d'ordre fini, ces conditions impliquent immédiatement que $\Delta = \mathbf{P}^1$. Si $A$ est conjuguée à une matrice du type $\left( \begin{smallmatrix} 1 & 0 \\ k & 1 \end{smallmatrix} \right)$ avec $k > 0$, elles entraînent que $\Delta$ contient la direction propre de $A$ et est donc aussi $\mathbf{P}^1$ tout entier. Pour terminer, il suffit donc de montrer que, si $A$ est hyperbolique ou conjuguée à une matrice du type $\left( \begin{smallmatrix} 1 & 0 \\ k & 1 \end{smallmatrix} \right)$ avec $k < 0$, l'ensemble $\Delta$ ne contient aucune droite propre. Or, si $D_{t_0}$ est une droite propre de $A$, il existe un point $t_1 \in \mathbf{R}$ tel que les droites $d_t$, $t \in [t_1, t_1 + 1]$, fassent plus d'un tour sur $\mathbf{P}^1$. La restriction de $\xi$ à $\mathbf{T}^2 \times [t_1, t_1 + 1]$ n'est donc pas tendue (proposition 3.36) et $\xi$ non plus. □

Soit maintenant $D \in \Delta$ une droite rationnelle et $C$ l'un des deux cônes convexes fermés bordés à gauche par $D$ et à droite par $A(D)$. Comme $D$ est rationnelle, l'enveloppe convexe $E$ des points entiers de $C$ ne contient sur son bord qu'un nombre fini de points primitifs. D'autre part, Pour toute structure de contact $\xi$ virtuellement tendue et sous forme normale sur $\mathbf{T}_A^3$, il existe un point $t \in \mathbf{R}$ tel que le feuilletage $\widetilde{\xi}F_t$ soit scindé par deux courbes et ait pour direction asymptotique la droite $D$. La proposition 3.36 n'offre alors qu'un nombre fini de lieux de retournement possibles pour $\widetilde{\xi}$ (ou pour la restriction de $\widetilde{\xi}$ à $\mathbf{T}^2 \times [t, t+1]$). Or le corollaire 3.28 montre que ce lieu détermine $\xi$ à isotopie près. Les classes d'isotopie de structures de contact virtuellement vrillées sont donc en nombre fini. □

*Calcul de $N(A)$ pour $A$ hyperbolique.* Soit $\widetilde{C}$ une composante connexe de l'ouvert

$$\left\{(x, y) \in \mathbf{R}^2 \mid \operatorname{tr}(A) A(x, y) \wedge (x, y) > 0\right\}$$

et $\widetilde{E}$ l'enveloppe convexe des points entiers de $\widetilde{C}$. On choisit une droite $D$ passant par un sommet de $\partial \widetilde{E}$ et on note $C$ le cône fermé dans $\widetilde{C}$ bordé par $D$ et $A(D)$.

**Lemme 4.9.** *Une structure de contact $\xi$ sous forme normale sur $\mathbf{T}_A^3$ est virtuellement vrillée si et seulement si les conditions suivantes sont satisfaites :*
- *$R(\widetilde{\xi})$ est non vide et inclus dans $\partial \widetilde{E}$ ;*
- *la direction asymptotique $D_t$ de chaque feuilletage $\widetilde{\xi}F_t$, $t \in \mathbf{R}$, rencontre le cône $\widetilde{C}$.*

*Démonstration.* Si $R(\widetilde{\xi})$ est inclus dans $\partial \widetilde{E}$ et si chaque droite $D_t$ rencontre $\widetilde{C}$, la structure $\widetilde{\xi}$ est tendue d'après le lemme 4.3. Par suite, $\xi$ est tendue, et donc virtuellement vrillée lorsque $R(\widetilde{\xi})$ est non vide. Inversement, si $\xi$ est virtuellement vrillée, les droites $D_t$ balaient $\widetilde{C}$ (lemme 4.8) et $R(\widetilde{\xi})$ est non vide (proposition 3.29). On choisit alors un point $t \in \mathbf{R}$ pour lequel $D_t = D$ et $\widetilde{\xi}F_t$ est scindé par deux courbes. Comme l'enveloppe convexe des points entiers du cône $C$ n'est autre que $\widetilde{E} \cap C$, il résulte de la proposition 3.36 que $R(\widetilde{\xi})$ est inclus dans $\partial \widetilde{E}$. □

On termine maintenant le calcul de $N(A)$. Soit $\xi$ une structure de contact virtuellement vrillée, sous forme normale et orientée. Quitte à déformer $\xi$ par isotopie parmi les structures sous forme normale, on suppose que, pour un certain $t \in \mathbf{R}$, le feuilletage orienté $\widetilde{\xi}F_t$ est une suspension scindée par deux courbes dont les cycles asymptotiques sont sur $D \cap \widetilde{C}$. On pose alors

$$\widehat{R}(\xi) = R\bigl(\widetilde{\xi}\bigl|_{\mathbf{T}^2 \times [t, t+1]}\bigr).$$

Selon que la trace de $A$ est positive ou négative, $\widehat{R}(\xi)$ a un nombre pair ou impair de points. Par ailleurs, comme $\xi$ n'est pas universellement tendue, $\widehat{R}(\xi)$ n'est ni vide ni formé des deux points primitifs de $\partial C$. Il résulte dès lors du corollaire 1.10 que le nombre de valeurs prises par la fonction $\gamma$ sur les parties permises et de cardinal impair (resp. pair) est $\prod_{i=0}^{n}(a_{m+i} - 1)^k$ (resp. $\prod_{i=0}^{n}(a_{m+i} - 1)^k - 2$). Reste à voir que $\gamma\bigl(\widehat{R}(\xi)\bigr)$ est un invariant d'isotopie.

Soit $\xi_0$ et $\xi_1$ deux structures de contact virtuellement vrillées, sous forme normale et orientées. On suppose que, pour un certain $t_i \in \mathbf{R}$, le feuilletage orienté $\widetilde{\xi}_iF_{t_i}$ est une suspension scindée par deux courbes dont les cycles asymptotiques sont sur $D \cap \widetilde{C}$. Si $\gamma\big(\widehat{R}(\xi_0)\big) \neq \gamma\big(\widehat{R}(\xi_1)\big)$, la restriction de $\widetilde{\xi}_0$ à $\mathbf{T}^2 \times [t_0, t_0+1]$ n'a aucun plongement dans $(\mathbf{T}^2 \times \mathbf{R}, \widetilde{\xi}_1)$ qui soit isotope à l'inclusion.

Quitte à passer au revêtement double $\mathbf{T}^3_{A^2} \to \mathbf{T}^3_A$, il suffit de prouver cette assertion lorsque la trace de $A$ est positive. S'il existe un tel plongement, son image $W$ est contenue dans un produit compact $\mathbf{T}^2 \times [t_1 - m, t_1 + m]$ et chaque composante du complémentaire de Int $W$ est un tore épais. En mettant sous forme normale la restriction de $\widetilde{\xi}_1$ à ces composantes puis en exprimant de deux manières la classe d'Euler relative de $\widetilde{\xi}_1 \mid_{\mathbf{T}^2 \times [t_1-m, t_1+m]}$, on obtient une contradiction avec le lemme 1.9-b). □

*Calcul de $N(A)$ pour $A$ d'ordre fini.* Comme $A$ est d'ordre fini, ou bien $A = \pm I$, ou bien $A$ est conjuguée à une matrice du type $J_{\varepsilon,k} = \varepsilon \left(\begin{smallmatrix} k & 1 \\ -1 & 0 \end{smallmatrix}\right)$ où $\varepsilon = \pm 1$ et $k \in \{-1, 0, 1\}$. Dans tous les cas, si $\xi$ est une structure de contact sous forme normale sur $\mathbf{T}^3_A$, la direction asymptotique $D_t$ de $\widetilde{\xi}F_t$ tourne indéfiniment sur $\mathbf{P}^1$ quand $t$ parcourt $\mathbf{R}$ (lemme 4.8).

Si $A = \pm I$ et si $R(\widetilde{\xi})$ est non vide, la preuve de la proposition 3.36 montre comment construire dans $(\mathbf{T}^3_A, \xi)$ un anneau à bord legendrien qui viole l'inégalité de Bennequin.

Si $A = J_{\varepsilon,k}$, soit $D$ le demi-axe des $x > 0$ et $C$ le cône convexe fermé dans $\mathbf{R}^2 \setminus \{0\}$ bordé à gauche par $D$ et à droite par $\varepsilon A(D)$. L'enveloppe convexe $E$ des points entiers de $C$ n'a sur son bord que deux points entiers primitifs, lesquels sont aussi sur $\partial C$. Pour toute structure de contact orientée, virtuellement vrillée et sous forme normale, on définit alors $\widehat{R}(\xi)$ comme dans le cas hyperbolique. Ici, le cardinal de $\widehat{R}(\xi)$ est pair ou impair selon que $\varepsilon$ vaut 1 ou $-1$. On voit ainsi immédiatement que $N(A)$ est nul si $\varepsilon = 1$ et vaut au plus 2 si $\varepsilon = -1$. Soit maintenant $D'$ le demi-axe des $y < 0$ et $C'$ le cône fermé dans $\mathbf{R}^2 \setminus \{0\}$ bordé à gauche par $D'$ et à droite par $\varepsilon A(D')$. Si $\varepsilon = -1$ et $k = -1$, le point primitif de $-A(D)$ est dans $C'$ mais n'est pas sur le bord de l'enveloppe convexe $E'$ des points entiers de $C'$. Un argument identique exclut le point entier primitif de $D$. □

*Calcul de $N(A)$ pour $A$ parabolique.* Soit $A$ la matrice $\left(\begin{smallmatrix} 1 & 0 \\ k & 1 \end{smallmatrix}\right)$, $k \neq 0$, $D$ le demi-axe des $x$ positifs, $C$ le cône bordé à gauche par $D$ et à droite par $\pm A(D)$ et $E$ l'enveloppe convexe des points entiers de $C$. Les points entiers primitifs de $\partial E$ sont tous alignés et il y en a $|k| + 1$ si $k$ est négatif, 3 si $k$ est positif et pair, 2 sinon. Le nombre de valeurs prises par la fonction $\gamma$ sur les sous-ensembles permis et de cardinal pair est donc, selon les cas, $|k| - 1$, 1 et 0. □

# Références

Emmanuel GIROUX  
*Unité de Mathématiques Pures et Appliquées,*  
*École Normale Supérieure de Lyon,*  
*46, allée d'Italie,*  
*69364, Lyon cedex 07, France*